\documentclass{amsart}
\usepackage{amssymb}
\usepackage{amscd}
\usepackage{amsmath}
\usepackage{enumerate}
\usepackage[dvips]{graphicx}
\newtheorem{theo}{Theorem}[section]
\newtheorem{lema}[theo]{Lemma}
\newtheorem{cor}[theo]{Corollary}
\newtheorem{prop}[theo]{Proposition}
\newtheorem{definition}[theo]{Definition}
\newtheorem{adendum}[theo]{Adendum}
\newtheorem{remark}[theo]{Remark}
\newtheorem{notation}[theo]{Notation}
\newtheorem{example}[theo]{Example}

\newcommand{\CC}{{\mathbb{C}}}
\newcommand{\NN}{{\mathbb{N}}}
\newcommand{\RR}{{\mathbb{R}}}
\newcommand{\ZZ}{{\mathbb{Z}}}
\newcommand{\calA}{{\mathcal{A}}}
\newcommand{\calD}{{\mathcal{D}}}

\newcommand{\calF}{{\mathcal{F}}}
\newcommand{\calG}{{\mathcal{G}}}
\newcommand{\calH}{{\mathcal{H}}}

\newcommand{\calK}{{\mathcal{K}}}
\newcommand{\calM}{{\mathcal{M}}}

\newcommand{\calO}{{\mathcal{O}}}
\newcommand{\calP}{{\mathcal{P}}}

\newcommand{\calT}{{\mathcal{T}}}
\newcommand{\calX}{{\mathcal{X}}}
\newcommand{\calY}{{\mathcal{Y}}}
\newcommand{\calZ}{{\mathcal{Z}}}
\newcommand{\comp}{{\circ}}
\newcommand{\ocn}{{\calO_{\CC^n,O}}}
\newcommand{\mm}{{\mathbf{m}}}
\newcommand{\orb}{{\mathrm{Orb}}}
\newcommand{\trans}{{\pitchfork}}
\newcommand{\codim}{{\mathrm{codim}}}

\newcommand{\rk}{{\mathrm{rk}}}
\newcommand{\Id}{{\mathrm{Id}}}
\newcommand{\Sing}{{\mathrm{Sing}}}
\begin{document}
\title[Relative Morsification Theory]{Relative Morsification Theory}
\author{Javier Fern\'andez de Bobadilla}
\address{Mathematisch Instituut. Universiteit Utrecht. Postbus 80010. 3508 TA Utrecht. The Netherlands.}
\email{bobadilla@math.uu.nl}
\thanks{Supported by the Netherlands Organisation for Scientific Research (NWO). Supported by the Spanish MCyT project BFM2001-1448-C02-01}
\date{}
\subjclass[2000]{Primary 32S30; 32S55}
\begin{abstract}
In this paper we develope a Morsification Theory for holomorphic functions 
defining a singularity of finite codimension with respect to an ideal, which recovers most previously known Morsification results for
non-isolated singularities and generalise them to a much wider context. We also show that deforming functions of finite codimension with 
respect to an ideal within the same ideal respects the Milnor fibration. Furthermore we present some applications of the theory: we introduce new numerical 
invariants for non-isolated singularities, which explain various aspects of the deformation of functions within an ideal; we define generalisations of the 
bifurcation variety in the versal unfolding of isolated singularities; applications of the theory to the topological study of the Milnor fibration 
of non-isolated singularities are presented. Using intersection theory in a generalised jet-space we show how to interpret the newly defined invariants as 
certain intersection multiplicities; finally, we characterise which invariants can be interpreted as intersection multiplicities in the above mentioned 
generalised jet space.
\end{abstract}
\maketitle

\setcounter{section}{-1}

\section{Introduction} 
Let $\ocn$ be the ring of germs of holomorphic functions at the origin $O$ of $\CC^n$. Two germs $f$ and $g$ are 
R-equivalent (right-equivalent) if there exists a germ of biholomorphism $\phi:(\CC^n,O)\to (\CC^n,O)$ such that 
$f\comp\phi=g$. One of the main aims of singularity theory is the classification of germs of holomorphic functions up to
R-equivalence; this includes giving normal forms for each of the equivalence classes and invariants to decide whether two
functions belong to the same class, but also studying the adjacencies (hierarchy) between equivalence classes (a class $C$ 
is adjacent to a class $C'$ if every function of $C'$ can be deformed into a function of $C$ by arbitrarily small 
deformation). When Arnold worked out the beginning of the classification of isolated singularities he observed that 
infinite series of classes of singularities occur (like $A_k$ or $D_k$), it appeared clear that the series where associated
with non-isolated singularities, and that the hierarchy of the series reflects the hierarchy of non-isolated singularities;
therefore the study of classification and hierarchy of non-isolated singularities is interesting not just by itself, but 
also in connection with the study of series of isolated singularities.

Many useful invariants for the classification and adjacency problems stem from the study of the Milnor fibration 
associated to a holomorphic function (vanishing (co)-homology, monodromy, intersection form, spectrum, homotopy type of the 
Milnor fibre, etc). This is quite well understood in case that the function has an isolated singularity, but when 
non-isolated singularities are present our knowledge is still very limited. See~\cite{Si4} for a recent survey of known results and open questions.

A fruitful way of studying isolated singularities is the so called morsification method:
Denote by $B_\epsilon$ the closed ball of radius $\epsilon$ centered at the 
origin of $\CC^n$; denote by $D_\eta$ the closed disk of radius $\eta$ centered at 
$0$ in $\CC$. Given a holomorphic germ $f\in\ocn$ with an isolated singularity, let $\epsilon$, $\eta$ be a pair of radii 
for which the Milnor fibration of $f$ is defined; then, any (small enough) generic perturbation $g$, has, as critical locus
inside $B_\epsilon$, as many Morse-type singularities as its Milnor number; moreover the Milnor fibration of 
$f$ is preserved by the perturbation in the following sense: the restriction 
\[g:B_\epsilon\cap g^{-1}(\partial D_\eta)\to\partial D_\eta\] 
is $C^\infty$-equivalent to the Milnor fibration of $f$. This gives a powerful method to study the Milnor fibration: 
for example it allows to compute the homotopy type of the Milnor fibre, and to relate the monodromy with the intersection 
form via the Picard-Lefschetz theory. 

The goal of this paper is to generalise this method to a wide class of non-isolated 
singularities. However, an arbitrary small perturbation of a function $f$ with non-isolated singularities does not preserve 
the Milnor fibration in the sense explained above; therefore we need to restrict the type of deformations that we will 
allow. The main idea is allow only deformations of $f$ within a suitable ideal $I$ of $\calO_{\CC^n,O}$ (which in some cases is the ideal of
functions that are singular where $f$ is singular, with the same generic transversal type). This point of view has been 
used in~\cite{Si1},~\cite{Si2},~\cite{Pe1},~\cite{Pe2},~\cite{Pe3},~\cite{Jo},~\cite{Za},~\cite{Ne}, to prove generalised 
Morsification Theorems and study the Milnor fibration of functions with smooth $1$-dimensional critical locus and simple 
transversal type, or with an i.c.i.s. (isolated complete intersection singularity) of dimension at most $2$ as critical locus and transversal type $A_1$.
Recently the Morsification Theorem has been generalised in~\cite{Bo1} 
for critical locus an i.c.i.s. of any dimension with transversal type $A_1$. In all these works the Morsification Theorem 
and the preservation of the Milnor fibration are proved exploiting special properties of the 
considered ideal (low dimensionality of its zero set or defining a complete intersection), which do 
not generalise easily. In all the cases the functions that can be morsified are precisely the functions of finite 
extended codimension with respect to the considered ideal $I$ in the sense of~\cite{Pe2} (see~Section~\ref{sect1} for a 
definition). Using a new conceptual approach, in this paper we generalise the morsification and preservation of Milnor 
fibration for functions of finite extended codimension with respect to an ideal, for any ideal of $\calO_{\CC^n,O}$. Moreover we introduce
and study new numerical invariants for non-isolated singularities. Below we summarise informally the main results of 
the paper.

In Section~\ref{sect2} we establish the preservation of the Milnor fibration: given {\em any} ideal 
$I$ of $\calO_{\CC^n,O}$ and a function $f$ of finite extended codimension with respect to $I$, we prove that any small 
enough perturbation of $f$ within the ideal preserves the Milnor fibration (see Theorem~\ref{topconst}).

The main and most difficult result of this paper is a Relative Morsification Theorem (see Theorem~\ref{morsificacion})  
that roughly 
states the following: given any function $f$ of finite extended codimension with respect to an ideal $I$, 
a generic small perturbation $g$ of $f$ within the ideal only has singularities that are 
simplest with respect to
$I$ in a certain natural sense. Moreover, in a small neighbourhood of the origin there are finitely many points in which
$g$ has positive extended codimension, with the property that the germ of $g$ at any of them is unsplittable, in the 
sense that it can not be splitted into several points of positive codimension by further perturbation. Furthermore, the 
points of positive extended codimension that are outside $V(I)$ are always Morse points of $g$. The proof uses the theory developed in Sections $1$ and
$3-8$. Let us sketch it here. Fix an ideal $I$.

Section~1 is preliminary. Some known definitions and results are recalled. Among them it is the concept of extended codimension of a function of $I$ with 
respect to $I$. This is an invariant that generalises the Milnor number to our setting. Indeed, the functions of finite extended codimension 
with respect to $I$ will be the ones that can be studied via the Relative Morsification Theorem (as isolated singularities are the functions whose 
properties can be studied via the usual morsification theory). In particular, by a theorem of Pellikaan (see Theorem~\ref{unfoldingtheo}), any function of 
finite extended codimension has an unfolding which is versal within $I$. On the other hand the extended codimension does not have all the good properties that
the Milnor number has. The Milnor number is conservative in the following 
sense: let $f$ be a function with finite Milnor number at $O$; there is a small neighbourhood $U$ of the origin such that, given any deformation $f+tg$, if $t$
is small enough, the sum of the Milnor numbers of $f+tg$ at the points of $U$ equals the Milnor number of $f$ at the origin. We show with an example that 
the extended codimension with respect to an ideal is not conservative in general. 

Consider a coherent ideal sheaf $\tilde{I}$ defined in a neighbourhood $U$ of the origin such that $\tilde{I}_O=I$, define 
\begin{equation}
\label{intro1}
J^\infty(U,\tilde{I}):=\coprod_{x\in U}\tilde{I}_x.
\end{equation}
In Section~3 we give a infinite-dimensional analytic structure to $J^\infty(U,\tilde{I})$, viewing it as a generalised $\infty$-jet space associated to 
$\tilde{I}$. For this we view $J^\infty(U,\tilde{I})$ as a the projective limit of 
\[J^m(U,\tilde{I}):=\coprod_{x\in U}\tilde{I}_x/(\mm^{m+1}\cap\tilde{I}_x),\]
when $m\geq 0$, and give a certain finite-dimensional generalised analytic structure to each $J^m(U,\tilde{I})$. In particular
we endow $J^\infty(U,\tilde{I})$ with a topology and define concepts as finite-determined set (the finite 
determinacy of a set implies that it can be understood as a subset of a finite-dimensional analytic space), analytic set, irreducibility, smoothness and 
codimension in $J^\infty(U,\tilde{I})$ of a closed analytic subset. We
prove the existence of decompositions of analytic subsets in irreducible components.

In Section~4 we consider the filtration of $J^\infty(U,\tilde{I})$ by sets containing germs of ascending extended codimension. We prove that the levels of the 
filtration are finite-determined closed analytic subsets of $J^\infty(U,\tilde{I})$. Despite this is more or less clear from the intuitive point of view, the 
proof gets rather technical. We also give a lower bound for the codimension in $J^\infty(U,\tilde{I})$ of the set of germs of 
positive and finite extended codimension.

In Section~5 we stratify the space $J^\infty(U,\tilde{I})$ in locally closed analytic strata containing germs of the same topological type. This is done 
combining the results of~\cite{Bo2} with the theory developed in Section~3.

In Section~6 the concept of Whitney regularity is defined for stratifications of $J^\infty(U,\tilde{I})$. Given any locally finite partition of 
$J^\infty(U,\tilde{I})$ by locally closed analytic subsets, it is proved the existence of a canonical Whitney 
stratification refining it. 

In Section~7 we define analytic mappings from an analytic set to $J^\infty(U,\tilde{I})$. After, the concept of transversality of a mapping to a submanifold
of $J^\infty(U,\tilde{I})$ is introduced. The Parametric Transversality Theorem is extended to our setting. It is shown that versality implies transversality
in the following sense: let $f\in I$ be a function of finite extended codimension and $F:(\CC^n,O)\times(\CC^r,O)\to\CC$ be a versal unfolding of $f$ within 
$I$, then mapping 
\[\rho_F:(\CC^n,O)\times(\CC^r,O)\to J^\infty(U,\tilde{I})\]
 assigning to $(x,s)$ the germ of $F_{|\CC^n\times\{s\}}$ at $x$ is analytic and transversal to any submanifold of $J^\infty(U,\tilde{I})$ which is invariant
by the natural action of holomorphic diffeomorphims preserving $\tilde{I}$.

Finally, in Section~8 the Relative Morsification Theorem (see Theorem~\ref{morsificacion}) is stated. A certain Whitney stratification of 
$J^\infty(U,\tilde{I})$ (which is built up
from the filtration by ascending extended codimension and the stratification considered in Section~5 using the geometric tools developed in Sections~3~and~6)
is needed. The proof is a transversality argument using the results developed in Section~7. As a by-product of our method we introduce generalizations of the 
bifurcation variety in the base of a versal unfolding.

The rest of the paper presents applications of our theory and several examples:

In Section~9 we present two immediate applications of the Morsification Theorem. In the first application we define new numerical invariants for functions 
of finite codimension with respect to an ideal, namely the {\em splitting function}, the {\em corrected extended codimension} and the {\em Morse number}. All 
of them are conservative in the sense explained above. The corrected extended codimension is close to the extended codimension, and could be thought as a 
conservative version of it. The splitting function is a finer invariant. A consequence of the Relative  Morsification Theorem is that the only singularities 
outside $V(I)$ of a generic deformation within $I$ of a function of finite extended codimension form a finite set of Morse points; the Morse number is 
the cardinality of this set. The second application, based also on the results of Section~2, shows how to study the 
topology of a function of finite extended codimension using a morsification. In particular, homology splitting and bouquet 
decomposition theorems are presented (see Theorem~\ref{bouquet}). 

In Section~10 we study further numerical invariants.
Given any $f\in I$ we have a jet-extension mapping $\rho_f:(\CC^n,O)\to J^\infty(U,\tilde{I})$. We develope a bit of intersection theory on the generalized 
jet space $J^\infty(U,\tilde{I})$ so that the intersection number at the origin of the mapping $\rho_f$ with any $n$-codimensional subvariety of 
$J^\infty(U,\tilde{I})$ is well defined and have some natural properties. This will enable to interpret the numerical invariants introduced above as 
intersection numbers, and give a formula of them in terms of dimensions of certain complex vector spaces (this becomes more explicit in the case of the Morse 
number). Using the developed intersection theory we are able to prove that any numerical invariant which is conservative and satisfies two other natural 
properties is actually a linear combination of intersection multiplicities with $n$-codimensional subvarieties of $J^\infty(U,\tilde{I})$ which are invariant
by the natural action of holomorphic diffeomorphims preserving $\tilde{I}$ (see Theorem~\ref{invchar}). A classification of such subvarieties yields in many 
examples on Section~11 an 
expression of any conservative invariant as a linear combination of well known invariants (for example in terms of the Morse number and number of 
$D_\infty$-points in the case of isolated-line singularities; see~\cite{Si1}). 

In the Section~11 we analize several applications and examples. In particular we point out previous morsification theorems that are now consequences of this 
theory. Among them is the Morsification Theorem for singularities whose critical locus is an $i.c.i.s.$ and have generic transversal type $A_1$; we show how
to easily extend this theorem to the case in which the critical locus is not necessarily a complete intersection (see Proposition~\ref{tipomorse}). 

\section{Functions of finite codimension with respect to an ideal} 
\label{sect1}
Let $\mm$ be the maximal ideal of $\ocn$, 
let $x_1,...,x_n$ be coordinate functions for $\CC^n_O$. The module of germs of vector fields at the origin $O$ of $\CC^n$ will
be denoted by $\Theta$; then $\mm\Theta$ are the vector fields vanishing at the origin. 
Denote by $\calD$ the group of germs of holomorphic diffeomorphims of $\CC^n$ fixing the origin, and by $\calD_e$ the set of germs of holomorphic 
diffeomorphims at the origin that not 
necessarily fix it, (we have no group structure because composition need not be defined). Let $I\subset \ocn$ be an ideal; given $U$, a small enough 
neighbourhood of the origin, there is a coherent sheaf of ideals $\tilde{I}$ whose stalk at $O$ is $I$. 
Following~\cite{Pe2} we define $\calD_I$ and $\calD_{I,e}$ to be respectively subgroup of 
$\calD$ and the subset of $\calD_e$ of elements preserving the ideal:
\begin{definition}
Define $\calD_{I,e}$ as the set of all $\phi\in\calD_e$ that have a representative $\phi:V\to W$, with $V$ and $W$ open subsets in $U$ and $O\in V$, such that
\[\phi^*(\Gamma(W,\tilde{I}))=\Gamma(V,\tilde{I}).\]
Define $\calD_I=\calD_{I,e}\cap\calD$.
\end{definition}    
Clearly, the action 
of $\calD$ on $\ocn$ by composition on the right restricts to an action $\sigma_I:I\times\calD_I\to I$. 
Given $f\in I$ we denote by $\orb(f)$ its orbit by $\sigma_I$.

Let $\phi_t$ be a 1-parameter family of holomorphic diffeomorphims of $\calD_e$ smoothly depending on $t$, such 
that $\phi_0=\Id_{\CC^n}$; let $\phi_{1,t},...,\phi_{n,t}$ be its components; consider $f\in\ocn$. The chain rule gives:
\begin{equation}
\label{reglacadena}
\frac{df\comp\phi_t}{dt}_{|t=0}=\sum_{i=1}^n\frac{\partial f}{\partial x_i}\frac{d\phi_{i,t}}{dt}_{|t=0}=X(f),
\end{equation}
where $X$ is the holomorphic vector field given by $X=\sum_{i=1}^nd\phi_{i,t}/dt_{|t=0}\partial/\partial x_i$. 
If $\phi_t\in\calD$ for any $t$ then, $X\in\mm\cap\Theta$. If $\phi_t\in\calD_{I,e}$ for any $t$ then $X(I)\subset I$.
Define $\Theta_{I,e}$ by the formula $\Theta_{I,e}:=\{X\in\Theta:X(I)\subset I\}$. 
If $\phi_t\in\calD_{I}$ for any $t$ then $X\in\mm\cap\Theta_{I,e}$; define 
$\Theta_I:=\mm\cap\Theta_{I,e}$. 

Conversely, integration associates to any $X\in\Theta$, a 1-parameter flow $\phi_t$ of 
holomorphic diffeomorphims of $\calD_e$, with $\phi_0=\Id_{\CC^n}$, such that if $X\in\Theta_{I,e}$ then $\phi_t\in\calD_{I,e}$, and if $X\in\Theta_I$
then $\phi_t\in\calD_I$.

Given a family $\phi_t\subset \calD_I$ as above, and $f\in I$, we may regard the family of functions $f\comp\phi_t$ as 
a smooth path in $\orb(f)$. This motivates the following:

\begin{definition}
\label{codimension}
Consider $f\in I$, define the {\em tangent space} and {\em extended tangent space} at $f$ to its orbit respectively by 
\[\tau_I(f):=\Theta_I(f)\quad\quad\tau_{I,e}:=\Theta_{I,e}(f)\]
Moreover we define the (possibly infinite) {\em $I$-codimension} and {\em extended $I$-codimension} respectively by
\[c_I(f):=\dim_\CC\frac{\int I}{\tau_I(f)}\quad\quad c_{I,e}(f):=\dim_\CC\frac{\int I}{\tau_{I,e}(f)}\]
\end{definition}

Notice that both $\tau_I(f)$ and $\tau_{I,e}(f)$ are ideals of $\ocn$. It is easy to see (cf.~\cite{Pe2}) that the 
$I$-codimension is finite if and only if the extended $I$-codimension is so. 

Our setting is more general that the one studied by Pellikaan in~\cite{Pe1},~\cite{Pe2}. There the ideal $I$ is asked to 
be the primitive of another ideal $I'\subset\ocn$: consider an ideal $I'\subset\ocn$, we define $\int I'$, the 
{\em primitive ideal of} $I'$, as 
\[\int I':=\{f\in\ocn: (f)+J_f\subset I'\},\]
where $J_f:=\{\partial f/\partial x_1,...,\partial f/\partial x_n\}$ is the 
Jacobian ideal of $f$. Heuristically $\int I'$ is the ideal of functions that vanish and ``are singular at'' the analytic 
space defined by $I'$. For many applications it is sufficient to work with the class of primitive ideals: for example in the study 
of non-isolated singularities with generic transversal type $A_1$. However the ideals considered by de Jong in~\cite{Jo} 
for the cases of line singularities with transversal types $A_3$, $E_7$, and $E_6$ (this last case for ambient space of
dimension at least $4$) are not primitive of any other ideal. 

Observe that when $I=\int I'$ two different definitions of (extended) $I$-codimension are possible: 
we can consider diffeomorphims
that preserve $I'$ (as Pellikaan does) instead of diffeomorphisms that preserve $\int I'$; hence we use
the modules $\Theta_{I',e}$ and $\Theta_{I'}$ instead of $\Theta_{\int I',e}$ and $\Theta_{\int I'}$. It is easy to see 
that $\Theta_{I',e}\subset \Theta_{\int I',e}$, but the equality is not known in general, up to the author's knowledge. This give rise to two, a priori, 
different morsification theories for the ideal $\int I$. The constructions and results of this paper are valid for any of them. 
Anyhow, the equality $\Theta_{I',e}=\Theta_{\int I',e}$ holds when $I'$ is a radical ideal (it is easy to prove that any
$X\in\Theta_{\int I',e}$ must be tangent to the analytic space $V(I')$ defined by $I'$ at each of its smooth points; this 
implies that $X\in\Theta_{I',e}$), and in all
the examples that we have considered.

\begin{notation}{\em
Let $\calF$ be a coherent sheaf on $U$. Denote by $\calF_x$ its stalk at $x$; given any section $\varphi\in\Gamma(\calF,U)$
denote by $\varphi_x$ its germ at $x$.}
\end{notation}

Let $\tilde{\Theta}$ be the free $\calO_U$-module of vector fields over $U$, and $\Theta_{\tilde{I},e}$
the coherent $\calO_U$-module of vector fields preserving $\tilde{I}$. The stalk of 
$\Theta_{\tilde{I},e}$ at any $x\in U$ is the $\calO_{U,x}$-module $\Theta_{\tilde{I}_x,e}$. Therefore, for any
$f$ holomorphic in $U$, we have that $\Theta_{\tilde{I},e}(f)$ is a coherent sheaf of ideals whose stalk at any 
$x\in U$ is $\tau_{I_x,e}(f_x)$. Define the coherent $\calO_U$-module
\begin{equation}
\label{Fcal}
\calF:=\frac{\tilde{I}}{\Theta_{\tilde{I},e}(f)}.
\end{equation}
Standard properties of coherent sheaves imply that $c_{I,e}(f)<\infty$ if and only if $\calF$ is concentrated at the 
origin, i.e. $\calF_x=0$ for $x\neq O$.

\begin{notation}{\em 
Let $\pi:\CC^n\times\CC^r\to\CC^r$ be the projection to the second factor; take any $s\in\CC^r$. 
If $\calG$ is a coherent analytic sheaf on $\CC^n\times\CC^r$ we denote by $\calG_{|s}$ the pullback of $\calG$ to 
$\pi^{-1}(s)$. 

If $F$ is an analytic function on $\CC^n\times\CC^r$ we denote by $F_{|s}$ the restriction of $F$ to 
$\pi^{-1}(s)$. Therefore $F_{|s,x}$ denotes the germ at $x$ of the restriction of $F$ to $\pi^{-1}(s)$. If $X$ is an 
analytic subset of $\CC^n\times\CC^r$ we denote by $X_{s}$ the fibre of $X$ over $s$.}
\end{notation}

\begin{definition}
\label{def:unfold}
Consider $f\in I$, a $r$-{\em parametric} $I$-{\em unfolding of} $f$ is an holomorphic germ 
$F:(\CC^n\times\CC^r,(O,O))\to\CC$ such that $F_{|O,O}=f$ and $F_{|s,O}$ belongs to $I$ for any $s\in\CC^r$. 
We denote by $I(r)$ the module formed by all the $r$-parametric $I$-unfoldings. Define 
$\Theta_{I,e}(r):=\pi^*\Theta_{I,e}=\calO_{\CC^n\times\CC^r}\Theta_{I,e}$. 
\end{definition}

\begin{lema}
\label{igualdadbasica} 
The following equalities hold: $I(r)=\pi^*I=I\calO_{\CC^n\times\CC^r,(O,O)}$.
\end{lema}
\begin{proof}
The only non-trivial statement is that $I(r)\subset\pi^*I$. Consider the coordinates $x_1,...,x_n$ for $\CC^n$ and fix
coordinates $s_1,...,s_r$ for $\CC^r$. Let $f_1,...,f_k$ be a set of generators for $I$. Take $F\in I(r)$; we have to find 
$G_1$,...,$G_k$, convergent power series in $x_1,...,x_r,s_1,...,s_r$, such that 
\begin{equation}
\label{primera}
F=\sum_{i=1}^kG_if_i
\end{equation}
By Artin's Aproximation Theorem it is enough to find formal power series $G_i$ satisfying the last equation.

Express each $G_i$ as $G_i=\sum g^i_{i_1,...,i_n}x_1^{i_1}...x_n^{i_n}$ and $F$ as 
$F=\sum A_{j_1,...,j_n}x_1^{j_1}...x_n^{j_n}$ where each $g^i_{j_1,...,j_n}$ and each $A_{j_1,...,j_n}$ is a power series
in $s_1,...,s_k$; express each $f_i$ as $f_i=\sum a^i_{j_1,...,j_n}x_1^{j_1}...x_n^{j_n}$ where each  
$a^i_{j_1,...,j_n}$ is a complex number. For any positive integer $N$ the truncation of Equality~(\ref{primera}) to its $N$-jet
with respect of $x_1,...,x_n$ may be seen as a linear system whose variables are $\{g^i_{j_1,...,j_n}:j_1+...+j_n\leq N\}$,
whose coefficients are $\{a^i_{j_1,...,j_n}:j_1+...+j_n\leq N\}$ (complex numbers), and whose independent terms are 
$\{A_{j_1,...,j_n}:j_1+...+j_n\leq N\}$ (holomorphic functions in $s_1,...,s_r$). The fact that $F$ is a $I$-unfolding 
implies
that for any value of $s_1,...,s_r$ close enough to the origin the system has a solution; using this and the fact that the 
rank of the fundamental matrix of the system does not depend on $s_1,...s_r$ (the $a^i_{j_1,...,j_n}$'s are complex 
numbers), we deduce that there exist a solution of the system depending holomorphically on $s_1,...s_r$ in a neighbourhood 
of the origin of $\CC^r$. This provides a solution for the truncation of Equality~\ref{primera} to its $N$-jet. Applying 
Krull's intersection theorem we deduce the existence of formal solutions.
\end{proof}

\begin{remark}{\em 
In the last lemma we have shown the following statement: consider $f_1,...,f_k,F\in\calO_{\CC^n\times\CC^r,(O,O)}$ such 
that $f_1,...,f_k$ are independent on $s_1,...,s_r$; if for any $s\in\CC^r$ close enough to the origin 
$F_{|s}\in (f_{1|s},...,f_{k|s})$ then $F\in (f_1,...,f_k)$. The independence of the $f_i$'s on $s_1,...,s_r$ is needed,
as the following example shows: $f(x,s)=xs^2$, $F(x,s)=xs$.}
\end{remark}

Given $U$ and $V$, neighbourhoods of the origin of $\CC^n$ and $\CC^r$ respectively, we define the coherent 
$\calO_{U\times V}$-modules $\tilde{I}(r)$ and $\Theta_{\tilde{I},e}(r)$ as $\tilde{I}(r):=\tilde{I}\calO_{\CC^n\times\CC^r}$ and 
$\Theta_{\tilde{I},e}\calO_{\CC^n\times\CC^r}$.  
 
\begin{lema}
\label{Gcalunfold}
Consider a representative $F:U\times V\to\CC$ of a $r$-parametric $I$-unfolding of a germ $f\in\int I$. Define a coherent 
$\calO_{U\times V}$-module by the formula
\begin{equation}
\label{Gcalcoh}
\calG:=\frac{\tilde{I}(r)}{(\Theta_{\tilde{I},e}(r))(F)}.
\end{equation}
Then, 
\begin{enumerate}
\item We have that $\calG_{(x,s)}=0$ if and only if $c_{\tilde{I}_x,e}(F_{|s,x})=0$. 
\item If $c_{I,e}(f)=0$ then $U$ and $V$ 
can be shrinked enough so that $c_{\tilde{I}_x,e}(F_{|s,x})=0$ for any $(x,s)\in U\times V$.
\item Let $Z\subset U\times V$ the support of $\calG$; let $p:U\times V\to V$ be the projection to the second factor; 
define $\varphi:=p_{|Z}$. If $c_{I,e}(f)<\infty$ then we can shrink $U$ and $V$ so that $\varphi$ is finite and 
$p_*\calG$ is a coherent $\calO_V$ module.
\end{enumerate} 
\end{lema}
\begin{proof}
Denote by $\mm_s$ the maximal ideal of $\calO_{\CC^r,s}$. A standard commutative algebra argument shows 
\begin{equation}
\label{formula1}
\frac{\calG_{(x,s)}}{\mm_s\calG_{(x,s)}}\cong\frac{\tilde{I}_x}{\tau_{\tilde{I}_x,e}(F_{|s,x})}.
\end{equation}
Therefore $\calG_{(x,s)}=0$ implies $c_{\tilde{I}_x,e}(F_{|s,x})=0$; conversely, if $c_{\tilde{I}_x,e}(F_{|s,x})=0$ then 
$\calG_{(x,s)}=\mm_s\calG_{(x,s)}$, and hence $\calG_{(x,s)}=0$ by Nakayama's Lemma. We have shown the first assertion

Suppose that $c_{I,e}(f)=0$; then, as 
$f=F_{|O,O}$ we have that $\calG_{(O,O)}=0$; as the support of a coherent $\calO_{U\times V}$-module is closed, 
the second assertion follows.

Suppose $c_{I,e}(f)<\infty$. 
Shrinking $U$ we can assume that $\calF$ is concentrated at the origin of $\CC^n$. Therefore, by the previous 
assertions, we
have that $\calG_{(x,O)}=0$ for any $x\in U\setminus\{O\}$. By the Projection Lemma of~\cite{GR},~page~62, we can shrink 
$U$ and $V$ so that the third assertion follows.
\end{proof}

A corollary of Lemma~\ref{Gcalunfold} is the upper semicontinuity of extended codimension:

\begin{cor}
\label{uppersem}
Let $f\in I$ with $c_{I,e}(f)<\infty$; consider $F:U\times V\to\CC$ a representative of a $I$-unfolding of $f$ such 
that the second statement of Lemma~\ref{Gcalunfold} holds. Let $\varphi$ be the mapping introduced in the last lemma. Then for $s$ close enough to the origin, 
\begin{equation}
\label{ineqcod}
\sum_{x\in \varphi^{-1}(s)}c_{\tilde{I}_x,e}(F_{s,x})\leq c_{I,e}(f)
\end{equation}
\end{cor}
\begin{proof}
After Lemma~\ref{Gcalunfold} we are reduced to prove that, for $s$ close enough to the origin,
$\dim_\CC(\varphi_*\calG)\otimes(\calO_{V,O}/\mm_O)\geq\dim_\CC(\varphi_*\calG\otimes(\calO_{V,s}/\mm_s)$,
which is standard for being $\varphi_*\calG$ coherent.
\end{proof}

When $I=\ocn$ then the extended codimension is equal to the Milnor number, in this case the inequality~\ref{ineqcod}
becomes an equality. Although this happens in many other cases such as isolated line singularities 
(see~\cite{Si1},~\cite{Pe3}) and most of the examples that we have computed, the following example shows that the equality 
does not hold in general.

\begin{example}
\label{ejem} 
{\em
Let $I=(x^2,y)\subset\calO_{\CC^2,O}$ (ideal of curves tangent to the $x$-axis); then $\Theta_{I,e}$ is generated by 
$x\partial/\partial x,y\partial/\partial x,x^2\partial/\partial y,y\partial/\partial y$. Consider $f(x,y)=y^2+x^3$; 
then $\tau_{I,e}(f)=(x^3,x^2y,y^2)$ and $c_{I,e}(f)=3$. Any unfolding $F:\CC^2\times\CC\to\CC$ of $f$ such that for any 
$s$ the curve $F_{|s}=0$ is tangent to the $x$-axis at the origin is a $I$-unfolding. Choose 
an $I$-unfolding such that $F_{|s}=0$ is smooth at the origin; then one checks easily that $c_{I,e}(F_{|s})=0$. Therefore 
the only terms that can contribute to the left hand side of the Inequality~\ref{ineqcod} are the extended codimensions
of $F_{|s}$ at points $x$ outside the support of $I$. As $I_x=\calO_{\CC^2,x}$ at these points, the extended codimension 
coincides with the Milnor number. As $\mu(f)=2$ the Inequality~\ref{ineqcod} is strict in this case.}
\end{example}
 
Let $F\in I(r)$, $G\in I(q)$ be $I$-unfoldings of $f\in I$. A morphism $\xi:F\to G$ of $I$-unfoldings is a
pair $(\Phi,\lambda)$ consisting of holomorphic germs $\lambda:(\CC^r,O)\to(\CC^q,O)$ and $\Phi:(\CC^n\times\CC^r,(O,O))\to \CC^n$ with the following 
properties: $\Phi_{|\CC^n\times\{s\}}\in\calD_{I,e}$ for any $s$, $\Phi_{|\CC^n\times\{O\}}=\Id_{\CC^n}$, and $G(\Phi(x,s),\lambda(s))=F$. 

\begin{definition}
Let $f\in I$ and $F$ an $I$-unfolding of it. We say that $F$ is versal if for any other $I$-unfolding $G$ of $f$ 
there exists a morphism of $I$-unfoldings from $G$ to $F$.
\end{definition}

The following theorem was proved by Pellikaan in~\cite{Pe2} as an application of the general results of Damon in~\cite{Da};
although Pellikaan proves it for primitive ideals, his arguments extend without change to our setting.

\begin{theo}[Unfolding Theorem]
\label{unfoldingtheo}
Let $F:\CC^n\times\CC^r\to\CC$ be an $I$-unfolding of a function $f\in I$. Let $s_1,...,s_r$ be the coordinates of the base
space $\CC^r$. The following statements are equivalent
\begin{enumerate}
\item $\tau_{I,e}(f)+\CC(\partial F/\partial s_{1})_{|s=0}+...+\CC(\partial F/\partial s_{r})_{|s=0}=I$.
\item $F$ is a versal $I$-unfolding of $f$.
\end{enumerate}
\end{theo}

A corollary of this is that $f\in I$ has a versal unfolding if and only if $c_{I,e}(f)<\infty$. Versality is open in 
the following sense:

\begin{prop}
\label{versalopen} 
Let $F$ be a versal $I$-unfolding of a function $f\in I$; let $F:U\times V\to\CC$ be a representative of the germ $F$.
Then $U$ and $V$ can be shrinked enough so that $F$ is a versal $\tilde{I}_x$-unfolding of $F_{|s,x}$ for any $(x,s)\in U\times V$.
\end{prop}
\begin{proof}
Let $\calF$ and $\calG$ be the sheaves associated to $f$ and $F$ by the formulae~(\ref{Fcal}) and~(\ref{Gcalcoh}) respectively.
As $f$ has a versal unfolding then $c_{I,e}(f)<\infty$. Hence, by the second assertion of Lemma~\ref{Gcalunfold}, we can 
shrink $U$ and $V$ so that $\varphi$ is finite and $p_*\calG=\varphi_*\calG$ is a coherent $\calO_V$-module.

Let $s_1,...,s_r$ the coordinates of $\CC^r$ and $(s)$ the ideal generated by them. The functions $\partial F/\partial s_1,...,\partial F/\partial s_r$, 
defined on the whole $U\times V$, can be seen as sections of $p_*\tilde{I}(r)$ over $V$; denote by $\partial_iF$ the image of 
$\partial F/\partial s_i$ by the natural homomorphism $p_*\tilde{I}(r)\to p_*\calG$. Define $\calM\subset p_*\calG$ to be the
coherent $\calO_V$-module generated by $\partial_1F,...,\partial_rF$. We claim that $\calM=p_*\calG$ if we shrink $V$ 
enough. To prove the claim we only need to show that $\calM_{O}=(p_*\calG)_{O}$. By Nakayama, this reduces to prove the 
equality $p_*\calG_{O}=\calM_{O}+(s)p_*\calG_O$. As $\varphi^{-1}(O)=\{(O,O)\}$, then $(p_*\calG)_O=\calG_{(O,O)}$; 
therefore $(p_*\calG)_{O}/(s)(p_*\calG)_{O}$ is equal to $\calG_{(O,O)}/(s)\calG_{(O,O)}$, which, by 
formula~\ref{formula1}, is isomorphic to $I/\tau_{I,e}(f)$. We have constructed an isomorphism
$\psi:(p_*\calG)_{O}/(s)(p_*\calG)_{O}\to I/\tau_{I,e}(f)$. It is easy to see that the image of $\partial_iF$ by $\psi$
is the class of $(\partial F/\partial s_{i})_{|O}$ in $I/\tau_{I,e}(f)$. As $F$ is a versal $I$-unfolding of $f$, 
by Theorem~\ref{unfoldingtheo} we conclude that the restriction of $\psi$ to 
$(\calM_{O}+(s)(p_*\calG)_{O})/(s)(p_*\calG)_{O}$ is surjective. This shows the claim.

Consider $(x,s)\in U\times V$, if $\calG_{(x,s)}=0$ then, by Lemma~\ref{Gcalunfold}, we have the equality 
$\tau_{\tilde{I}_x,e}(F_{|s,x})=\tilde{I}_x$;  
hence, by Theorem~\ref{unfoldingtheo} $F$ is a versal $\tilde{I}_x$-unfolding of $F_{|s,x}$. Suppose that $(x,s)\in\mathrm{Supp}(\calG)$. Let $\mm_s$ be the 
maximal ideal of $\calO_{\CC^r,s}$. By the 
finiteness of $\varphi$ we have an equality $(p_*\calG)_s=\bigoplus_{y\in\varphi^{-1}(s)}\calG_{(y,s)}$. Using 
formula~(\ref{formula1}) we obtain an isomorphism  
\[\psi_s:(p_*\calG)_s/\mathbf{m}_s(p_*\calG)_s=\bigoplus_{y\in\varphi^{-1}(s)}\calG_{(y,s)}/\mathbf{m}_s\calG_{(y,s)}\cong\bigoplus_{y\in\varphi^{-1}(s)} \tilde{I}_{y}/\tau_{\tilde{I}_y,e}(F_{|s,y}).\]
Using that $\calM=p_*\calG$, noting that the image of $\partial_i F$ by $\psi_s$ has $(\partial F/\partial x_i)_{|s,y}$ 
as component in $\tilde{I}_{y}/\tau_{\tilde{I}_y,e}(F_{|s,y})$, and taking into account that $x\in\varphi^{-1}(s)$ we obtain that 
$\tilde{I}_x=\tau_{\tilde{I}_x,e}(F_{s,x})+\CC(\partial F/\partial x_1)_{|s,x}+...+\CC(\partial F/\partial x_r)_{|s,x}$. Then, by 
Theorem~\ref{unfoldingtheo} the mapping $F$ is a versal unfolding of $F_{|s,x}$.
\end{proof}

\section{Topology of unfoldings of functions of finite $I$-codimension.}
\label{sect2}

Denote by $\dot{D}_\eta$ the punctured disk of radius $\eta$ centered at the origin and by $B_\epsilon$, $\overline{B}_\epsilon$ and $S_\epsilon$ the open 
ball, closed ball and sphere of radius $\epsilon$ centered at the origin of $\CC^n$;
let $f\in\ocn$. L\^e proved in~\cite{Le} that if $\epsilon>0$ is small enough and $\epsilon>>\eta>0$ then
\begin{equation}
\label{fibration}
f_{|\overline{B}_\epsilon\cap f^{-1}(\dot{D}_\eta)}:\overline{B}_\epsilon\cap f^{-1}(\dot{D}_\eta)\to\dot{D}_\eta
\end{equation}
is a locally trivial fibration, and, moreover, if $(\epsilon',\eta')$ is another pair with $\epsilon'\leq\epsilon$ and 
such that $f_{|\overline{B}_{\epsilon'}\cap f^{-1}(\dot{D}_{\eta'})}$ is also a 
locally trivial fibration, then both fibrations are equivalent. Moreover, due to Hironaka~\cite{Hi},~\S5, if we consider
in $\CC$ the stratification $\{\CC\setminus\{0\},\{0\}\}$, then there exist an analytic Whitney stratification of a 
neighbourhood $U$ of the origin of $\CC^n$ containing $\overline{B}_\epsilon$ such that $U\cap f^{-1}(\CC\setminus\{0\})$ is a 
stratum, the mapping $f:U\to\CC$ satisfy the Thom $A_f$-condition respect to this stratifications and, for each stratum 
$X\subset U$ and each point $x\in X\cap S_\epsilon$, we have $T_xX\trans S_\epsilon$. 
\begin{definition}
\label{gsr}
A pair $(\epsilon,\eta)$
with all the properties above is called a {\em good system of radii} for $f$.
\end{definition}
The fibration
\begin{equation}
\label{Milnorfib}
f_{|\overline{B}_\epsilon\cap f^{-1}(\partial D_\eta)}:\overline{B}_\epsilon\cap f^{-1}(\partial D_\eta)\to\partial D_\eta
\end{equation}
is called the {\em Milnor fibration} of $f$. 

The main result of this section is the preservation of the transversality with the Milnor sphere for unfoldings of functions of finite 
codimension:

\begin{theo}
\label{topconst}
Let $f\in I$ such that $c_{I,e}(f)<\infty$;
consider $F$ a $1$-parametric $I$-unfolding of $f$; let $(\epsilon,\eta)$ be a good system of radii for $f$. Given a value 
$s$ of the parameter consider the restriction
\begin{equation}
\label{qwer}
F_{|s}:X_{D_\eta,s}:=F_{|s}^{-1}(D_\eta)\cap\overline{B}_\epsilon\to D_\eta.
\end{equation}
Then $\delta$ can be chosen small enough so that, given $s,s'\in D_\delta$
\begin{enumerate}
\item If $t\in \dot{D_\eta}\setminus\{0\}$ then $F_{|s}^{-1}(t)\trans S_\epsilon$.
\item The locally trivial fibrations that $F_{|s}$ and $F_{|s'}$ induce over $\partial D_\eta$ are equivalent. 
\item $X_{D_\eta,s}$ and $X_{D_\eta,s'}$ (where $X_{D_\eta,s}:=F_{|s}^{-1}(D_\eta)\cap\overline{B}_\epsilon$) are homeomorphic.
\end{enumerate}
\end{theo}

This result is a generalisation of analogous statements in~\cite{Si1},~\cite{Si2},~\cite{Jo},~\cite{Za}, and~\cite{Bo1}; in those papers the result
is proved in the case in which the dimension of $V(I)$ is at most $2$, and/or $I$ is of a rather particular type, in all the 
cases the idea is to use the special properties of $I$ to control explicitly the geometry of $F$ at $S_\epsilon$. 
Next we show how to control the geometry of $F$ in a neighbourhood of $S_\epsilon$ for any $I$, without making use of any geometric property of $I$.

\begin{lema}
\label{controlext}
Let $F:\CC^n\times\CC\to\CC$ be a $1$-parametric $I$-unfolding, consider a smooth path $\gamma:(-\nu,\nu)\to\CC$ such
that $\gamma(0)=0$. There exist a positive number $\mu<\nu$, a neighbourhood $V_0$ of $S_\epsilon$ in $\CC^n$, a 
neighbourhood $W$ of $S_\epsilon\times (-\mu,\mu)$ in $\CC^n\times(-\mu,\mu)$, and a $C^\infty$-diffeomorphism
\begin{equation}
\label{controlfron}
\Psi:V_0\times (-\mu,\mu)\to W
\end{equation}
of the form $\Psi(x,t)=(\psi_t(x),t)$ (i.e. a $C^\infty$-family $\Psi$ of diffeomorphisms $\psi_t$ depending on $t$), such 
that $\psi_0=\Id_{\CC^n}$ and $F_{|\gamma(t)}\comp\psi_{t}=f_{|V_0}$ for any $t\in (-\mu,\mu)$.
\end{lema} 
\begin{proof}
As $c_{I,e}(f)<\infty$, the sheaf $\calF$ defined in Equation~\ref{Fcal} is concentrated at the origin of $\CC^n$. This 
implies that, if $\epsilon$ is small enough, for any $x\in S_\epsilon$, we have $\calF_x=0$ and hence 
$c_{\tilde{I}_x,e}(f_x)=0$. 
As $f_x=F_{|0,x}$, we have $c_{\tilde{I}_x,e}(F_{|0,x})=0$, which by Lemma~\ref{Gcalunfold}, implies 
$\calG_{(x,0)}=0$ for any $x\in S_\epsilon$. 
By coherence, there is an open neighbourhood $U$ of 
$S_\epsilon\times\{0\}$ in $\CC^n\times\CC$ such that $\calG_{|U}=0$. 

Let $s$ be a coordinate for the base space of the unfolding $F$; as the ideal sheaf
$\tilde{I}(1)$ is closed
under differentiation respect to $s$ we have $\partial F/\partial s\in\tilde{I}(1)$. As $\calG_{|U}=0$, 
for any $x\in U$ there exists an open neighbourhood $U^x$ of $x$ in $U$ and a vector 
field $X^x\in\Gamma(U^x,\Theta_{\tilde{I},e}(1))$
such that $X^x(F_{|U^x})=\partial F/\partial s_{|U^x}$. Choose $U^1,...,U^k$ a finite collection of the above open subsets 
such that $\cup_{i=1}^kU^i \supset S_\epsilon\times\{0\}$. Let $X^i$ be the vector field associated to $U^i$. 
Redefine $U$ to be the union of the chosen subsets. Consider a $C^\infty$ 
partition of unity $\{\rho^1,...,\rho^k\}$ subordinated to the cover $\{U^1,...,U^k\}$ and define 
$X:=\sum_{i=1}^k\rho^iX^i$. The vector field $X$ is of the form $X=\sum_{i=1}^n \alpha_i\partial/\partial x_i$ where 
$\alpha_i$ are $C^\infty$ functions over $U$, and satisfies 
$X(F_{|U})=\sum_{i=1}^k\rho^iX^i(F_{|U^i})=\partial F/\partial s_{|U}$. 

Choose $0<\mu<\nu$ and a neighbourhood $V_0$ of $S_\epsilon$ in $\CC^n$ such that $V_0\times\gamma(-\mu,\mu)$ is contained in $U$. 
Consider a $C^\infty$ vector field tangent to $V_0$ and smoothly depending on $t$ defined by $Y(x,t):=-(d\gamma/dt)X(x,\gamma(t))$. 
Express $Y$ as $Y=\sum_{i=1}^n\beta_i\partial/\partial x_i$, then $\beta_i(x,t)=-(d\gamma/dt)\alpha_i(x,\gamma(t))$. 
Integrating $Y$ (and may be shrinking $V_0$ and $\mu$) 
we obtain a $C^\infty$-family of diffeomorphisms $\psi:V_0\times (-\mu,\mu)\to U$ such that $\psi_0=\Id_{\CC^n}$; define $\Psi(x,t):=(\psi(x,t),t)$ and
$W:=\Psi(V_0\times (-\mu,\mu))$. We need to check that $F_{|\gamma(t)}\comp\psi_{t}=f_{|V_0}$ for any 
$t\in (-\mu,\mu)$. This is obvious for $t=0$. Hence if $G(x,t):=F(\psi(x,t),\gamma(t))$, it suffices to show that
$\partial G/\partial t=0$. By the chain rule,
\[\frac{\partial G}{\partial t}(x,t)=\sum_{i=1}^n\frac{\partial F}{\partial x_i}(\psi(x,t),\gamma(t))\beta_i(\psi(x,t),t)+
\frac{\partial F}{\partial s}(\psi(x,t),\gamma(t))\frac{d\gamma}{dt}(t)=\]
\[\sum_{i=1}^n\frac{\partial F}{\partial x_i}(\psi(x,t),\gamma(t))(-\frac{d\gamma}{dt}(t))\alpha_i(\psi(x,t),\gamma(t))+\frac{\partial F}{\partial s}(\psi(x,t),\gamma(t))\frac{d\gamma}{dt}(t)=\]
\[\frac{d\gamma}{dt}(t)[(-X(F)+\frac{\partial F}{\partial s})(\psi(x,t),\gamma(t))]=0\]
\end{proof}
\begin{proof}[Proof of Proposition 1]
The pair $(\epsilon,\eta)$ is chosen so that $\Sing(f)\cap\overline{B}_\epsilon\subset f^{-1}(0)$, therefore, $\delta$ can be
chosen small enough so that $F_s$ has no critical points at $F_{|s}^{-1}(\partial D_\eta)\cap\overline{B}_\epsilon$.
Define $\bar{F}:\CC^n\times\CC\to\CC\times\CC$ by $\bar{F}(x,s):=(F(x,s),s)$. The second assertion follows from the first 
applying Ehresmann fibration theorem to $\bar{F}_{|\bar{F}^{-1}(\partial D_\eta\times D_\delta)}$. 
Define $X_{D_\eta}:=\bar{F}^{-1}(D_\eta\times D_\delta)\cap (\overline{B}_\epsilon\times D_\delta)$ and consider 
$p:X_{D_\eta}\to D_\delta$, the restriction to $X_{D_\eta}$ of the projection to the second factor; then
$p^{-1}(s)=X_{D_\eta,s}$ for any $s\in D_\delta$. Using the first assertion, it is easy to show that $X_{D_\eta}$ is a
manifold with corners and that $p$ is a proper differentiable map, whose restriction to the boundary and corners is 
submersive; using a version of Ehresmann fibration theorem for manifolds with corners the third assertion follows.  

Define $Y\subset S_\epsilon\times D_\delta$ to be the set of pairs $(x,s)$ such that $F_{|s}(x)\neq 0$ and either $F_{|s}$
is critical at $x$ or $F_{|s}^{-1}(F_{|s}(x))$ is not transversal to $S_\epsilon$ at $x$. Identifying $\CC^n\times\CC$ with
$\RR^{2n+2}$ it is easy to see that the closure $\bar{Y}$ of $Y$ is a real analytic subset. The following claim 
implies the proposition (perhaps shrinking $\eta$): the intersection $\bar{Y}\cap F^{-1}(0,0)$ is empty. Suppose the 
contrary: 
let $(x,0)\in \bar{Y}\cap F^{-1}(0,0)$; by the Curve Selection Lemma (see~\cite{Mi}, \S3) there is an analytic path 
$\alpha:(-\nu,\nu)\to\bar{Y}$ such that $\alpha(-\nu,\nu)\subset Y$ and $\alpha(0)=(x,0)$. Define $\gamma:(-\nu,\nu)\to\CC$
to be the second component of the composition $\overline{F}\comp\alpha$, and choose $\mu<\nu$ so that the statement of 
Lemma~\ref{controlext} holds; let $\Psi$ be the family of diffeomorphisms predicted by Lemma~\ref{controlext}. Consider
a sequence $\{t_n\}\subset (\mu,\mu)$ convergent to $0$; define $x_n$ to be the first component of $\alpha(t_n)$ in $\CC^n\times\CC$,
by the formula $(x_n,\gamma(t_n))=\alpha(t_n)$, and 
$y_n:=\psi_{t_n}^{-1}(x_n)$. As $\{x_n\}$ and $\{\psi_{t_n}\}$ converge to $x$ and $\mathrm{Id}_{\CC^n}$ respectively, we deduce 
that $\{y_n\}$ converges to $x$. If $F_{|\gamma(t_n)}$ is singular at $x_n$ then $f=F_{|\gamma(t_n)}\comp\psi_{t_n}$ is 
singular at $y_n$; then $f(y_n)=F_{|\gamma(t_n)}(x_n)=0$, which is not true; therefore $F_{|\gamma(t_n)}$ is not singular 
at $x_n$. Then $F_{|\gamma(t_n)}^{-1}(F_{|\gamma(t_n)}(x_n))$ is not transversal to $S_\epsilon$ at $x_n$, which means that
$T_{x_n}F_{|\gamma(t_n)}^{-1}(F_{|\gamma(t_n)}(x_n))\subset T_{x_n}S_\epsilon$; this implies that 
$T_{y_n}f^{-1}(f(y_n))\subset d\psi_{t_n}^{-1}(x_n)(T_{x_n}S_\epsilon)$. Taking a subsequence we can assume that the 
sequence 
$T_{y_n}f^{-1}(f(y_n))$ converges to a linear subspace $T\subset T_xS_\epsilon$. On the other hand we have fixed Whitney
stratifications of $\CC$ and of an open neighbourhood $U$ of the origin of $\CC^n$ (containing $U\setminus f^{-1}(0)$ as an stratum)
such that $f$ satisfies the Thom 
$A_f$-condition respect to them and that $X\trans S_\epsilon$ for any stratum $X$ of $U$. 
Let $X$ be the stratum containing $x$, as $f(x)=0$ and $U\setminus f^{-1}(0)$ is an stratum we have that 
$X\subset f^{-1}(0)$; hence $\mathrm{ker}(df_{|X}(x))=T_xX$. By Thom $A_f$-condition 
$T_xX\subset \mathrm{lim}T_{y_n}f^{-1}(f(y_n))=T$; as $T$ is included in $T_xS_\epsilon$, we contradict the transversality
$X\trans S_\epsilon$. This proves the claim.
\end{proof}

\section{Generalized Jet-Spaces}

The Morsification Theorem for isolated singularities can be proved in the following way (which is not the easiest but
has the virtue of generalising to our setting):                     
the subset $Z_1\subset J^2(\CC^n,\CC)$ consisting of singular $2$-jets (jets whose linear part vanish) is a closed analytic
subset of codimension $n$; the subset $Z_2\subset Z_1$ of singular $2$-jets germs with degenerate Hessian is closed 
analytic of codimension bigger than $n$. On the other hand any germ $f\in\ocn$ can be approximated by germs $g$ such that
the $2$-jet extension $j^2g:\CC^n\to J^2(\CC^n,\CC)$ is transversal to $Z_1$ and $Z_2$ in a neighbourhood of the origin; 
therefore the singularities that an approximation $g$ can have close to the origin should have non-degenerate Hessian 
(which means being of Morse type).

Consider an ideal $I\subset\ocn$. When trying to generalise 
the classical morsification method to functions of finite codimension with respect to $I$
we meet (among others) the following new difficulties: consider $f\in I$ such that 
$c_{I,e}(f)=\infty$; it is not clear a priori
which singularities may have a generic deformation within the ideal. On the other hand,
due to the presence of non-isolated singularities, if we work with ordinary 
$m$-jet spaces the spaces parametrising singularity types (analogous to $Z_1$ and 
$Z_2$) will be of arbitrarily large codimension as $m$ increases, which would make 
the above method collapse. The idea to overcome these problems will be to use a presentation of the ideal $I$
by generators and relations to define a sort of generalized $m$-jet spaces, in which the
codimension the varieties parametrising the relevant singularity types remain stable
as $m$-increases.
 
We need to introduce the concept of $\infty$-jet spaces and give it an infinite dimensional analytic 
structure. 

Let $U\in\CC^n$ be an open subset. For any $m<\infty$ the $m$-th jet-space $J^m(U,\CC^r)$ is a vector bundle over $U$ with projection mapping 
$pr_m:J^m(U,\CC^r)\to U$. There is a natural analytic vector bundle epimorphism $pr^m_l:J^m(U,\CC^r)\to J^l(U,\CC^r)$ for any $m\geq l$. 
The set $J^\infty(U,\CC^r):=\coprod_{x\in U}\calO_{U,x}^r$ is clearly the projective limit of the system formed by the $J^m(U,\CC^r)$'s and the $\pi^m_l$'s.
For any $m$ there is a projection mapping $\pi^\infty_m:J^\infty(U,\CC^r)\to J^m(U,\CC^r)$.  

Fix a coordinate system $\{x_1,...,x_n\}$ in $\CC^n$; denote by $\CC\{x\}$ the ring of convergent power series in 
$\{x_1,...,x_n\}$. There is a bijection
\[\tau_\infty:U\times\CC\{x\}^r\to J^\infty(U,\CC^r)\] which assigns to $(x,(f_1,...,f_r))$ the unique $r$-tuple of germs 
$(g_1,...,g_r)\in\calO_{\CC^n,x}^r$ such that the Taylor expansion of $g_i$ at $x$ is $f_i$. Passing to $m$-jets this defines analytic vector bundle 
trivialisation 
\begin{equation}
\label{r1trivi}
\tau_m:U\times (\CC\{x\}/\mm^{m+1})^r\to J^m(U,\CC^r)
\end{equation}
for any $m<\infty$.

A subset $C\subset J^\infty(U,\CC^r)$ is said
$k$-{\em determined} if it satisfies $(\pi^\infty_k)^{-1}(\pi^\infty_k(C))=C$; the subset $C$ is said to be 
a {\em $k$-determined closed analytic subset of} $J^\infty(U,\CC^r)$ if it is $k$-determined 
and $\pi^\infty_k(Y)$ is a closed analytic subset of $J^k(U,\CC^r)$. Any $k$-determined (closed analytic) subset is also $m$-determined 
(closed analytic) for any $m\geq k$. A $k$-determined {\em locally closed analytic} subset is the difference between two 
$k$-determined closed analytic subsets.

Consider in each $J^m(U,\CC^r)$ the transcendental topology; we endow $J^\infty(U,\CC^r)$ with the final topology for the family of projections 
$\pi^\infty_m$. A family 
$\{C_j\}_{j\in J}$ of subsets of $J^\infty(U,\CC^r)$ is locally finite if for any $x\in J^\infty(U,\CC^r)$ there exists a positive integer $m$ and a 
neighbourhood $U$ of 
$\pi^\infty_m(x)$ in $J^m(U,\CC^r)$ such that $(\pi^\infty_m)^{-1}(U)$ only meets finitely many $C_j$'s. Therefore if each $C_j$ is finite-determined, choosing
$m$ high
enough, the union $\cup_{j\in J}C_j$ looks locally a $m$-determined subset. This motivates the following definition: 
a {\em closed analytic} subset of $J^\infty(U,\CC^r)$ is 
any union of a locally finite collection of finite-determined closed analytic subsets of $X^\infty$; a {\em locally closed analytic} subset of 
$X^\infty$ is the difference between two closed analytic subsets.

Let $\tilde{I}$ be an ideal sheaf defined on a neighbourhood $U$ of the origin, having $I$ as stalk at the origin. The 
natural candidates to substitute the ordinary jet-spaces considered in the isolated-singularity case are the following 
spaces: define the sets
\[J^\infty(U,\tilde{I}):=\coprod_{x\in U}\tilde{I}_x\quad J^m(U,\tilde{I}):=\coprod_{x\in U}\tilde{I}_x/\mm_x^{m+1}\cap\tilde{I}_x\] 
for any positive integer $m$. For 
any $\infty\geq m\geq k\geq 0$ we consider the projection mapping $\pi^m_k:J^m(U,\tilde{I})\to J^k(U,\tilde{I})$. The set of spaces 
$\{J^m(U,\tilde{I})\}_{m<\infty}$ together with the mappings $\{pr^m_k\}_{k<m<\infty})$ forms a projective system of sets whose limit is 
$J^\infty(U,\tilde{I})$. There are natural projection mappings $\pi^\infty_m:J^\infty(U,\tilde{I})\to J^m(U,\tilde{I})$ satisfying 
$\pi^\infty_l=\pi^m_l\comp\pi^\infty_m$.
The concept of $k$-determined subset of $J^\infty(U,\tilde{I})$ is defined analogously to the case
of systems of analytic varieties. For any $m\leq\infty$ we consider the natural projection $pr_m:J^m(U,\tilde{I})\to U$; its fibre over $x\in U$ is the vector 
space $\tilde{I}_x/\tilde{I}_x\cap\mm_x^{m+1}$. 

For any $x\in U$ we define the function $\mu_x:\ZZ_{>0}\to\ZZ$ by the formula
\begin{equation}
\label{hilbsam}
\mu_x(m):=\dim_\CC(\tilde{I}_x/\tilde{I}_x\cap\mm_x^{m+1}).
\end{equation}
The function $\mu_x$ is Zarisky-lower semicontinous in $W$ (see~\cite{Bo2}). There exists a stratification of $U$ (which is called the {\em Zarisky-Samuel 
stratification} with respect to $\tilde{I}$) by Zarisky locally closed analytic subsets, which is the minimal partition such that $\mu_x=\mu_y$ for any 
$x$, $y$ in the same stratum. We will refer to the strata as the $\tilde{I}$-strata of $U$, and we will denote them by $\Sigma_0,...,\Sigma_s$, where 
$\Sigma_s$ is the stratum containing the origin, the stratum $\Sigma_0$ is the complement of the zero set of $\tilde{I}$ and, if 
$\Sigma_i\subset\overline{\Sigma}_j$ then $i\geq j$.

\begin{notation}{\em
Consider an analytic function $f:V\to\CC$, with $V$ an open subset of $\CC^n$. Given any $m\leq\infty$ there is an associated {\em jet extension}
$j^mf:J^m(V,\CC)$ assigning to $x\in V$ the class of $f_x$ in $\calO_{\CC^n,x}/\mm_x^{m+1}$ (we adopt the convention $\mm_x^\infty=(0)$).}
\end{notation}
 
We define a certain kind of analytic atlas on $J^m(U,\tilde{I})$ in the following way:

Consider an open subset $V$ of $U$ and a set $\calH=\{h_1,...,h_s\}\subset\Gamma(U,\tilde{I})$ which generate $\tilde{I}_{|V}$. Define 
\begin{equation}
\label{puresen1} 
\varphi_{\calH}:\calO_V^s\to\calO_{|V}
\end{equation}
by the formula $\varphi_\calH(f_1,...,f_s):=\sum_{i=1}^sf_ih_i$. For any $m\leq\infty$, taking $m$-jets we obtain a mapping
\begin{equation}
\label{puresen2}
j^m\varphi_\calH:J^m(V,\CC^s)\to J^m(V,\CC).
\end{equation}
where $j^m\varphi_\calH(j^mf_1(x),...,j^mf_s(x)):=\sum_{i=1}^sj^m(f_ih_i)(x)$ for any $(f_1,...,f_s)\in(\calO_{V,x})^s$ and $x\in V$. The mapping $j^m\varphi$ 
has $J^m(V,\tilde{I})$ as image, and it is a homomorphism of analytic vector bundles if $m<\infty$. 

\begin{definition}
For any $m\leq\infty$,  a {\em chart} of 
$J^m(U,\tilde{I})$ is a surjective mapping of the form
\begin{equation}
\label{puresen3}
j^m\varphi_\calH:J^m(V,\CC^s)\to J^m(V,\tilde{I}).
\end{equation} 
for a certain open subset $V$ (which will be called the {\em base subset} of the chart) and a set $\calH$ of holomorphic functions on $V$ generating $\calH$ of
$\tilde{I}_{|V}$. The {\em canonical analytic atlas} of $J^m(U,\tilde{I})$ is the collection of all the charts of $J^m(U,\tilde{I})$.
\end{definition}

Given another open subset $V'$ and a set of generators $\calH'=\{h'_1,...,h'_{s'}\}$ of $\tilde{I}_{|V'}$, there exists an open covering $\{V''_l\}_{l\in L}$ 
of $V\cap V'$ such that, given any $l\in L$, for any $i\leq s'$ we have analytic functions 
$\theta_{i,1},...,\theta_{i,s}$ defined on $V''_l$ such that $h'_i:=\sum_{j=1}^s\theta_{i,j}h_j$. Define the $\calO_{V''_l}$-module homomorphism 
$\theta:\calO_{V''_l}^{s'}\to\calO_{V''_l}^{s}$ by the formula $\theta(f_1,...,f_{s'})=(f_1,...,f_{s'})M_\theta$ where $M_\theta$ is the matrix whose 
$(i,j)$-entry is $\theta_{i,j}$. For any $m\leq\infty$ let
\begin{equation}
\label{transi}
j^m\theta:J^m(V''_l,\CC^{s'})\to J^m(V''_l,\CC^{s})
\end{equation}
be the associated mapping of jet-spaces, which is analytic for $m<\infty$. We have clearly the compatibility relation
\begin{equation}
\label{compati}
j^m\varphi_{\calH'}=j^m\varphi_{\calH}\comp j^m\theta.
\end{equation}
The mapping~(\ref{transi}) is called {\em a transition function} between the charts $j^m\varphi_\calH$ and $j^m\varphi_{\calH'}$. Notice that unlike in the 
case of manifolds transition functions between charts are not globally defined in the intersection of the two charts and need not be unique. 

As we work locally around the origin of $\CC^n$ we can take $U$ small enough so that there is a set $\calG=\{g_1,...,g_r\}$ of analytic functions defined on 
$U$ which generate $\tilde{I}_x$ for any $x\in U$. We will denote the homomorphism $\varphi_\calG$ simply by  
\begin{equation}
\label{presen1} 
\varphi:\calO_U^r\to\calO_U,
\end{equation}
and, for any $m\leq\infty$, the mapping $j^m\varphi_\calG$ by
\begin{equation}
\label{presen2}
J^m(U,\CC^r)\stackrel{j^m\varphi}{\longrightarrow}J^m(U,\CC).
\end{equation}

\begin{notation}
{\em Let $X\subset U$ be any subset. For any $m\leq\infty$, we will denote by $J^m(X,\tilde{I})$ or $J^m(X,\tilde{I})_{|X}$ the inverse image of $X$ under 
$pr_m:J^m(U,\tilde{I})\to U$. An analogous notation works for the jet space $J^m(U,\CC^r)$.}
\end{notation} 

\begin{remark}
\label{restricciones} 
For $m<\infty$ and any $\tilde{I}$-stratum, the restriction $j^m\varphi_{|J^m(\Sigma_i,\CC^r)}$ is a constant rank homomorphism of trivial analytic vector 
bundles over $\Sigma_i$, whose image is $J^m(\Sigma_i,\tilde{I})$. This gives a natural structure of trivial analytic vector bundle to 
$J^m(\Sigma_i,\tilde{I})$, for $m<\infty$ (see~\cite{Bo2} for an application). 
\end{remark} 

For our purposes it is convenient to consider subsets of $J^\infty(U,\tilde{I})$ parametrising germs with certain geometric properties, just as the subset
$Z_1\subset J^2(U,\CC)$ considered at the beginning of the section parametrises Morse singularities. To have a geometric understanding of these subsets we 
will look at their inverse image by the charts of the canonical analytic atlas of $J^m(U,\tilde{I})$.

For any $m<\infty$ we give to $J^m(U,\tilde{I})$ the final topology for the set of all charts of the canonical 
analytic atlas. As the transition functions are continuous the topology on $J^m(U,\tilde{I})$ is just the final topology for the mapping 
$j^m\varphi:J^m(U,\CC^r)\to J^m(U,\tilde{I})$. For any $s$ we give to $J^\infty(U,\tilde{I})$ the topology 
obtained viewing it as projective limits of the systems of topological spaces $\{J^m(U,\tilde{I})\}_{m\in\NN}$.
It is easy to check that the charts for $m=\infty$ are continuous. 

\begin{remark}
As $j^m\varphi_{|J^m(\Sigma_i,\CC^r)}$ is a submersion for any $m<\infty$ and any $\tilde{I}$-stratum, the restriction of the topology of $J^m(U,\tilde{I})$ 
to $J^m(\Sigma_i,\tilde{I})$ coincides with the restriction of the topology of $J^m(\Sigma_i,\CC)$ to $J^m(\Sigma_i,\tilde{I})$. Therefore the topology 
considered in~\cite{Bo2} for 
$J^m(\Sigma_i,\tilde{I})$ is the restriction of the topology considered here for $J^m(U,\tilde{I})$.
\end{remark}

\begin{definition}
Consider $0<k\leq m\leq\infty$.
A $k$-determined subset $C\subset J^m(U,\tilde{I})$ is {\em closed analytic} if for any chart $\varphi_\calH$ with base subset $V$ we have that 
$(j^k\varphi_\calH)^{-1}(\pi^m_k(C))$ is 
$m$-determined closed analytic in $J^k(V,\CC^m)$. A {\em locally closed} $k$-determined analytic subset of $J^m(U,\tilde{I})$ is the difference 
between two
$m$-determined closed analytic subsets. A {\em closed analytic} subset of $J^\infty(U,\tilde{I})$ is a locally finite union of finite determined closed 
analytic subsets of $J^\infty(U,\tilde{I})$. A {\em locally closed analytic} subset of $J^\infty(U,\tilde{I})$ is the difference between two closed 
analytic subsets.
\end{definition} 

Consider a $k$-determined subset of $J^m(U,\tilde{I})$. Given any $k\leq k'\leq m$ it is easy to check that it is a $k$-determined (locally) closed analytic 
subset if and only if it is a $k'$-determined (locally) closed analytic subset. Let $Z$ be a closed $m$-determined subset of $J^m(U,\tilde{I})$, and let 
$j^m\varphi_\calH:J^m(V,\CC^s)\to J^m(U,\tilde{I})$ and 
$j^m\varphi_{\calH'}:J^m(V',\CC^{s'})\to J^m(U,\tilde{I})$ be two charts. By the compatibility~(\ref{compati}), if $(j^k\varphi_{\calH})^{-1}(Z)$ is a  
closed analytic subset, then $(j^k\varphi_{\calH'})^{-1}(Z_{|V''_l})$ is closed analytic in $J^m(V''_l,\CC^{s'})$. Therefore to prove the analyticy of a 
$m$-determined subset it is enough to check the condition for a set of charts whose base subsets cover $U$. In particular it is enough to check that 
$(j^m\varphi)^{-1}(Z)$ is analytic. 

\begin{definition}
A subset $Y\subset J^m(U,\CC^r)$ satisfying $(j^m\varphi)^{-1}(j^m\varphi(Y))=Y$ is called $j^m\varphi$-{\em saturated}. 
\end{definition}
 
Closed analytic subsets of $J^m(U,\tilde{I})$ are in a bijective correspondence with $j^m\varphi$-saturated closed analytic subsets of $J^m(U,\CC^r)$. 

\begin{definition} 
A closed analytic subset of $J^\infty(U,\tilde{I})$ is {\em irreducible} if it cannot be expressed as the union of two closed analytic 
subsets of $J^\infty(U,\tilde{I})$ not containing it.  
\end{definition}

It follows that any irreducible closed analytic subset is finite-determined.
 
We shrink $U$ so that its closure is compact and contained in an open subset where $\tilde{I}$ is defined. Then the following uniform Artin-Rees theorem 
holds (see~\cite{BM} for a proof): there exists $\lambda\in\ZZ_{\geq 0}$ such that 
\begin{equation}
\label{artrees}
\tilde{I}_x\cap\mm_x^{m+\lambda}\subset\mm_x^{m}\tilde{I}_x
\end{equation}
for any $x\in U$ and any $m$. The minimal $\lambda$ so that the last inclusion holds is called the {\em Uniform Artin-Rees constant}.

\begin{lema}
\label{irraux}
Let $C\subset J^\infty(U,\tilde{I})$ be a $k$-determined closed analytic subset. For any irreducible component $K$ of $(j^{k}\varphi)^{-1}(\pi^\infty_{k}(C))$ 
we have that $(\pi^{\lambda+k}_k)^{-1}(K)$ is $j^{k+\lambda}\varphi$-saturated.
\end{lema}
\begin{proof}
Set $m=k+\lambda$ and 
\[A:=(\pi^m_k)^{-1}((j^k\varphi)^{-1}(\pi^\infty_k(C))).\] 
Let $K':=(\pi^m_k)^{-1}(K)$. Clearly $A$ is closed analytic and $K'$ is one of its irreducible 
components. We have to show that 
\begin{equation}
\label{wish}
K'_x+\mathrm{ker}(j^m\varphi_x)\subset K'_x
\end{equation}
for any $x\in U$ (where $K'_x=K'\cap pr_m^{-1}(x)$, and $j^m\varphi_x$ is the restriction of $j^m\varphi$ to $J^m(U,\CC^r)_x$). 

An element of $\mathrm{ker}j^m\varphi_x$ is the class in $J^m(U,\CC^r)_x$ of a $r$-tuple
$(f_1,...,f_r)\in\calO_{U,x}^r$ such that $\varphi_x(f_1,...,f_r)$ belongs to $\tilde{I_x}\cap\mm_x^{m+1}$. By Uniform Artin-Rees we have that 
$\varphi_x(f_1,...,f_r)$ belongs to $\mm_x^{k+1}\tilde{I}_x$. Therefore there exists 
$(h_1,...,h_r)\in\mm_x^{k+1}\calO_{U_x}^r$ so that $\varphi_x(h_1,...,h_r)=\varphi_x(f_1,...,f_r)$. We deduce that $(f_1,...,f_r)$ belongs to 
$\mm_x^{k+1}\calO_{U,x}^r+\mathrm{ker}(\varphi_x)$. Hence 
\begin{equation}
\label{conclusion0}
\mathrm{ker}(j^m\varphi_x)\subset\mm_x^{k+1}\calO_{U,x}^r/\mm_x^{m+1}\calO_{U,x}^r+\pi^\infty_m(\mathrm{ker}(\varphi_x)).
\end{equation}

As $\varphi^{-1}(C)$ is $k$-determined we have 
\begin{equation}
\label{inva1}
A_x+\mm_x^{k+1}\calO_{U,x}^r/\mm_x^{m+1}\calO_{U,x}^r\subset A_x
\end{equation}
for any $x\in U$. The set 
$E:=\coprod_{x\in U}\mm_x^{k+1}\calO_{U,x}^r/\mm_x^{m+1}\calO_{U,x}^r$ is a trivial vector sub-bundle of $J^m(U,\CC^r)$. Let $r_1$ be the rank of $E$. Consider
an analytic trivialisation $\tau_1:G_1\times U\to E$, where $G_1\cong\CC^{r_1}$. View $G_1$ as an additive analytic group. We have an analytic action
\[G_1\times J^m(U,\CC^r)\stackrel{\sigma_1}{\longrightarrow} J^m(U,\CC^r)\]
defined by $\sigma_1(g,f):=\tau_1(g)+f$. By~(\ref{inva1}) the subset $A$ is left invariant by the action.

For being $\tilde{I}$ coherent we can take $U$ small enough so that we have an exact sequence as follows 
\begin{equation}
\label{resolution}
\calO_U^s\stackrel{\psi}{\longrightarrow}\calO_U^r\stackrel{\varphi}{\longrightarrow}\tilde{I}\to 0.
\end{equation}
Taking $m$-jets we obtain a homomorphism $j^m\psi:J^m(U,\CC^s)\to J^m(U,\CC^r)$ of trivial analytic vector bundles. The exactness of (\ref{resolution}) implies
\begin{equation}
\label{exactitud}
\pi^\infty_m(\mathrm{ker}\varphi_x)=\pi^\infty_m(\mathrm{im}\psi_x)=\mathrm{im}j^m\psi_x
\end{equation}
for any $x\in U$. This, together with the $j^m\varphi$-saturation of $A$ implies
\begin{equation}
\label{inva2}
A_x+\mathrm{im}(j^m\psi_x)\subset A_x.
\end{equation}
Let $r_2$ be the rank of the vector bundle $J^m(U,\CC^s)$. Consider an analytic trivialisation $\tau_2:G_2\times U\to J^m(U,\CC^s)$, where $G_1\cong\CC^{r_2}$.
View $G_2$ as an additive analytic group. We have an analytic action
\[G_2\times J^m(U,\CC^r)\stackrel{\sigma_2}{\longrightarrow} J^m(U,\CC^r)\]
defined by $\sigma_2(g,f):=j^m\psi\comp\tau_2(g)+f$. By~(\ref{inva2}) the subset $A$ is left invariant by the action.

Define the additive analytic group $G:=G_1\oplus G_2$ and the action 
\begin{equation}
\label{accion}
\sigma:G\times J^m(U,\CC^r)\to J^m(U,\CC^r)
\end{equation}
by $\sigma:=\sigma_1\oplus\sigma_2$. The subset $A$ is invariant by the action $\sigma$. As $G$ is irreducible we conclude that each of the irreducible 
components of $A$ are also invariant by the action. This, together with~(\ref{conclusion0})~and~(\ref{exactitud}) imply~(\ref{wish}).
\end{proof}

\begin{prop}
Let $C\subset J^\infty(U,\tilde{I})$ be a $k$-determined closed analytic subset. The following are equivalent:
\begin{enumerate}
\item The set $C$ is irreducible.
\item The set $C$ is finite determined and $(j^m\varphi)^{-1}(\pi^\infty_m(C))$ is irreducible for a certain $m\geq k$.
\item The set $C$ is finite determined and $(j^m\varphi)^{-1}(\pi^\infty_m(C))$ is irreducible for any $m\geq k$.
\end{enumerate}
Moreover, for any closed analytic subset $C$ of $J^\infty(U,\tilde{I})$ there exists a unique irredundant locally finite decomposition of $C$ in 
irreducible closed analytic subsets. If $C$ is $k$-determined, then each of its irreducible components is $k+\lambda$-determined. 
\end{prop}
\begin{proof}
For any $m\geq k$ the morphism $\pi^m_k:J^m(U,\CC^r)\to J^k(U,\CC^r)$ is a vector bundle epimorphism. As 
$(j^m\varphi)^{-1}(\pi^\infty_m(C))=(\pi^m_k)^{-1}((j^k\varphi)^{-1}(\pi^\infty_k(C)))$, the irreducibility of $(j^m\varphi)^{-1}(\pi^\infty_m(C))$ and 
$(j^k\varphi)^{-1}(\pi^\infty_k(C))$ are equivalent. This shows $(2)\Leftrightarrow (3)$.

Let us prove  $(3)\Rightarrow (1)$. Let $C$ satisfying the property $(3)$. Suppose $C=C_1\cup C_2$, where $C_1$ and $C_2$ are closed analytic. Let $C$ be 
$k$-determined. Using the ideas of Lemma~\ref{irraux} and the fact that $C_1$ and $C_2$ are locally finite-determined closed analytic sets, it is easy to show 
that $C_1$ and $C_2$ are $k+\lambda$-determined. From here proving $(1)$ is straightforward. 

Suppose $C$ is $k$-determined. Let $\{K_j\}_{j\in J}$ be the decomposition in irreducible components of $(j^k\varphi)^{-1}(\pi^\infty_k(C))$. For any $j\in J$
the set $K'_j:=(\pi^{k+\lambda}_k)^{-1}(K_j)$ is a irreducible component of $(j^{k+\lambda}\varphi)^{-1}(\pi^\infty_{k+\lambda}(C))$. By Lemma~\ref{irraux} 
the set $K'_j$ is $j^{k+\lambda}\varphi$-saturated. Therefore each $C_j:=(\pi^\infty_{k+\lambda})^{-1}j^{k+\lambda}\varphi(K')$ is a $k+\lambda$-determined 
closed subset of $J^\infty(U,\tilde{I})$ and $\cup_{j\in J}C_j=C$. As we have already shown $(2)\Rightarrow (1)$ we have that each $C_j$ is irreducible.
This shows the existence of a unique irredundant decomposition in irreducible components for finite-determined subsets, and also proves $(1)\Rightarrow (2)$.

The existence and unicity of irredundant decompositions in irreducible components for non-necessarily finite-determined  closed analytic subsets is easily 
deduced from the same property in the special case of finite-determined subsets. 
\end{proof}

Let $C$ be a $k$-determined closed analytic subset of $J^\infty(U,\tilde{I})$. The irreducible components of $C$
are not $k$-determined in general. This can be seen already in the simplest examples:

\begin{example}{\em
Consider $I\subset\CC\{z\}$ generated by $z^2$. Let $\varphi:\calO_\CC\to\calO_\CC$ be defined by $\varphi(f):=z^2f$. The set $\coprod_{x\in\CC}\mm_x$ is a 
$1$-determined closed subset of $J^\infty(\CC,\tilde{I})$. It has two irreducible components $C_1=J(\CC,\tilde{I})_0$ and 
$C_2=\mm_0^3\cup\coprod_{x\neq 0}\mm_x$. Clearly $C_2$ is $2$-determined, but not $1$-determined.}    
\end{example}

\begin{definition}
The {\em codimension} $\mathrm{codim}(C)$ of a $k$-determined irreducible locally closed subset of $J^\infty(U,\tilde{I})$ is the 
codimension of $(j^k\varphi)^{-1}(pr^\infty_k(C))$ in $J^k(U,\CC^r)$.
\end{definition}

In the situation of the last definition, if $m\geq k$, then clearly $\mathrm{codim}(C)$ equals the codimension of $(j^m\varphi)^{-1}(pr^\infty_m(C))$ in
$J^m(U,\CC^r)$. It is also clear that if $C\subset C'$ are closed analytic subsets of $J^\infty(U,\tilde{I})$, and $C'$ is irreducible, then either $C=C'$ or 
$\mathrm{codim}(C_i)>\mathrm{codim}(C')$ for any irreducible component $C_i$ of $C$.

Actually the codimension of $C$ does not depend on the chosen chart:

\begin{lema}
\label{indepcarta}   
Let $C$ be a $k$-determined irreducible locally closed subset of $J^\infty(U,\tilde{I})$. Given any open subset $V\subset U$ and any system of $s$ generators 
$\calH$ of $\tilde{I}_{|V}$, the codimension of $(j^k\varphi_\calH)^{-1}(pr^\infty_k(C))$ in $J^k(V,\CC^r)$ equals $\mathrm{codim}(C)$.
\end{lema}
\begin{proof}
Consider another set of generators $\calH'=\{h'_1,...,h'_{s'}\}$ of $\tilde{I}_{|V'}$. If $\calH'$ contains $\calH$, in order to provide transition functions 
between the charts associated to $\calH$ and $\calH'$, the functions $\theta_{i,j}$ can be chosen to be $\theta_{i,j}=0$ if $j\neq i$ and 
$\theta_{i,i}=1$ for $i\leq s$. Note that the $\theta_{i,j}$'s are defined in the whole $V\cap V'$ (the set where both $\calH$ and $\calH'$ are defined). 
In this case the transition function $j^k\theta$ is a submersion for any $k<\infty$. Therefore, as 
$(j^k\varphi_{\calH'})^{-1}(\pi^\infty_k(C))=(j^k\theta)^{-1}((j^k\varphi_{\calH})^{-1}(\pi^\infty_k(C)))$,
the codimension of $(j^k\varphi_{\calH'})^{-1}(\pi^\infty_k(C))$ in $J^k(V,\CC^{s'})$ equals the codimension of $(j^k\varphi_{\calH})^{-1}(\pi^\infty_k(C))$ 
in $J^k(V,\CC^{s})$. 

Given any two sets of generators of $\tilde{I}_{|V'}$ their union gives a set of generators containing both of them. Therefore the codimension does not depend
on the set of generators giving rise to the chart.
\end{proof}

The following natural fact is easily deduced as a consequence of the existence of decomposition in irreducible components:

\begin{lema}
\label{topclo}
The topological closure of a locally closed analytic subset of $J^\infty(U,\tilde{I})$ is closed analytic.
\end{lema}
\begin{proof}
Let $C$ be a locally closed analytic subset of $J^\infty(U,\tilde{I})$. We have $C=A\setminus B$, with $A$ and $B$ closed analytic subsets of 
$J^\infty(U,\tilde{I})$. Let $A=\cup_{j\in J}A_j$ be a decomposition of $A$ in irreducible components. Let $J'\subset J$ be the subset of consisting of the 
indices $j$ such that $A_j$ is not included in $B$. By the locally finiteness of the decomposition in irreducible components we have 
$\overline{C}=\cup_{j\in J'}\overline{A_j\setminus B}$. Therefore we are reduced to the case in which $A$ is irreducible. We can cover $J^\infty(U,\tilde{I})$
by open subsets $\{U_l\}_{l\in L}$ such that $A$, $U_l$ and $B\cap U_l$ are $m_l$-determined for a certain integer $m_l$. This means that for any $l\in L$ 
there exists an open subset $V_l\subset U$, a $j^{m_l}\varphi$-saturated open subset  $W_l\subset J^{m_l}(V_l,\CC^r)$ and $j^{m_l}\varphi$-saturated closed 
analytic subsets $A_l$ and $B_l$ of $J^{m_l}(V_l,\CC^r)$ such that 
\[(\pi^\infty_{m_l})^{-1}(j^{m_l}\varphi(W_l))=U_l,\]
\[(\pi^\infty_{m_l})^{-1}(j^{m_l}\varphi(A_l))=A\cap U_l,\]
\[(\pi^\infty_{m_l})^{-1}(j^{m_l}\varphi(B_l))=B\cap U_l.\]
As $A$ is irreducible any irreducible, component of $B\cap A$ has strictly bigger codimension. Therefore $\overline{A_l\setminus B_l}=A_l$, and consequently 
(by the definition of the topology in $J^\infty(U,\tilde{I})$) the set $A\cap U_l$ equals $\overline{C}\cap U_l$. It follows that $A=\overline{C}$.
\end{proof}

After this lemma it makes sense to define: 

\begin{definition}
A locally closed analytic subset is {\em irreducible} if its closure is irreducible.
\end{definition}

\begin{lema}
\label{singlocus}
Let $C$ be a $k$-determined irreducible locally closed analytic subset of $J^\infty(U,\tilde{I})$. Given any open subset $V\subset U$ 
and any system of generators $\calH$ of $\tilde{I}_{|V}$ we consider 
$A:=(j^{k+\lambda}\varphi_\calH)^{-1}(\pi^\infty_{k+\lambda}(C))$. Then the set $\mathrm{Sing}(A)$ is $j^{k+\lambda}\varphi_\calH$-saturated.
\end{lema} 
\begin{proof}
In the proof of Lemma~\ref{irraux} we have shown that $A$ is invariant by the group action $\sigma:G\times J^{k+\lambda}(V,\CC^s)\to J^{k+\lambda}(V,\CC^s)$ 
(where $s$ is the number of generators of $\calH$). Therefore for 
any $g\in G$ we have $\sigma(g,\mathrm{Sing}(A))=\mathrm{Sing}(A)$. This, together with~(\ref{conclusion0})~and~(\ref{exactitud}) implies  
\begin{equation}
\mathrm{Sing}(A)+\mathrm{ker}(j^{k+\lambda}\varphi_{\calH,x})\subset\mathrm{Sing}(A)
\end{equation}
for any $x\in V$.
\end{proof}

\begin{definition}
\label{smoothness}
Let $C$ be a locally closed analytic in $J^\infty(U,\tilde{I})$ and $f\in C$. We say that $C$ {\em irreducible at} $f$ there is only one irreducible component 
of $C$ containing $f$. We say that $C$ is {\em smooth} at $f$ if
\begin{enumerate}
\item It is irreducible at $f$.
\item Let $C'$ be the unique irreducible component of $C$ containing $f$; let $C'$ be $k$-determined. For any open subset $V\subset U$ containing 
$pr_\infty(f)$, any system of generators $\calH$ of $\tilde{I}_{|V}$ and any $m\geq k+\lambda$ the locally closed analytic 
$(j^{m}\varphi_\calH)^{-1}(\pi^\infty_{m}(C'))$ is smooth at any $h$ such that $j^{m}(h)=\pi^\infty_{m}(f)$. 
\end{enumerate}
\end{definition}

Using an argument similar to the proof of Lemma~\ref{indepcarta} it is easy to check that it is enough to check smoothness at a single chart.
\begin{remark}{\em
Taking into account Lemma~\ref{singlocus} it is easy to check that it is enough to check the second condition for smoothness for a particular $m\geq k+\lambda$
and a particular $h$ such that $j^{m}(h)=\pi^\infty_{m}(f)$. Moreover, for any locally closed analytic subset $C$ the set $\mathrm{Sing}(C)$ of singular points
is closed analytic in $C$, and it is $k+\lambda$-determined if $C$ is $k$-determined.}
\end{remark}

\section{The filtration by extended codimension} 
The natural generalisation of the subvariety of $Z_1\subset J^2(U,\CC)$ parametrising Morse singularities is the set of 
germs in $J^m(U,\tilde{I})$ parametrising sigularities of extended codimension equal to $1$, this motivates the following:

\begin{definition}
\label{prisa} 
Suppose that $m$ is either a non-negative integer or $\infty$; define 
\[C_m:=\{f\in J^\infty(U,\tilde{I}):c_{\tilde{I}_{\pi^\infty(f)},e}(f)\geq m\}\]
\[\dot{C}_m:=\{f\in J^\infty(U,\tilde{I}):c_{\tilde{I}_{\pi^\infty(f)},e}(f)=m\}\]
\[K_m:=(j^\infty\varphi)^{-1}(C_m)\quad\quad\dot{K}_m:=(j^\infty\varphi)^{-1}(\dot{C}_m)\]
\end{definition}
Clearly $C_\infty=\dot{C}_\infty=\cap_{m\in\NN}C_m$ and $K_\infty=\dot{K}_\infty=\cap_{m\in\NN}K_m$.

\begin{lema}
\label{deter1}
For any $x\in U$ and any $f\in \tilde{I}_x$ such that $c_{\tilde{I}_x,e}(f)=m$ we have $\mm_x^{m}\tilde{I}_x\subset\Theta_{\tilde{I}_x,e}(f)$.
\end{lema}
\begin{proof}
The $i$-th graded piece of $\tilde{I}_x/\Theta_{\tilde{I}_x,e}(f)$ by the $\mm_x$-adic filtration is the module
\[M_i:=\frac{\mm_x^i\tilde{I}_x+\Theta_{\tilde{I}_x,e}(f)}{\mm_x^{i+1}\tilde{I}_x+\Theta_{\tilde{I}_x,e}(f)}.\] 
As $m=c_{\tilde{I}_x,e}(f)=\dim_\CC(\tilde{I}_x/\Theta_{\tilde{I}_x,e}(f))=\sum_{i=0}^\infty\dim_\CC(M_i)$ we deduce that there exists $l\leq m$ such 
that $M_l=0$, which is the same that $\mm_x^l\tilde{I}_x\subset\Theta_{\tilde{I}_x,e}(f)+\mm_x^{l+1}\tilde{I}_x$. Applying Nakayama's Lemma to the module
$\mm_x^l\tilde{I}_x +\Theta_{\tilde{I}_x,e}(f)/\Theta_{\tilde{I}_x,e}(f)$ we conclude $\mm_x^l\tilde{I}_x\subset\Theta_{\tilde{I}_x,e}(f)$.
\end{proof}

\begin{lema}
\label{deter2}
For any $x\in U$ the subset 
\[A_m:=\{f\in \tilde{I}_x:\mm_x^m\tilde{I}_x\subset\Theta_{\tilde{I}_x,e}(f)\}\] 
satisfies $A_m+\mm_x^{m+2}\tilde{I}_x=A_m$.
\end{lema}
\begin{proof}
Consider $f\in A_m$ and $g\in\mm_x^{m+2}\tilde{I}_x$. Clearly 
\[\Theta_{\tilde{I}_x,e}(f)\equiv\Theta_{\tilde{I}_x,e}(f+g)\quad(\text{mod}\quad\mm_x^{m+1}\tilde{I}_x).\]
Then, as $\mm_x^m\tilde{I}_x\subset\Theta_{\tilde{I}_x,e}(f)$, we have that $\mm_x^m\tilde{I}_x\subset\Theta_{\tilde{I}_x,e}(f+g)+\mm_x^{m+1}\tilde{I}_x$;
by Nakayama's Lemma we have $\mm_x^{m}\tilde{I}_x\subset\Theta_{\tilde{I},e}(f+g)$, and hence $f+g\in A_m$. 
\end{proof}

\begin{lema}
\label{deter2.5}
The subsets $C_m$ and $K_m$ are $\lambda+m$-determined; the subsets $\dot{C}_m$ and $\dot{K}_m$ are $\lambda+m+1$-determined, where $\lambda$ is the Uniform
Artin-Rees constant.
\end{lema}
\begin{proof}
Obviously it is enough to prove the lemma for $C_m$ and $\dot{C}_m$.
To prove that $\dot{C}_m$ is $\lambda+m+1$-determined we have to check the following claim: consider $f\in\dot{C}_m$, let 
$pr_\infty(f)=x$, consider $g\in\mm^{\lambda+m+2}_x\cap\tilde{I}_x$, then $f+g\in\dot{C}_m$. By the Uniform Artin-Rees theorem $g\in\mm^{m+2}_x\tilde{I}_x$;  
Lemma~\ref{deter1} implies  
$\Theta_{\tilde{I}_x,e}(f)\supset\mm^m_x\tilde{I}_x$; then, by Lemma~\ref{deter2}, we have 
$\Theta_{\tilde{I}_x,e}(f+g)\supset\mm_x^m\tilde{I}_x$, and hence
\[\frac{\tilde{I}_x}{\Theta_{\tilde{I}_x,e}(f+g)}=\frac{\tilde{I}_x/\mm_x^m\tilde{I}_x}{\Theta_{\tilde{I}_x,e}(f+g)/\mm_x^m\tilde{I}_x}.\] 
On the other hand, as $g\in\mm^{m+2}_x\tilde{I}_x$ we have 
$X(g)\in\mm_x^{m+1}\tilde{I}_x$ for any $X\in\Theta_{\tilde{I}_x,e}$. Therefore 
$\Theta_{\tilde{I}_x,e}(f+g)/\mm_x^m\tilde{I}_x=\Theta_{\tilde{I}_x,e}(f)/\mm_x^m\tilde{I}_x$, 
and hence $c_{\tilde{I}_x,e}(f+g)=c_{\tilde{I}_x,e}(f)=m$; this shows the claim. 
As $C_m=I\setminus\cup_{i<m}\dot{C}_i$, and each $\dot{C}_i$ is $\lambda+i+1$-determined, we conclude that $C_m$ is $\lambda+m$-determined.
\end{proof}

\begin{prop}
\label{deter3}
\mbox{\quad}
\begin{enumerate}
\item The subsets $C_m$ and $K_m$ are $\lambda+m$-determined closed analytic subsets of $J^\infty(U,\tilde{I})$ and $J^\infty(U,\CC^r)$ respectively.
\item The subset $\dot{C}_m$ and $\dot{K}_m$ are $\lambda+m+1$-determined open subsets in the analytic Zariski topology of $C_m$ and $K_m$ respectively.
\end{enumerate}
\end{prop}
\begin{proof}
It is enough to prove the statements for $K_m$ and $\dot{K}_m$.
The determinacy statements are proved in Lemma~\ref{deter2.5}; for the rest we work by induction on $m$. We take as 
initial step $m=-1$, in this case everything is 
trivial. Assume that the proposition is true for any $k<m$. Then $K_{m-1}$ is closed analytic in $J^\infty(U,\CC^r)$ and 
$\dot{K}_{m-1}\subset K_{m-1}$ is an open inclusion in the analytic Zariski topology of $K_{m-1}$; as 
$K_m:=K_{m-1}\setminus\dot{K}_{m-1}$ the subset $K_m$ is closed analytic in $K_{m-1}$, and hence in $J^\infty(U,\CC^r)$.
This shows the first assertion. 

Let $\{\xi_1,...,\xi_k\}$ be a system of generators of $\Theta_{I,e}$ as a $\ocn$-module; we can assume (perhaps shrinking 
$U$) that each $\xi_i$ is defined on $U$ and that the germs $\{\xi_{1,x},...,\xi_{k,x}\}$ generate 
$\Theta_{\tilde{I}_x,e}$ for any $x\in U$. Let $h_1,...,h_l\in\ocn$ be the set of monomials in $x_1,...,x_n$ of degree 
lower or equal than $\lambda+m+2$; then $J^{\lambda+m+1}(U,\CC)_x$ is generated by $\{j^{\lambda+m+1}h_1(x),...,j^{\lambda+m+1}h_s(x)\}$ for any 
$x\in U$. Therefore, defining $\{\theta_1,...,\theta_d\}:=\{h_i\xi_j:1\leq i\leq l,\ 1\leq j\leq k\}$, the set  
$\{\theta_{1,x},...,\theta_{s,x}\}$ generates $\Theta_{\tilde{I}_x,e}/\mm^{\lambda+m+2}_x\Theta_{\tilde{I}_x,e}$ as a complex 
vector space for any $x\in U$; hence it also generates $\Theta_{\tilde{I}_x,e}/\mm^{\lambda+m+2}_x\cap\Theta_{\tilde{I}_x,e}$.

Recall that we have fixed a set of functions $g_1,...,g_r$ generating $\tilde{I}$ in $U$. For any $k\leq d$ 
define an analytic vector bundle homomorphism 
\[\alpha_k:J^{\lambda+m+2}(U,\CC^r)\to J^{\lambda+m+1}(U,\CC)\]
 by the formula
\[\alpha_k(j^{\lambda+m+2}f_{1}(x),...,j^{\lambda+m+2}f_{r}(x)):=\]
\[\sum_{i=1}^rj^{\lambda+m+1}f_{i}(x)j^{\lambda+m+1}\theta_k(g_i)(x)+j^{\lambda+m+1}\theta_k(f_i)(x)j^{\lambda+m+1}g_i(x),\]
where $x\in U$, $f_1,...,f_r\in\calO_{\CC^n,x}$. The mapping $\alpha_k$ is defined so that 
\begin{equation}
\label{veintisiete}
j^{\lambda+m+1}\theta_{k}(\varphi(f_1,...,f_r))(x)=\alpha_k(j^{\lambda+m+2}f_{1}(x),...,j^{\lambda+m+2}f_{r}(x)).
\end{equation}

Given any $f\in J^{\lambda+m+2}(U,\CC^r)$ such that $pr_{\lambda+m+2}(f)=x$, denote by $S(f)$ the subspace of $J^{\lambda+m+1}(U,\CC)_x$ spanned
by $\{\alpha_k(f):1\leq k\leq d\}$. For any $f\in J^{\lambda+m+2}(U,\CC^r)$ define $\mathrm{rk}(f):=\dim_\CC(S(f))$.
We can consider $J^{\lambda+m+1}(U,\tilde{I})_x$ as a subspace of $J^{\lambda+m+1}(U,\CC)_x$ via the canonical isomorphism
$\tilde{I}_x/\mm_x^{\lambda+m+2}\cap\tilde{I}_x\cong\tilde{I}_x+\mm_x^{\lambda+m+2}/\mm_x^{\lambda+m+2}$. Then clearly 
$S(f)\subset J^{\lambda+m+1}(U,\tilde{I})_x$. Taking into account~(\ref{veintisiete}) we deduce easily
\begin{equation}
\label{isomor}
J^{\lambda+m+1}(U,\tilde{I})_x/S(f)=\tilde{I}_x/(\Theta_{\tilde{I}_x,e}(\varphi_x(h))+\mm_x^{\lambda+m+2}\cap\tilde{I}_x)
\end{equation}
for any $h\in\calO_{\CC^n,x}^r$ such that $\pi^\infty_{\lambda+m+2}(h)=f$.  

Let $\Sigma_i$ be an $\tilde{I}$-stratum; consider the positive integer $N_i:=\mu_x(\lambda+m+1)$, where $x$ is any point in $\Sigma_i$ and $\mu_x$ is the 
function defined in Formula~\ref{hilbsam}. Define the closed analytic subset $T_i\subset J^{\lambda+m+2}(U,\CC^r)$ as follows:
\[T_i:=\{f\in J^{\lambda+m+2}(\overline{\Sigma}_i,\CC^r):\mathrm{rk}(f)<N_i-m.\}\]
 
Suppose that $h\in\dot{K}_m$, consider the $\tilde{I}$-stratum $\Sigma_i$ such that $x:=pr_\infty(h)\in\Sigma_i$. Lemma~\ref{deter1} implies 
that $\Theta_{\tilde{I}_x,e}(\varphi(h)_x)$ contains $\mm_x^{m}\tilde{I}_x$, which, by Uniform Artin-Rees, contains $\mm^{\lambda+m+2}_x\cap\tilde{I}_x$.
This, together with equality~(\ref{isomor}), implies that 
\begin{equation}
\label{veintiocho}
N_i-\rk(j^{\lambda+m+2}h(x))=c_{\tilde{I}_x,e}(\varphi(h)_x)=m.
\end{equation}
Therefore 
\begin{equation}
\label{impor}
\dot{K}_{m|\Sigma_i}\subset K_m\setminus (\pi^\infty_{\lambda+m+2})^{-1}(T_i).
\end{equation}

The union 
\[M:=\coprod_{x\in U}(\mm_x^m/\mm_x^{\lambda+m+2})^r\]
is clearly an analytic sub-bundle of $J^{\lambda+m+1}(U,\CC^r)$. Consider the vector bundle homomorphism 
$j^{\lambda+m+1}\varphi:J^{\lambda+m+1}(U,\CC^r)\to J^{\lambda+m+1}(U,\CC)$. The image of $j^{\lambda+m+1}\varphi_{|M}$
is the union
\[L:=\coprod_{x\in\Sigma_i}L_x,\]
where $L_x:=\mm_x^m\tilde{I}_x+\mm_x^{\lambda+m+2}/\mm_x^{\lambda+m+2}$.
We consider the stratification of $U$ in locally closed analytic subsets $\Sigma_{i,j}$ defined to be the minimal common refinement of the Hilbert-Samuel 
stratification with respect to $\tilde{I}$, and the minimal stratification such that the rank of $j^{\lambda+m+1}\varphi_{|M}$ restricted to each stratum is 
constant. Given any $\Sigma_{i,j}$, the restriction $L_{|\Sigma_{i,j}}$ is an analytic vector bundle over it.

Consider any stratum $\Sigma_{i,j}$. Let $Y$ be an irreducible component of the locally closed $\lambda+m+2$-determined subset 
$K_{m|\Sigma_{i,j}}\setminus (\pi^\infty_{\lambda+m+2})^{-1}(T_i)$.\\

\noindent
\textbf{Claim 1}: {\em if $Y$ and $\dot{K}_m$ have non-empty intersection then $Y\subset\dot{K}_m$}. 

Let us prove the claim: as $\rk(f)\geq N_i-m$ for any $f\in J^{\lambda+m+2}(\Sigma_i,\CC^r)\setminus T_i$, the $\lambda+m+2$ determined 
set 
\[Y':=\{f\in Y:\rk(\pi^\infty_{\lambda+m+2}(f))=N_i-m\}\]
is closed analytic in $Y$ (for being the locus where 
$\rk(pr^\infty_{\lambda+m+2}(f))$ is minimal in $Y$). By~(\ref{veintiocho}) we have $\dot{K}_m\cap Y\subset Y'$.
 
By Lemma~\ref{deter1}, if $f\in\dot{K}_m$ and $x=pr_{\lambda+m+2}(f)$ then $S(\pi^\infty_{\lambda+m+2}(f))$ contains the linear subspace $L_x$. 
Conversely, suppose that $f\in Y'$ satisfies 
\[S(\pi^\infty_{\lambda+m+2}(f))\supset L_x.\] 
Using equality~(\ref{isomor}), we deduce 
\[\Theta_{\tilde{I}_x,e}(\varphi(f)_x)+\mm_x^{\lambda+m+2}\cap\tilde{I}_x\supset\mm_x^m\tilde{I}_x.\]
By
Uniform Artin-Rees $\mm_x^{\lambda+m+2}\cap\tilde{I}_x\subset\mm_x^{m+2}\tilde{I}_x$, and therefore  
\[\mm_x^m\tilde{I}_x\subset\Theta_{\tilde{I}_x,e}(\varphi(f)_x)+\mm_x^{m+2}\tilde{I}_x.\]
Using Nakayama's Lemma we conclude that 
$\Theta_{\tilde{I}_x,e}(\varphi(f)_x)$ contains $\mm_x^m\tilde{I}_x$, which in turn contains $\mm_x^{\lambda+m+2}\cap\tilde{I}_x$. Therefore 
\[c_{\tilde{I}_x,e}(\varphi(f)_x)=\dim_\CC(\frac{\tilde{I}_x}{\Theta_{\tilde{I}_x,e}(\varphi(f)_x)+\mm_x^{\lambda+m+2}\cap\tilde{I}_x}),\]
which, by the equality~(\ref{isomor}) is equal to $N_i-\rk(f)=N_i-(N_i-m)=m$; consequently $f\in\dot{K}_m$. We have shown
\[\dot{K}_m\cap Y=\{f\in Y':S(\pi^\infty_{\lambda+m+2}(f))\supset L_x.\]

We consider the restriction 
$pr_{\lambda+m+2}:\pi^\infty_{\lambda+m+2}(Y'_{|\Sigma_{i,j}})\to\Sigma_{i,j}$, the pullback vector bundle 
$pr_{\lambda+m+2}^*(J^{\lambda+m+1}(U,\CC)_{|\Sigma_{i,j}})$, and the sub-bundle $pr_{\lambda+m+2}^*(L_{|\Sigma_{i,j}})$. Taking into account that
$S(f)$ has constant dimension when $f$ ranges in $Y'$ it is easy to check that the union 
\[S:=\coprod_{f\in\pi^\infty_{\lambda+m+2}(Y'_{|\Sigma_{i,j}})}S(f)\]
is an analytic vector sub-bundle of 
$pr_{\lambda+m+2}^*(J^{\lambda+m+1}(U,\CC)_{|\Sigma_{i,j}})$. Let $Y''$ be the 
closed analytic subset of $Y'$ formed by points $f$ such that the fibre
$S_{\pi^\infty_{\lambda+m+2}(f)}$ contains the fibre $(pr_{\lambda+m+2}^*L_{|\Sigma_{i,j}})_{\pi^\infty_{\lambda+m+2}(f)}$. As $\dot{K}_{m}\cap Y$ is the 
subset of $Y'$ formed by germs such that 
$S(\pi^\infty_{\lambda+m+2}(f))$ contains $L_x$, we obtain that $\dot{K}_{m|\Sigma_{i,j}}=Y''$. We have shown that $\dot{K}_m\cap Y$ is a closed 
analytic subset of $Y'$, and hence of $Y$.

For being $Y$ included in $K_m$, if we have $f\in Y\setminus \dot{K}_m$ and $x:=pr_\infty(f)$ then $c_{\tilde{I}_x,e}(\varphi(f)_x)>m$.
If $Y\not\subset\dot{K}_m$ then $\dot{K}_m\cap Y$ is a $\lambda+m+2$-determined proper analytic subset of the irreducible analytic subset 
$Y$. Then it is possible to find an analytic path $\gamma:D\to\pi^\infty_{\lambda+m+2}(Y)$, from a small disk $D$ to $Y$ such that
$\gamma^{-1}(\pi^\infty_{\lambda+m+2}(\dot{K}_m))=\{0\}$. 

Consider the analytic trivialisation  
\[\tau_{\lambda+m+2}:U\times (\CC\{x\}/\mm^{\lambda+m+3})^r\to J^{\lambda+m+2}(U,\CC^r)\]
defined in~(\ref{r1trivi}). There is an analytic path $\xi:D\to U$ and polynomials $h_{1,t},...,h_{r,t}$ of degree $\lambda+m+2$ in $x_1,...,x_n$, with 
coefficients depending 
analytically on $t$ such that $\gamma(t)=\tau_{\lambda+m+2}(\xi(t),(h_{1,t},...,h_{r,t}))$ for any $t\in D$. For any $t,i$ let $\tilde{h}_{i,r}$ be the unique 
function in $\ocn$ whose Taylor expansion at $\xi(t)$ equals to $h_{i,t}$; clearly $\tilde{H}_{i,t}$ is a polynomial whose the coefficients 
depend analytically on $t$.
 The $I$-unfolding 
\[F:\CC^n\times D\to\CC\]
defined by 
$F(x,t):=\varphi(\tilde{h}_{1,t},...,\tilde{h}_{r,t})(x)$ satisfies 
\[c_{\tilde{I}_O,e}(F_{|0,O})=n\quad\text{and}\quad c_{\tilde{I}_{\xi(t)},e}(F_{|t,\xi(t)})>n,\] 
which contradicts the upper-semicontinuity of the extended codimension (Corollary~\ref{uppersem}). This shows Claim 1.

Decompose each $K_m\cap(\pi^\infty_{\lambda+m+2})^{-1}(T_i)$ (for  $1\leq i\leq s$) in irreducible components; we say that an irreducible component $Z$ 
is {\em relevant} if 
$pr_{\lambda+m+2}(Z)\cap\Sigma_i\neq\emptyset$. Consider any stratum $\Sigma_{i,j}$, decompose $K_{m|\overline{\Sigma}_{i,j}}$ in irreducible components;
we say that an irreducible component is {\em excessive} if it is contained in $K_{m+1}$. Define the closed analytic subset of $T\subset K_m$ to be the union of 
all the excessive components of each of the $K_{m|\overline{\Sigma}{i,j}}$'s and all the relevant components
of each of the $K_m\cap(\pi^\infty_{\lambda+m+2})^{-1}(T_i)$'s.\\

\noindent
\textbf{Claim 2}: $T=K_{m+1}$. 

If the claim is true the induction step is complete, and hence the proposition is proved.

Let us prove the claim: consider $f\in K_m\setminus T$, let $\Sigma_{i,j}$ be the stratum such that $x=pr_\infty(f)\in\Sigma_{i,j}$.
As $f$ does not belong to $T$, it cannot be in $(\pi^\infty_{\lambda+m+2})^{-1}(T_i)$. Let $Y$ be the irreducible component of $K_{m|\overline{\Sigma}_{i,j}}$ 
to which $f$ belongs. As 
$f\not\in T$, the component $Y$ cannot be excessive, and therefore $Y\cap\dot{K}_m$ is non empty. If $Y_{|\Sigma_{i,j}}\cap\dot{K}_m\neq\emptyset$ then, by 
Claim~1 we have $Y_{|\Sigma_{i,j}}\subset\dot{K}_m$; therefore $f\in\dot{K}_m$. If $Y_{|\Sigma_{i,j}}\cap\dot{K}_m=\emptyset$ then
$Y_{|\Sigma_{i,j}}\subset K_{m+1}$. Hence $Y\cap\dot{K}_m$ is 
included in the proper analytic subset $Y_{|\partial\Sigma_{i,j}}$, where $\partial\Sigma_{i,j}$ is the difference 
$\overline{\Sigma}_{i,j}\setminus\Sigma_{i,j}$. As $Y_{|\Sigma_{i,j}}\subset K_{m+1}$ is Zarisky open in the finite-determined irreducible closed analytic 
subset $Y$, it is possible to find an analytic path $\gamma:D\to\pi^\infty_{\lambda+m+2}(Y)$, from a disk to $Y$ such that $\gamma^{-1}(\dot{K}_m)=\{0\}$. 
As we have seen previously this leads to a contradiction; we conclude $Y_{|\Sigma_{i,j}}\cap\dot{K}_m\neq\emptyset$, and, hence $f\in\dot{K}_m$. 
This shows $K_{m+1}\subset T$. 

Suppose $f\in T$. If $f$ is in a excessive component of $K_{m|\overline{\Sigma}{i,j}}$'s for certain $i,j$ then $f\in K_{m+1}$. Otherwise 
$\pi^\infty_{\lambda+m+2}(f)$ is in a relevant
component $Y$ of $K_m\cap(\pi^\infty_{\lambda+m+2})^{-1}(T_i)$ for a certain $i$. 
By~(\ref{impor}) we have $K_m\cap (\pi^\infty_{\lambda+m+2})^{-1}(T_{i|\Sigma_i})\subset K_{m+1}$; therefore 
$Y\cap\dot{K}_m$ is contained in the proper analytic subset $Y_{|\partial\Sigma_i}$. As $Y$ is irreducible we conclude like previously that $Y\subset K_{m+1}$.
This completes the proof of the claim. 
\end{proof}

The next proposition gives a codimension bound for the irreducible components of $K_1$ or $C_1$ that will be interesting for us 
(those corresponding to functions of finite extended codimension).

\begin{prop}
\label{codimbound}
If $Y$ is an irreducible component of $K_1$ (resp. $C_1$) of codimension
strictly smaller than $n$, then $Y$ is included in $K_\infty$ (resp. $C_\infty$).
\end{prop}
\begin{proof}
It is enough to work with components of $K_1$.
Let $Y\subset K_1$ be such a component. Suppose that $Y\not\subset K_\infty$; then there is $h\in Y$ such that
$c_{\tilde{I}_x,e}(j^\infty\varphi(h)(x))$ is finite, where $x=pr_\infty(h)$. The element
$h$ is the germ at $x$ of an analytic function from a neighbourhood $W$ of $x$ to $\CC^r$, with components
 $h_1,...,h_r$; the function $g:=\varphi(h)$ is a section of $\tilde{I}$ over $W$, whose germ $g_x$ at $x$ equals 
$j^\infty\varphi(h)$. Consider the analytic mapping $j^\infty h:W\to J^\infty(W,\CC^r)$; the locus $(j^\infty h)^{-1}(Y)$
consists of points in which $g$ has positive extended codimension. As $c_{\tilde{I}_x,e}(g_x)$ is finite, the 
associated sheaf $\calF$, defined in formula~(\ref{Fcal}), is skyscraper at $x$, and hence 
$\dim_x((j^\infty h)^{-1}(Y))=0$. On the other hand, the following easy statement implies that the last dimension is 
positive, giving a contradiction.

\noindent
$(\dagger)$ {\em Let $\rho:(X,x)\to (Z,z)$ be an analytic morphism of germs of complex spaces, with $(Z,z)$ smooth. Consider $(Y,z)$ an 
analytic subgerm of $(Z,z)$ such that $\codim_z(Y,Z)<\dim_x(X)$, then $\dim_x(\rho^{-1}(V))>0$}.

If $\dim_x(\rho^{-1}(z))>0$ we are done. Assume $\dim_x(\rho^{-1}(z))=0$; by the Proposition~of~page~63~of~\cite{GR}, 
shrinking $X$ and $Z$ we can assume that $\rho$ is finite, and therefore $(\rho(X),z)$ is a closed analytic subspace of 
dimension equal to $\dim_x(X)$.
As $(Z,z)$ is smooth, $\dim_z(Y\cap\rho(X))\geq\dim_z(\rho(X))-\codim_z(Y,Z)>0$. Then $\dim_x(\rho^{-1}(Y))>0$.
\end{proof} 

\section{The topological partition}

Functions of the same extended codimension do not need to have the same topological type, even if they lie in a family. For example, when $I$ defines a smooth
subvariety of dimension $d$, all the singularity types $D(d,k)$ defined in~\cite{Pe2} have zero extended codimension 
and pairwise different topological types. We need to subdivide the subvarieties $\dot{K}_m$ so that
the functions on the (relevant) resulting pieces have constant topological type. As we will show there is a canonical way
to do it.

Any representative $\phi:V\to U$ of a germ $\phi\in\calD_{\tilde{I}_x,e}$ induces by push-forward a bijection 
\begin{equation}
j^\infty\phi_*:J^\infty(V,\tilde{I})\to J^\infty(\phi(V),\tilde{I}), 
\end{equation}
defined by $j^\infty\phi_*(f_y):=(f_y\comp\phi^{-1})_{\phi(y)}$ for any $y\in V$ and $f_y\in\tilde{I}_y$.
As the $\mm_{y}$-adic filtration is transformed by push-forward into the $\mm_{\phi(y)}$-adic filtration, for any $m<\infty$, the mapping $j^\infty\phi_*$ 
descends to a bijection
\begin{equation}
j^m\phi_*:J^m(V,\tilde{I})\to J^m(\phi(V),\tilde{I}).
\end{equation}

The functions $g_1\comp\phi^{-1},...,g_r\comp\phi^{-1}$ are sections of $\tilde{I}$ defined over $\phi(V)$ (where $g_1,...,g_r$, the fixed set of generators 
for $\tilde{I}$ at $U$) . As $\tilde{I}$ is generated over $U$ by 
$g_1,...,g_r$, if we shrink $V$ enough we can assume that for any $i\leq r$ we have an expression of the form
\[g_i\comp\phi^{-1}=\sum_{j=1}^rh_{i,j}g_j,\]
where each $h_{i,j}$ is an analytic function defined on $\phi(V)$. Consequently, if $f\in\Gamma(V,\tilde{I})$ is of the form 
$f=\sum_{i=1}^rf_ig_i$, then 
\begin{equation}
\label{levad}
f\comp\phi^{-1}=\sum_{i=1}^r(f_i\comp\phi^{-1})g_i\comp\phi^{-1}=\sum_{i=1}^r(\sum_{j=1}^rf_j\comp\phi^{-1}h_{j,i})g_i.
\end{equation}
For any $m\leq\infty$ the mapping 
\begin{equation}
j^m\tilde{\phi}_*:J^k(V,\tilde{I})\to J^m(\phi(V),\tilde{I}), 
\end{equation}                                      
defined by the formula 
\[j^m\tilde{\phi}_*(j^mf_1(x),...,j^mf_r(x)):=(....,j^m(\sum_{j=1}^rf_j\comp\phi^{-1}h_{j,i})(\phi(x)),...)\]
defines a mapping which is analytic for $m<\infty$, and such that, by Equation~(\ref{levad}), satisfies 
$j^m\varphi\comp j^m\tilde{\phi}_*=j^m\phi_* j^m\varphi$. Therefore we say that 
$j^m\tilde{\phi}_*$ is an {\em analytic local lifting} of $j^m\phi_*$. 

\begin{definition}
\label{calD_I-inv}
Let $T\subset J^\infty(U,\tilde{I})$ be a (locally) closed analytic subset. We say that $T$ is $\calD_{\tilde{I}_e}$-\emph{invariant} if for any analytic 
diffeomorphism $\phi:V\to U$ preserving $\tilde{I}$ and any $x\in V$, we have that
$\phi_{t*}(T_x)$ is equal to $T_{\phi_t(x)}$.
\end{definition}

Clearly the set $C_m$ is $\calD_{\tilde{I},e}$-invariant for any $m\leq\infty$. 

Now we recall them briefly the results of~\cite{Bo2}:

\begin{definition}
Two germs $f:(\CC^n,x)\to\CC$ and $g:(\CC^n,y)\to\CC$ are called {\em topologically equivalent} if there exits germs of homeomorphisms 
$\phi:(\CC^n,x)\to (\CC^n,y)$ and $\alpha:(\CC,f(x))\to (\CC,g(x))$ such that $\alpha\comp f=g\comp\phi$.
\end{definition}

A closed subset $A\subset J^\infty(U,\tilde{I})$ is called {\em residual} if for any positive integer $c$ there is a closed analytic subset $T'$ 
containing $A$, and such that all its irreducible components are of codimension at least $c$.  

Let $T$ be an irreducible locally closed subset of $J^{\infty}(U,\tilde{I})$, due to Proposition~10 of~\cite{Bo2}, there exists a unique closed subset 
$\Gamma\subset T$ with the following properties: 
 \begin{enumerate}[(i)]
\item We have a decomposition $\Gamma=\Gamma^{(a)}\cup\Gamma^{(r)}$
where $\Gamma^{(a)}$ is a closed analytic subset of $T$, and $\Gamma^{(r)}$ is a residual closed subset.
\item Any $f,g\in T$ in the same path-connected component of $T\setminus\Gamma$ are topologically equivalent. 
\item The subset $\Gamma$ is minimal among the subsets of $T$ satisfying Properties~(i)~and~(ii).
\end{enumerate}

The subset $\Gamma$ is called the {\em topological discriminant} of $T$. The decomposition $\Gamma=\Gamma^{(a)}\cup\Gamma^{(r)}$ is unique provided that 
$\Gamma^{(r)}$ is chosen minimal; we say that $\Gamma^{(a)}$ and $\Gamma^{(r)}$ are respectively the {\em analytic and residual parts} of $\Gamma$. 

In~\cite{Bo2} it was proved that if $V$ is an open subset of $U$ then $\Gamma_{|V}$ 
is the topological discriminant of $T_{|V}$.  

\begin{lema}
\label{invariant}
If $T$ is $\calD_{\tilde{I},e}$-invariant then $\Gamma$ is also 
$\calD_{\tilde{I},e}$-invariant. 
\end{lema}
\begin{proof}
Let $\phi:V\to\phi(V)$ be an analytic diffeomorphism preserving $\tilde{I}$. It is 
enough to show
\begin{equation}
\label{enough}
\phi_*(\Gamma_{|V})\supset\Gamma_{|\phi(V)}
\end{equation}
(for the 
opposite inclusion we consider the analogous statement for $\phi^{-1}_*$ and apply 
$\phi_*$). The set $\Gamma':=\phi_*(\Gamma_{|V})\cap \Gamma_{|\phi(V)}$ is a closed 
subset which is union of the analytic subset 
$\Gamma^{(a)}_{\phi(V)}\cap\phi_*(\Gamma^{(a)}_{|V})$ and the residual subset 
\[[\Gamma^{(r)}_{|\phi(V)}\cap\phi_*(\Gamma_{|V})]\cup [\Gamma_{|\phi(V)}\cap\phi_*(\Gamma^{(r)}_{|V})].\] 
If  any two germs that can be connected by a continuous path in 
$T_{|\phi(V)}\setminus\Gamma'$ have the same topological type then $\Gamma'$ must contain the 
topological discriminant $\Gamma_{|\phi(V)}$ of $T_{|\phi(V)}$, and therefore 
inclusion~(\ref{enough}) holds.

Let $\gamma:[0,1]\to T_{|\phi(V)}\setminus\Gamma'$ be a continuous path. The interval
$[0,1]$ is the union of the open subsets 
$\gamma^{-1}(T_{|\phi(V)}\setminus\Gamma_{|\phi(V)})$ and 
$\gamma^{-1}(T_{|\phi(V)}\setminus\phi_*(\Gamma_{|V})$. As over each of these subsets
the topological type clearly remains constant, then it also does along $\gamma$.
\end{proof}

Shrink $U$ so that its closure is contained in an open subset where $\Theta_{\tilde{I},e}$ is generated by global sections, then the Main Theorem 
of~\cite{Bo2}, applied to $T=J^{\infty}(U,\tilde{I})$, gives:

\begin{theo} 
\label{Varchenko} 
There exist a unique filtration (which we call the {\em filtration by successive discriminants}) 
\[J^{\infty}(U,\tilde{I})=A_0\supset A_1\supset ...\supset A_i\supset ...\]
by closed analytic subsets, and a residual subset $\Gamma^{(r)}$ (called the {\em cumulative residual topological discriminant} of 
$J^{\infty}(U,\tilde{I})$),  with the following properties: 
\begin{enumerate}
\item We have $\cap_{i\geq 0}A_i\subset\Gamma^{(r)}$. 
\item For any $i\geq 0$ the set $A_{i+1}\cup (\Gamma^{(r)}\cap A_i)$ is the topological discriminant of $A_i$.
\item Any irreducible component of $A_i$ has codimension at least $i$.  
\item If $T$ is $\calD_{\tilde{I},e}$-invariant then $A_i$ is $\calD_{\tilde{I},e}$-invariant for any $i\geq 0$. The set $\Gamma^{(r)}$ is contained in a 
residual subset which is an intersection of $\calD_{\tilde{I},e}$-invariant closed analytic subsets of $C$.
\end{enumerate}
\end{theo}

We need to make a remark concerning the level of generality of this paper in comparison to~\cite{Bo2}. Instead of working with analytic subsets of 
$J^\infty(U,\tilde{I})$, in~\cite{Bo2} our attention was restricted to closed analytic subsets of $J^\infty(\Sigma_i,\tilde{I})$, where $\Sigma_i$ is a 
certain $\tilde{I}$-stratum of $U$. Once that the existence of decompositions in irreducible components has been proved in our setting, all the arguments 
of~\cite{Bo2} can be translated with minimal changes to prove the statements given here. Another 
difference is in Property~$4$: here we state the $\calD_{\tilde{I},e}$-invariance of
certain subsets, and in~\cite{Bo2} we state flow-invariance, which is a weaker 
property. The stronger property holds because here we work over $\CC$, and 
in~\cite{Bo2} we considered also the real case; actually Lemma~\ref{invariant} 
provides the needed additional arguments to those given in~\cite{Bo2}.       

Because of the fourth assertion of Theorem~\ref{Varchenko} there exists a $\calD_{\tilde{I},e}$-invariant closed analytic subset $\Delta_1$ 
containing $\Gamma^{(r)}$, which has all of its irreducible components of codimension at least $n+1$. We can suppose $\Gamma^{(r)}=\cap_{k\in\NN}\Gamma_k$, 
being each $\Gamma_k$ closed analytic with the properties  $\Delta_1\supset\Gamma_k$ and $\Gamma_k\supset\Gamma_{k+1}$. 
For any $i\leq n$ let $\{A_{i,j}\}_{j\in J_i}$ be the set of irreducible components of $A_i\setminus A_{i+1}$ of codimension smaller or equal 
than $n$. The subset $A_{i,j}\setminus\Gamma^{(r)}$ is path-connected for any $i,j$, as it is the union $\cup_{k\in\NN}A_{i,j}\setminus\Gamma_k$, forming the 
$A_{i,j}\setminus\Gamma_k$'s an increasing sequence of irreducible (and hence connected) locally closed subsets. We conclude that all the germs in 
$A_{i,j}\setminus\Gamma^{(r)}$ have the same topological type.

We say that two components $A_{i,j}$ and $A_{i,j'}$ are equivalent if their respective generic 
germs have the same topological type. 
Non-equivalent $A_{i,j}$ and $A_{i,j'}$ are disjoint. If this were not the case let $Y$ be an irreducible component of their intersection. As $\Gamma^{(r)}$ 
is contained in closed subsets of arbitrarily large codimension we deduce that $Y$ is not contained in $\Gamma^{(r)}$. Therefore the topological type of a 
generic germ in $Y$ coincides with the generic topological type in $A_{i,j}$ and $A_{i,j'}$. This contradicts the non-equivalence of $A_{i,j}$ and $A_{i,j'}$.
Let $L_i$ be the set of equivalence classes
and for any $l\in L_i$ define $B_{i,l}$ as the union of all the subsets of the class $l$. Consider the decomposition
\[\bigcup_{j\in J_i}A_{i,j}=\coprod_{l\in L_i}B_{i,l}.\]
Clearly the sets $B_{i,l}$ are closed analytic in $A_i\setminus A_{i+1}$ and the elements of 
$B_{i,l}\setminus\Gamma^{(r)}$ and $B_{i,l}\setminus\Delta_1$ are germs which 
have pairwise the same topological type.

Define $\Delta_{2,i}$ to be the union of all the irreducible components of $A_i$ of codimension at least $n+1$. Define 
$\Delta_2:=\cup_{j\in\NN}\Delta_{2,j}$. As $\Delta_2$ contains $A_{n+1}$ the union is easily seen to be locally finite, and hence $\Delta_2$ is an analytic 
closed subset. 

We have a locally finite partition
\begin{equation}
\label{partition}
J^\infty(U,\tilde{I}):=\big[\coprod_{i\leq n,l\in L_i}B_{i,l}\setminus (\Gamma^{(r)}\cup\Delta_2)\big]\coprod (\Gamma^{(r)}\cup\Delta_2)
\end{equation}
in disjoint subsets such that for any $i\leq n$ and $l\in L_i$ the set $B_{i,l}$ is a $\calD_{\tilde{I},e}$-invariant locally closed analytic subset 
in $J^\infty(U,\tilde{I})$ such that any two germs in $B_{i,l}\setminus\Gamma^{(r)}$ have the same topological type.
All the subsets of the partition are canonically defined.

\begin{definition}
\label{toppar}
The partition introduced above is called the {\em Topological Partition of} $J^\infty(U,\tilde{I})$ up to codimension $n$, and is canonically defined.
\end{definition}

\section{Whitney stratifications in $J^\infty(U,\tilde{I})$}

Let $C$ be a locally closed analytic subset of $J^{\infty}(U,\tilde{I})$. An {\em stratification of} $C$ is a partition of $C$ in a locally finite family
$\{X_j\}_{j\in J}$ of disjoint smooth irreducible locally closed analytic subsets of $J^{\infty}(U,\tilde{I})$. Given two smooth irreducible 
locally closed analytic subsets $X$ and $Y$ of $J^{\infty}(U,\tilde{I})$, we say that $X$ is {\em Whitney-regular} over $Y$ if for any $k$ such that both 
$X$ and $Y$ are $k$-determined, and for any open subset $V\subset U$
and any system of generators $\calH$ of $\tilde{I}_{|V}$, the subset $(j^k\varphi_\calH)^{-1}(X)$ is Whitney-regular over 
$(j^k\varphi_\calH)^{-1}(Y)$ (this makes sense as $(j^k\varphi_\calH)^{-1}(X)$ 
and $(j^k\varphi_\calH)^{-1}(Y)$ are submanifolds of $J^k(U,\CC^r)$). 

It is easy to show that given two smooth locally closed analytic subsets $X, Y\subset J^k(U,\CC^r)$, a positive integer $m\geq k$, and a point 
$y\in (\pi^m_k)^{-1}(Y)$ then $(\pi^m_k)^{-1}(X)$ is Whitney-regular over $(\pi^m_k)^{-1}(Y)$ at a point $y$ if and only if $X$ is Whitney-regular over 
$Y$ at $\pi^m_k(y)$. Therefore it is enough to check the Whitney regularity condition for a particular $k$. An argument similar to the proof of 
Lemma~\ref{indepcarta} shows that to prove Whitney regularity at a point it is enough to check it at a single chart containing the point.

A {\em Whitney} stratification of a locally closed analytic subset of $J^{\infty}(U,\tilde{I})$ is an stratification of it such that any stratum is Whitney
regular over any other stratum.

We will make use of the following fact (see~\cite{W}): let $X$ and $Y$ be two irreducible locally closed analytic subsets of a complex manifold, such 
that $\dim(X)>\dim(Y)$. Denote by $X_{sm}$ and $Y_{sm}$ the set of smooth points of $X$ and $Y$. There exists a unique minimal proper closed analytic subset 
$W(X,Y)$ of $Y$ containing $\mathrm{Sing}(Y)$ such that $X_{sm}$ is Whitney regular over $Y_{sm}\setminus W(X,Y)$. Moreover the set of points $y\in Y_{sm}$ 
such that $X$ is not Whitney regular over $Y$ at $y$ is dense in $W(X,Y)\cap Y_{sm}$.

\begin{lema}
\label{whitneylemma}
Let $X$ and $Y$ be two $k$-determined irreducible locally closed analytic subsets of $J^\infty(U,\tilde{I})$ such that 
$\codim(X,J^\infty(U,\tilde{I}))<\codim(Y,J^\infty(U,\tilde{I}))$. There exists a unique minimal $k+\lambda$-determined proper closed analytic subset $W(X,Y)$ 
of $Y$ such that $X_{sm}$ is Whitney regular over $Y_{sm}\setminus W(X,Y)$. Moreover the set of points $y\in Y_{sm}$ such that $X$ is not Whitney regular over 
$Y$ at $y$ is dense in $W(X,Y)\cap Y_{sm}$. In addition, if $X$ and $Y$ are both of them $\calD_{\tilde{I},e}$-invariant, 
then $W(X,Y)$, $\mathrm{Sing}(X)$ and $\mathrm{Sing}(Y)$ are $\calD_{\tilde{I},e}$-invariant.
\end{lema}

\begin{proof}
Let $k$ be such that both $X$ and $Y$ are $k$-determined. The sets $X':=(j^k\varphi)^{-1}(\pi^\infty_k(X))$ and 
$Y':=(j^k\varphi)^{-1}(\pi^\infty_k(Y))$ are locally closed analytic 
subsets of $J^k(U,\CC^r)$ such that $\dim(X')>\dim(Y')$; consider the set $W(X',Y')$. If we prove that for a certain $m\geq k$ the set 
$W:=(\pi^m_k)^{-1}(W(X',Y'))$ is $j^m\varphi$-saturated, then the set $W(X,Y):=(\pi^\infty_m)^{-1}(j^m\varphi(W))$ clearly satisfies all the desired 
properties, with the exception of the one concerning $\calD_{\tilde{I},e}$-invariantness, which will be proved later. 

Let us prove the $j^m\varphi$-saturation for 
$m=k+\lambda$, where $\lambda$ is the uniform Artin-Rees constant. We have to check 
\begin{equation}
\label{wish2}
W_x+\mathrm{ker}(j^m\varphi_x)\subset W_x
\end{equation}
for any $x\in U$. Because of inclusion~(\ref{conclusion0}) it is enough to show that $W$ is invariant by the group action~(\ref{accion}). The subset $W$ is
clearly equal to $W((\pi^m_k)^{-1}(X'),(\pi^m_k)^{-1}(Y'))$; therefore, the subset $Z$ of $(\pi^m_k)^{-1}(Y')$ formed by the points where $(\pi^m_k)^{-1}(Y')$
is singular and the points at which 
$(\pi^m_k)^{-1}(X')$ is not Whitney-regular over $(\pi^m_k)^{-1}(Y')$ is dense in $W$. As $(\pi^m_k)^{-1}(X')$ and $(\pi^m_k)^{-1}(Y')$ are invariant by the 
action~(\ref{accion}) and Whitney regularity is preserved by diffeomorphisms, then $Z$ is also invariant by the action. This implies in turn the 
invariance of $W$.  

It only remains to be proved the $\calD_{\tilde{I},e}$-invariance of $W(X,Y)$.

Consider an open subset $V$ of $U$ and a set $\calH=\{h_1,...,h_s\}$ which generate $\tilde{I}_{|V}$. As Whitney regularity can be checked with respect to any 
chart, we have that that $W(X,Y)\cap J^\infty(V,\tilde{I})=(\pi^\infty_{m})^{-1}(j^m\varphi_\calH(W'))$, where 
$W'=(\pi^{k+\lambda}_k)^{-1}(W(X'_\calH,Y'_\calH))$, for $X'_\calH:=(j^k\varphi_\calH)^{-1}(\pi^\infty_k(X))$ and 
$Y'_\calH:=(j^k\varphi_\calH)^{-1}(\pi^\infty_k(Y))$. 

Suppose that $X$ and $Y$ are $\calD_{\tilde{I},e}$-invariant. Consider an analytic 
diffeomorphism $\phi:V\to U$ preserving $\tilde{I}$. We consider the associated 
push-forward mapping 
\begin{equation}
j^\infty\phi_*:J^\infty(V,\tilde{I})\to J^\infty(\phi(V),\tilde{I}). 
\end{equation}
As $\tilde{I}$ is generated over $U$ by $g_1,...,g_r$, the functions $g_1\comp\phi^{-1},...,g_r\comp\phi^{-1}$ generate $\tilde{I}$ over $\phi(V)$. 
Define $\alpha:\calO_{\phi(V)}^r\to\tilde{I}_{|\phi(V)}$ by the formula $\alpha(f_1,...,f_r):=\sum_{i=1}^rf_i(g_i\comp\phi^{-1})$. For any $k\leq\infty$ we let
\[j^k\alpha:J^k(\phi(V),\CC^r)\to J^k(\phi(V),\tilde{I})\]
be the associated mapping of $k$-jets. For any $k\leq\infty$ the mapping 
\begin{equation}
j^k\tilde{\phi}_*:J^k(V,\CC^r)\to J^k(\phi(V),\CC^r), 
\end{equation}                                      
defined by the formula 
\[j^k\tilde{\phi}_*(j^kf_1(x),...,j^kf_r(x)):=(j^k(f_1\comp\phi^{-1})(\phi(x)),...,j^k(f_r\comp\phi^{-1})(\phi(x)))\]
satisfies $j^k\alpha\comp j^k\tilde{\phi}_*=j^k\phi_*\comp j^k\varphi$, and, if $k<\infty$, defines an analytic isomorphism. Letting $k=m$ and taking into 
account the fact that the definition of the sets $\mathrm{Sing}(X_{|\phi(V)})$, $\mathrm{Sing}(Y_{|\phi(V)})$ and $W(X_{|\phi(V)},Y_{|\phi(V)})$ does not 
depend on the chosen set of generators of $\tilde{I}_{\phi(V)}$, we obtain easily that  invariance of $\mathrm{Sing}(X)$, 
$\mathrm{Sing}(Y)$ and $W(X,Y)$ by $\phi_*$.
\end{proof}

\begin{theo}
\label{Whitneystr}   
Let $X$ be a closed analytic subset of $J^\infty(U,\tilde{I})$, consider a locally finite partition 
\begin{equation}
\label{parinit}
X:=\coprod_{j\in J}X_j
\end{equation}
by closed analytic subsets of $J^\infty(U,\tilde{I})$. There is a canonical Whitney stratification of $X$ such that each stratum
is a locally closed analytic subset contained in one of the sets $X_i$. Moreover, if each of the subsets $X_j$ is 
$\calD_{\tilde{I},e}$-invariant, then the strata are $\calD_{\tilde{I},e}$-invariant.
\end{theo}

\begin{proof}
The construction is inductive: suppose that for a certain $N$ we have constructed a locally finite partition in disjoint locally closed analytic subsets
\begin{equation}
\label{parNth}
X:=(\coprod_{i\in I_N}Z_i)\coprod (\coprod_{i\in I'_N}Y_i)
\end{equation}
refining the original partition~(\ref{parinit}) and with the following properties: for any $i\in I_N$ the set $Z_i$ is smooth and irreducible of codimension at
most $N$, for any $i_1,i_2\in I_N$ the stratum 
$Z_{i_1}$ is Whitney regular over $Z_{i_2}$, for any $i\in I_N$ the set $Y_i$ is irreducible, and the union $\cup_{i\in I'_N}Y_i$ is closed analytic with 
all the irreducible components of codimension at least $N+1$. Moreover if each of the subsets $X_j$ is $\calD_{\tilde{I},e}$-invariant, then the $Z_i$'s and 
$Y_i$'s are also $\calD_{\tilde{I},e}$-invariant.

Let $L\subset I'_N$ be the set of indices parametrising components of codimension precisely 
$N+1$. For any $i\in L$ we define 
\[A_i:=\mathrm{Sing}(\overline{Y_i})\bigcup (\overline{Y}_i\setminus Y_i)\bigcup [\bigcup_{j\neq i}\overline{Y_j}\cap\overline{Y}_i]\bigcup [\bigcup_{j\in I_N}W(Z_j,\overline{Y}_i)].\]
To check that it is a closed analytic subset of $\overline{Y}_i$ we only have to show that the two possibly infinite unions involved in its definition are 
locally finite; this follows easily from the local finiteness of the family~(\ref{parNth}). 
We define $I_{N+1}:=I_N\cup L$ and $Z_i:=Y_i\setminus A_i$ for any $i\in L$. Let $\{Y_i\}_{i\in I'_{N+1}}$ be the set of disjoint subsets consisting of the 
union of $\{Y_i\}_{i\in I'_N\setminus L}$ and the sets of irreducible components of $A_i\cap Y_i$ for any $i\in L$. The partition 
\begin{equation}
\label{parN+1th}
X:=(\coprod_{i\in I_{N+1}}Z_i)\coprod (\coprod_{i\in I'_{N+1}}Y_i)
\end{equation}
has the properties of the partition~(\ref{parNth}) substituting $N$ by $N+1$. Repeating this process infinitely we obtain the desired Whitney stratification.
\end{proof}
 
\section{Transversality in Generalized Jet-Spaces}

\begin{definition}
Let $W$ be an open subset of $\CC^m$. A mapping $\alpha:W\to J^\infty(U,\CC^r)$ is {\em analytic} if the composition $\pi^\infty_k\comp\alpha$ is analytic for
any positive integer $k$. A mapping $\alpha:W\to J^\infty(U,\tilde{I})$ is {\em analytic} if there exists a neighbourhood $V$ around each point 
$x\in W$, and an analytic mapping $\tilde{\alpha}:V\to J^\infty(U,\CC^r)$ such that $\varphi\comp\tilde{\alpha}=\alpha_{|V}$. We say that 
$\tilde{\alpha}$ is a {\em local lifting} of $\alpha$ at $x$.
\end{definition}

\begin{definition}
\label{transdef}
Let $C\subset J^\infty(U,\tilde{I})$ be a closed analytic subset endowed with a stratification $\calX=\{X_j\}_{j\in J}$ by smooth irreducible locally closed 
subsets. Consider an analytic mapping $\alpha:W\to J^\infty(U,\tilde{I})$. We say that $\alpha$ is transversal to the stratification at a point $x\in W$, and 
we denote it by $\alpha\pitchfork_x\calX$ if either $\alpha(x)\not\in C$, or, when $\alpha(x)\in X_j$ for a certain $j\in J$, there exists a local lifting 
$\tilde{\alpha}:V\to J^\infty(U,\CC^r)$ around $x$, and a positive integer $m$ such that $X_j$ is $m$-determined and
\begin{equation}
\label{transcond}
d(\pi^\infty_m\comp\tilde{\alpha})_x(T_xV)+T_{\pi^\infty_m\comp\tilde{\alpha}(x)}\pi^\infty_m(\varphi^{-1}(X_i))=T_{\pi^\infty_m\comp\tilde{\alpha}(x)}J^m(U,\CC^r).
\end{equation}
\end{definition}

It is easy to check that Condition~{(\ref{transcond}) holds for a certain $m\geq k$ if and only if it holds for any $m\geq k$. As usually, transversality to a 
Whitney stratification is an open condition: 

\begin{lema}
\label{opentrans}
Suppose that the stratification $\calX$ considered in the last definition is a Whitney stratification. If $\alpha\pitchfork_x\calX$ then there is a 
neighbourhood $V$ of $x$ in $W$ such that $\alpha\pitchfork_y\calX$ for any $y\in V$.
\end{lema}
\begin{proof}
As $\calX=\{X_j\}_{j\in J}$ is locally finite there is a neighbourhood $\Omega$ of $\alpha(x)$ in $J^\infty(U,\tilde{I})$ meeting only finitely many strata 
$\{X_j\}_{j\in J'}$. Let $k$ be a positive integer such that all these strata are $k$-determined. Consider a local lifting 
$\tilde{\alpha}:V_1\to J^\infty(U,\CC^r)$ with $V_1\subset\alpha^{-1}(\Omega)$ and such that condition~(\ref{transcond}) is satisfied for $m=k$. 
This precisely means that the mapping $\pi^\infty_m\comp\tilde{\alpha}:V_1\to J^k(U,\CC^r)$ is transversal to the complex manifold 
$\pi^\infty_k(\varphi^{-1}(X_i))$ at $x$. For being $\calX$ a Whitney stratification the partition $\calX':=\{\pi^\infty_k(\varphi^{-1}(X_j))\}_{j\in J'}$ 
forms a Whitney stratification in the usual sense, and therefore there exists an open neighbourhood $V$ of $x$ in $V_1$ such that 
$\pi^\infty_m\comp\tilde{\alpha}$ is transversal to any stratum of $\calX'$. Clearly $\alpha\pitchfork_y\calX$ for any $y\in V$. 
\end{proof}

\begin{theo}[Generalized Parametric Transversality]
\label{transpar}
Let $M$ and $S$ be complex manifolds. Let $\varphi:M\times S\to J^\infty(U,\tilde{I})$ be a analytic mapping.  Consider a closed analytic subset $C$ of 
$J^\infty(U,\tilde{I})$ endowed with a Whitney stratification $\calX$. Suppose that $\varphi$ is transversal to $\calX$. Denote by $\Delta\subset S$ the 
set of points where $\varphi_{|s}:=\varphi_{|M\times\{s\}}$ is not transversal to $\calX$. Then $S\setminus\Delta$ is dense in $S$. Moreover, if there exists 
a compact subset $K\subset M$ such that $\varphi_s$ is transversal to $\calX$ at any point of $(M\setminus K)\times S$, then $\Delta$ is a proper closed 
analytic subset of $S$.
\end{theo}
\begin{proof}
Using local liftings it is easy to show that for being $\varphi$ analytic, the 
subset $C':=\varphi^{-1}(C)$ is closed analytic in $M\times S$, and for being 
$\varphi$ transversal to $\calX$ the stratification 
$\calY:=\{Y_j\}_{j\in J}$ (where $Y_j=\varphi^{-1}(X_j)$) is a Whitney stratification
of $\varphi^{-1}(C)$.  

Let $\pi:M\times S\to S$ be the projection to the second factor. Take 
$(x,s)\in M\times S$; if $(x,s)\in C'$ let $Y_j$ be the stratum of $\calY$ to which 
it belongs. An straightforward argument shows that the mapping $\varphi_{|s}$ is 
transversal to $\calX$ at $(x,s)$ if and only either $(x,s)\notin C'$ or $\pi_{|Y_j}$
is a submersion in $(x,s)$. For any $j\in J$ we let $Z_j$ be the set of points where
$\pi_{|Y_j}$ fails to be a submersion. By Sard's Theorem $\pi(Z_j)$ is a set of 
measure $0$ in $S$. The set $\Delta$ is equal to the union 
$\cup_{j\in J}\pi(Z_j)$, and has measure $0$ for being $J$ denumerable. This shows 
that $S\setminus\Delta$ is dense.

Suppose that there exists a compact $K$ with the property stated in the theorem. We
claim that the union $Z=\cup_{j\in J}Z_j$ is a closed analytic subset of $M\times S$.
Then, as $Z\subset K\times S$ the restriction $\pi_{|Z}$ is proper and 
$\Delta=\pi(Z)$ is a closed analytic subset. It only remains to show the claim. Using
the fact that $\calY$ is a {\em Whitney} stratification it is easy to show that $Z$ 
is closed. Therefore $Z=\cup_{j\in J}\overline{Z}_j$. Due to the locally finiteness 
of $\{Y_j\}_{j\in J}$, if $\overline{Z}_j$ is a closed analytic subset for any 
$j\in J$ the claim is true. For any $j\in J$ we consider the closure of the stratum 
$\overline{Y}_j$. Consider $(x,s)\in \overline{Y}_j$; choose local coordinates 
$(y_1,...,y_k)$ of $S$ around $s$, let $(\pi_1,...,\pi_k)$ be the components of 
$\pi$ with respect to it; let $f_1,...,f_N$ be a set of analytic equations in a 
neighbourhood $V$ of $(x,s)$ in $M\times S$ defining $\overline{Y}_j\cap V$. If the 
codimension of $\overline{Y}_i$ is $c$ then the rank of the set 
$\{df_1,...,df_N,d\pi_1,...,d\pi_k\}$ of $1$-forms defined over $V$ is at most $c+k$.
Define 
\[Z'_j:=\{z\in\overline{Y}_j\cap V:\mathrm{rank}\{df_1(z),...,df_N(z),d\pi_1(z),...,d\pi_k(z)\}<c+k\}.\]
Clearly $Z'_j$ is a closed analytic subset of $V$ such that $Z'_j\cap Y_j=Z_j\cap V$.
Therefore $\overline{Z}_j\cap V$ is analytic in $V$.
\end{proof}

Next we show how the versality of an unfolding implies the transversality of its associated jet extension with respect to any 
$\calD_{\tilde{I},e}$-invariant subset. 

Let $F:(\CC^n\times\CC^s,(O,O))\to\CC$ be an $s$-parametric $I$-unfolding. Due to Lemma~\ref{igualdadbasica} there are open neighbourhoods $V$ and $W$ of the 
origin in $\CC^n$ and $\CC^s$ respectively, such that $V\subset U$ and we can write $F=\sum_{i=1}^rg_iF_i$, being $\{g_1,...,g_r\}$ our fixed set of 
generators of $\tilde{I}$, and each $F_i$ an analytic function on $V\times W$. Then the jet extension
\begin{equation}
\label{jetextension}
\rho_F:V\times W\to J^\infty(U,\tilde{I})
\end{equation}
defined by $\rho_F(x,s)=j^\infty F_{|s}(x)$ is an analytic mapping, as it admits the analytic lifting 
\begin{equation}
\label{lifting}
\tilde{\rho}_F:V\times W\to J^\infty(U,\CC^r)
\end{equation}
defined by $\tilde{\rho}_F(x,s)=(j^\infty F_{1|s}(x),...,j^\infty F_{r|s}(x))$.

\begin{prop}
\label{versalimplytrans} 
Consider $F:(\CC^n\times\CC^s,(O,O))\to\CC$ as above. If $F$ is versal at $(x,s)\in V\times W$, then $\rho_F$ is transversal at $(x,s)$ to any $\calD_{\tilde{I},e}$-invariant
smooth locally closed subset of $J^\infty(U,\tilde{I})$.
\end{prop}
\begin{proof}
Let $C$ be a $k$-determined $\calD_{\tilde{I},e}$-invariant locally smooth closed subset. Suppose that $F$ is versal at $(x,s)$. If $\rho_F(x,s)\notin C$ 
there is nothing to prove. Suppose $\rho_F(x,s)\in C$. The set $A:=(j^k\varphi)^{-1}(\pi^\infty_k(C))$ is a smooth analytic subset of 
$J^k(U,\CC^r)$.

To get shorter formulas we denote the function $F_{|s}$ by $f$, the point $\pi^\infty_k(\rho_F(x,s))\in J^k(U,\tilde{I})$ by $p$, and the point 
$\pi^\infty_k(\tilde{\rho}_F(x,s))\in J^k(U,\CC^r)$ by $q$. Notice that $\rho_F(x,s)=f_x=j^\infty f(x)$, and $p=j^kf(x)=\pi^\infty_k(f_x)$.     

We have to show that
\begin{equation}
\label{want}
d(\pi^\infty_k\comp\tilde{\rho}_F)_{(x,s)}(T_{(x,s)}V\times W)+T_qA=T_qJ^k(U,\CC^r).
\end{equation}

As $J^k(U,\CC^r)$ is a trivial 
vector bundle, its tangent space at $q$ splits as 
the direct sum
\begin{equation}
\label{directsum}
T_{q}J^k(U,\CC^r)=T_xU\oplus E,
\end{equation}
where $E$ denotes $(\calO_{U,x}/\mm_x^{k+1})^r$, the fibre of the vector bundle $J^k(U,\CC^r)$ over $x$. 
Then the differential $d(pr_k)_q$ is the projection homomorphism to
the first summand. Denote by $\beta$ the 
projection homomorphisms to the second summand. We 
have an epimorphism $j^k\varphi_x:E\to\tilde{I_x}/(\tilde{I}_x\cap\mm^{k+1}_x)$. 
As $A$ is $j^k\varphi$-saturated we have
\begin{equation}
\label{primertrozo}
\{0\}\oplus\mathrm{ker}(j^k\varphi_x)\subset T_qA.
\end{equation}
 
Consider any germ of vector field $X\in\Theta_{\tilde{I}_x,e}$. Let 
$\phi:V_x\times (-\delta,\delta)\to U$ be a flow obtained integrating a 
representative of $X$ in  a neighbourhood $V_x$ of $x$ in $U$. As $f_x$ belongs to the $\calD_{\tilde{I},e}$-invariant subset $C$ we have that
$\phi_{t*}f_x$ belongs to $C$ for any $t\in (-\delta,\delta)$. Therefore the mapping
\[\gamma:(-\delta,\delta)\to J^k(U,\CC)\]
 defined by $\gamma(t):=\pi^\infty_k(\phi_{t*}f_x)=j^k(\phi_{t*}f)(\phi_t(x))$ is a differentiable curve whose image 
lies in $\pi^\infty_k(C)\subset J^k(U,\tilde{I})$ and such that $\gamma(0)=p$. Let $\Sigma_i$ be the
$\tilde{I}$-stratum to which $x$ belongs. Clearly if $\tau:=pr_k\comp\gamma$ then $\tau(t)=\phi_t(x)$. 
As any diffeomorphism preserving 
$\tilde{I}$ leaves invariant the $\tilde{I}$-strata, we have that the integral curve
$\tau$ maps $(-\delta,\delta)$ into $\Sigma_i$, and that the image
of $\gamma$ is contained in $J^k(\Sigma_i,\tilde{I})$. By Remark~\ref{restricciones}, the mapping
$j^k\varphi$ restricts to an epimorphism of analytic vector bundles
\[j^k\varphi_{|\Sigma_i}:J^k(\Sigma_i,\CC^r)\to J^k(\Sigma_i,\tilde{I}).\]
As $j^k\varphi_{|\Sigma_i}(q)=p$ there exists a differentiable lifting 
\begin{equation}
\label{elevacion}
\tilde{\gamma}:(-\delta,\delta)\to J^k(\Sigma_i,\CC^r)
\end{equation}
such that $\gamma=j^k\varphi\comp\tilde{\gamma}$ and $\tilde{\gamma}(0)=q$. As any lifting of $\gamma$ has its image contained in
$A$, the tangent vector $\tilde{\gamma}'(0)$ belongs to $T_qA$. Clearly $d(pr_k)_q(\tilde{\gamma}'(0))=\tau'(0)=X(x)$. 

Consider the pullback
vector bundle $\tau^*J^k(\Sigma_i,\tilde{I})\to (-\delta,\delta)$. Its fibre over
$0$ is $\tilde{I}_x/\tilde{I}_x\cap\mm_x^{k+1}$. Denote by 
$\eta:\tau^*J^k(\Sigma_i,\tilde{I})\to J^k(\Sigma_i,\tilde{I})$ the  
mapping defined by $\eta(h,t):=h$ for any $t\in (-\delta,\delta)$ and 
$h\in (\tau^*J^k(\Sigma_i,\tilde{I}))_t=J^k(\Sigma_i,\tilde{I})_{\tau(t)}$. We consider the trivialisation  
\begin{equation}
\label{latrivia}
\psi:(-\delta,\delta)\times (\tilde{I}_x/\tilde{I}_x\cap\mm_x^{k+1})\to\tau^*J^k(\Sigma_i,\tilde{I})
\end{equation}
defined by $\psi(t,\pi^\infty_k(h_x)):=(t,\pi^\infty_k(\phi_{t*}h_x))$ for any $h_x\in\tilde{I}_x$ and any $t\in (-\delta,\delta)$. 
Define the curve
\[\alpha_1:(-\delta,\delta)\to (-\delta,\delta)\times (\tilde{I}_x/\tilde{I}_x\cap\mm_x^{k+1})\]
by the formula $\alpha_1(t):=(t,\pi^\infty_k(f_x))$. Observe that $\gamma=\eta\comp\psi\comp\alpha_1$. 

On the other hand we consider the curve 
\[\sigma:(-\delta,\delta)\to V\times W\]
defined by $\sigma(t):=(\tau(t),s)$. Then 
$\pi^\infty_k\comp\rho_F\comp\sigma(t)=\pi^\infty_k(f_{\tau(t)})=j^kf(\tau(t))$. Defining 
\[\alpha_2:(-\delta,\delta)\to (-\delta,\delta)\times (\tilde{I}_x/\tilde{I}_x\cap\mm_x^{k+1})\]
as $\alpha_2(t):=(t,\pi^\infty_k(\phi_{-t*}f_{\tau(t)}))$ we have $\pi^\infty_k\comp\rho_F\comp\sigma=\eta\comp\psi\comp\alpha_2$.

Consider the following direct sum decomposition:
\begin{equation}
\label{sumadirecta}
T_{(0,p)}[(-\delta,\delta)\times\tilde{I}_x/\tilde{I}_x\cap\mm_x^{k+1}]=\RR\oplus\tilde{I}_x/\tilde{I}_x\cap\mm_x^{k+1}.
\end{equation} 
With respect to it we have $\alpha_1'(0)=(1,0)$ and $\alpha'_2(0)=(1,\pi^\infty_k(X(f_x)))$. Then, decomposing the tangent space of the vector bundle 
$J^k(U,\CC)$ at $\pi^\infty_k(f_x)$ as 
\[T_p J^k(U,\CC)=T_xU \times\calO_{U,x}/\mm^{k+1}\]
we have
\[(\pi^\infty_k\comp\rho_F\comp\sigma)'(0)-\gamma'(0)=d(\eta\comp\psi)_{(0,p)}(\alpha'_2(0)-\alpha'_1(0))=(0,\pi^\infty_k(X(f_x))).\]
As $\pi^\infty_k\comp\rho_F=j^k\varphi\comp\pi^\infty_k\comp\tilde{\rho}_F$ and $\gamma=j^k\varphi\comp\tilde{\gamma}$ we conclude that
\[(0,\pi^\infty_k(X(f_x)))\in d(j^k\varphi)_q\big(d(\pi^\infty_k\comp\tilde{\rho}_F)_{(x,s)}(T_{(x,s)}V\times W)+T_qA\big).\]
Taking into account inclusion~(\ref{primertrozo}) we obtain
\begin{equation}
\label{tercertrozo}
\{0\}\oplus (j^k\varphi_x)^{-1}(\pi^\infty_k(\Theta_{\tilde{I}_x,e}(f_x))\subset d(\pi^\infty_k\comp\tilde{\rho}_F)_{(x,s)}(T_{(x,s)}V\times W)+T_qA.
\end{equation} 

Let $\partial/\partial w_1,...,\partial/\partial w_d$ be a basis of the tangent space $T_{s}W$. The versality of $F$ at $(x,s)$ means
\begin{equation}
\label{consecversalidad}
\CC(\partial F/\partial w_1)_{|s,x}+...+(\partial F/\partial w_d)_{|s,x}+\Theta_{\tilde{I}_x,e}(f_x)=\tilde{I}_x.
\end{equation}
As $d(j^k\varphi\comp\pi^\infty_k\tilde{\rho}_F)_{(x,s)}(0,\partial/\partial w_i)=(0,\pi^\infty_F((\partial F/\partial w_i)_{|s,x})$, 
equality~(\ref{consecversalidad}) together with inclusion~(\ref{tercertrozo}) imply
\begin{equation}
\label{cuartotrozo}
\{0\}\oplus E\subset d(\pi^\infty_k\comp\tilde{\rho}_F)_{(x,s)}(T_{(x,s)}V\times W)+T_qA.
\end{equation}

After this it is sufficient to prove that $d(pr_k)_q$ maps the space 
\[d(\pi^\infty_k\comp\tilde{\rho}_F)_{(x,s)}(T_{(x,s)}V\times W)\]
surjectively over $T_xU$, but this is trivial because
$pr_k\comp\pi^\infty_k\comp\tilde{\rho}_F$ is the projection of $V\times W$ to its first factor.
\end{proof}

\section{The Relative Morsification Theorem}

Let $\Delta_1$ be a $\calD_{\tilde{I},e}$-invariant closed analytic subset of $\Gamma(U,\tilde{I})$ with all its irreducible components of codimension at least
$n+1$, and containing $\Gamma^{(r)}$ (where $\Gamma^{(r)}$ was introduced in Theorem~\ref{Varchenko}). Recall the set $\Delta_2$, appearing in the Topological 
Partition~(\ref{partition}). Define $\Delta:=\Delta_1\cup\Delta_2$. 
The partition $\calP_0$ defined by:
\begin{equation}
\label{primer}
J^\infty(U,\tilde{I}):=\big[\coprod_{i\leq n,l\in L_i}B_{i,l}\setminus\Delta\big]\coprod \Delta,
\end{equation}
is a partition of $J^\infty(U,\tilde{I})$ by locally closed analytic $\calD_{\tilde{I},e}$-invariant subsets, which is closely related with the topological 
partition: they only differ in subsets of codimension strictly bigger than $n$.
It has the advantage that all the terms involved in its definition are locally closed analytic subsets. On the other hand $\Delta$ is not canonically
defined, but its irreducible components have codimension at least $n+1$. As we will see, subsets of codimension bigger than $n+1$ in 
are too small to affect the singularity types appearing in a generic deformation of any function of $I$ of finite extended codimension.

By Proposition~\ref{deter3} the level sets $C_m$ of the filtration by extended codimension are finite-determined closed analytic subsets if $m<\infty$. Let 
$\{C_{1,j}\}_{j\in J}$ be the set of irreducible components of $C_1$. Given $C_{1,j}$ we let $m(j)$ be the extended codimension of a generic element in it, 
that is, then minimal $m$ such that the intersection $C_{i,j}\cap\dot{C}_m$ is not empty. Define $\partial C_{1,j}:=C_{1,j}\cap C_{m(j)+1}$. Observe that if
$m(j)=\infty$ then $\partial C_{1,j}=C_{1,j}$. The local finiteness of $\{\partial C_{1,j}\}_{j\in J}$ follows from the locally finiteness of 
$\{C_{1,j}\}_{j\in J}$. Therefore $\partial C_i:=\cup_{j\in J}\partial C_{1,j}$ is a closed analytic subset. We have a canonically defined partition by 
$\calD_{\tilde{I},e}$-invariant locally closed subsets
\begin{equation}
\label{segund}
J^\infty(U,\tilde{I})=\dot{C}_0\coprod (C_1\setminus\partial C_1)\coprod\partial C_1
\end{equation}
satisfying that the irreducible components of $C_1\setminus\partial C_1$ are at least $n$-codimensional and 
the irreducible components of $\partial C_1$ are either of codimension strictly bigger than $n$ or are contained in $C_\infty$.
Indeed: the irreducible components of $C_1\setminus\partial C_1$ are of the form $C_{1,j}\setminus\partial C_1$ for $j\in J$. If $\codim(C_{1,j})<n$, by 
Proposition~\ref{codimbound} we have $C_{1,j}\subset C_\infty$, and hence $m(j)=\infty$; in this case $\partial C_{1,j}=C_{1,j}$ and therefore 
$C_{1,j}\setminus\partial C_1$ is empty. Suppose that we have a component $C'$ of $\partial C_1$ not contained in $C_\infty$; then $C'$ is a component of 
$\partial C_{1,j}$ for a certain $j$ such that $C_{1,j}$ is not contained in $C_\infty$; in this case $\codim(C_{1,j})\geq n$ and $C'$ is a proper
closed analytic subset of the irreducible $C_{1,j}$; we deduce $\codim(C')\geq n+1$.

The partition
\begin{equation}
\label{panumerarla}
J^\infty(U,\tilde{I})=\coprod_{k\in K}Z_k,
\end{equation}
whose strata are the subsets of the form $X\cap Y$, where $X$ and $Y$ are respectively strata of the partitions~(\ref{primer}) and (\ref{segund}), is a locally
finite partition by locally closed analytic $\calD_{\tilde{I},e}$-invariant subsets refining the partitions~(\ref{primer}) and (\ref{segund}). It is 
canonically defined up to codimension $n$: the only non-canonical element involved in its definition is the subset $\Delta_1$, whose irreducible components are
at least $n+1$-codimensional. Therefore, given a different choice $\Delta'_1$ for this set, and an stratum $Z'_i$ of codimension smaller or equal than $n$
of the resulting partition, it is easy to see
that there exists an stratum $Z_k$ in the original partition with the same closure $\overline{Z}$ than $Z'_i$ such that 
\[\codim(\overline{Z}\setminus (Z_k\cap Z'_i)\geq n+1.\]
If $Z_k$ is an stratum such that $\codim(Z_k)\leq n$ then any two germs in $Z_k$ are topologically equivalent and have the same extended codimension. 

\begin{definition}
\label{canwhit}
The canonical Whitney stratification $\calX:=\{X_j\}_{j\in J}$ associated (by 
Theorem~\ref{Whitneystr}) to the partition~(\ref{panumerarla}) is called the 
{\em $n$-canonical Whitney stratification} of $J^\infty(U,\tilde{I})$. 
\end{definition}

\begin{remark}
The $n$-\emph{canonical Whitney stratification} of $J^\infty(U,\tilde{I})$ is canonically defined up to codimension $n$, in the sense explained above. 
If $X_j$ is an stratum such that $\codim(X_j)\leq n$ then any two germs in $X_j$ are topologically equivalent and have the same extended codimension.
\end{remark}

\begin{definition}
\label{unavoidable}
A stratum of the Topological Partition or of the $n$-canonical Whitney stratification of 
$J^\infty(U,\tilde{I})$ is called {\em thick} if it has a irreducible component
of codimension at most $n$.
\end{definition}

Consider $f\in I$ defined on a neighbourhood $V$ of the origin $O$ in $U$ such that 
$c_{I,e}(f)<\infty$. Let $\rho_f:V\to J^\infty(U,\tilde{I})$ be its associated 
jet-extension. By the local finiteness of $\calX$ we can shrink $V$ such that 
$\rho_f(V)$ only meets finitely many strata. Therefore there exist a radius 
$\epsilon$ such that $\rho_f(B_{2\epsilon})$ only meets the strata of $\calX$ 
whose closure contain $\rho_f(O)$. Let $X_1$,...,$X_s$ be such strata. If $\epsilon$ 
is small enough we know that the only point where $f$ has positive extended 
codimension in $B_{2\epsilon}$ is the origin. Therefore $f$ is its own versal unfolding for any $x\neq O$, and, by Proposition~\ref{versalimplytrans}, the jet 
extension $\rho_f$ is transversal to each stratum $X_i$ at $x$. Consequently $\rho_f^{-1}(X_i)$ is smooth 
outside the origin and $\{\rho_{f|B_{2\epsilon}}^{-1}(X_i)\}_{i\leq s}$ defines a 
Whitney stratification in $B_{2\epsilon}\setminus\{O\}$. Hence, if $\epsilon$ is small 
enough we can assume that $S_\epsilon$ is transversal to any stratum 
$\rho_{f|B_\epsilon}^{-1}(X_i)$. Any such radius $\epsilon$ satisfying the above 
properties is called a {\em good radius} for $f$. From this moment a {\em good system of radii} for $f$ is a pair $(\epsilon,\eta)$ satisfying the 
conditions imposed in Definition~\ref{gsr}, and such that $\epsilon$ is a good radius for $f$. 

\begin{definition}
\label{unavoby}
A stratum of the topological partition or of the $n$-canonical Whitney stratification of 
$J^\infty(U,\tilde{I})$ is called {\em unavoidable by} $f$ if it is 
thick and $\rho_f(O)$ belongs to the closure of one of its components of 
codimension at most $n$.
\end{definition}

\begin{definition}
\label{unsplittable}
Let $f\in\tilde{I}_x$ such that $c_{\tilde{I}_x,e}<\infty$. We say that $f_x$ is 
{\em codimensionally irreducible} if there not exists an unfolding 
$F:(\CC^n,0)\times (\CC,0)$ and a sequence $\{(x_n,s_n)\}_{n\in\NN}$ converging to 
$(x,0)$ such that $c_{\tilde{I}_{x_n},e}(F_{|s_n,x_n})<c_{\tilde{I}_x,e}$.
\end{definition}
   
Let $(A_1,...,A_s)$ and $(A'_1,...,A'_s)$ be two tuples of topological spaces such
that $A_i\subset A_1$ and $A'_j\subset A'_1$ for any $i,j\leq s$. We say that the two tuples are {\em 
topologically equivalent} if there exists an homeomorphism $h:A_1\to A'_1$ such that
$h(A_i)=A'_i$ for any $i$.

\begin{theo}[The Relative Morsification Theorem]
\label{morsificacion}
Let $f\in I$  such that $c_{I,e}(f)<\infty$, and let $F:U\times V\to\CC$ be a 
representative of a versal unfolding of $f$. Let $\{X_1,...,X_s\}$ and 
$\{B_1,...,B_k\}$ be respectively the strata of the $n$-canonical Whitney stratification 
and of the topological partition that are unavoidable by $f$. 
Let $\epsilon$ be a good radius for $f$. The subset $V$ can be chosen small enough so that there exists a 
closed analytic subset $\Delta\subset V$ (called the {\em discriminant} of $F$), and a positive number $\eta$ such 
that 
\[\rho_{F_{|s}}:B_{\epsilon+\eta}\times\{s\}\to J^{\infty}(U,\tilde{I})\]
is transversal to the canonical Whitney stratification $\calX$ for any 
$s\in V\setminus\Delta$. As a consequence 
\begin{enumerate} 
\item If $s\in V\setminus\Delta$ then the image of $\rho_{F_{|s}}$ meets only strata 
of $\calX$ or of the Topological Partition that are unavoidable by $f$.
\item The topological type of the tuples 
\begin{equation}
\label{tuplefina}
(\overline{B}_\epsilon,\rho_{F_s}^{-1}(X_1),...,\rho_{F_s}^{-1}(X_s)),
\end{equation}
\begin{equation}
\label{tuplegruesa}
(\overline{B}_\epsilon,\rho_{F_s}^{-1}(B_1),...,\rho_{F_s}^{-1}(B_k))
\end{equation} 
does not depend on $s\in V\setminus\Delta$. Moreover the topological type of these 
tuples does not depend on the chosen versal unfolding.
\item The germ of $F_{|s,x}$ is codimensionally irreducible for any 
$x\in B_{\epsilon}$ and any $s\in V\setminus\Delta$.
\item The number of points of $\overline{B}_\epsilon$ where $F_{|s}$ has a fixed extended 
codimension does not depend neither on $s\in V\setminus\Delta$ nor on the choice of 
the versal unfolding. The points of $B_\epsilon\setminus V(I)$ where $F_{|s}$ has positive extended 
codimension are $A_1$-singularities. Moreover $F_{|s}$ has 
extended codimension $0$ along $\partial B_\epsilon$. 
\item The analytic type of the germ $(V,\Delta,O)$ (where $O$ is the origin
of $V$) does not depend on the choice of the unfolding as long it is universal (that is versal with $\dim(V)=c_{I,e}(f)$).
\end{enumerate}
\end{theo}
\begin{proof}
As $F$ is versal, by Proposition~\ref{versalimplytrans}, the mapping $\rho_{F}$ is 
transversal to $\calX$. 
Using the same argument than in Lemma~\ref{controlext} choosing $V$ small enough we
can assume that $c_{\tilde{I}_x,e}(F_{|s,x})=0$ for any 
$x\in B_{2\epsilon}\setminus B_{\epsilon/2}\times V$; therefore for such points 
the germ $F_{|s,x}$ is its own versal unfolding. Hence, by 
Proposition~\ref{versalimplytrans}, the mapping $\rho_{F_{|s}}$ is transversal to
$\calX$ for any $x\in B_{2\epsilon}\setminus B_{\epsilon/2}$. Choose $\eta<\epsilon$. Applying 
Theorem~\ref{transpar} we conclude the existence of a closed analytic subset $\Delta$
such that $\rho_{F_s}\pitchfork\calX$ for any $s\in V\setminus\Delta$. Now we derive 
the stated consequences. 

As the $n$-canonical Whitney stratification is a subdivision of the topological partition
the assertions concerning the former are imply the analogous assertions 
concerning the later. 

Let $s\in V\setminus\Delta$; as $\dim(B_\epsilon)=n$, by transversality, the image of
$\rho_{F_{|s}}$ only can meet strata of $\calX$ of codimension at most $n$, hence 
thick. By local finiteness of $\calX$ and the fact that the radius $\epsilon$
is a {\em good radius} for $f$, we can choose $V$ small enough so that the closure
of any stratum of $\calX$ which is met by $\rho_F(B_{\epsilon+\eta}\times V)$ necessarily 
contains the origin. Therefore the image of $\rho_{F_{|s}}$ only can meet strata that
are unavoidable by $f$. This shows assertion~(1). 

Let $X_1,...,X_s,X_{s+1},...,X_{l}$ be the strata of $\calX$ whose closure meet 
$\rho_f(O)$. By the transversality $F\pitchfork\calX$ we have that 
\[\calY=\{Y_i:=\rho_F^{1}(X_i)\}_{i\leq l}\]
is a Whitney stratification of $B_{\epsilon+\eta}\times V$. Let $\pi:B_{\epsilon+\eta}\times V$ be the projection to the second factor. Define
$Y_{i|V\setminus\Delta}:=Y_i\cap\pi^{-1}(V\setminus\Delta)$. For any $s\in V\setminus\Delta$, the transversality $\rho_{F_{|s}}\pitchfork\calX$ implies that 
either $Y_{i|V\setminus\Delta}$ is empty or it projects submersively to  $V\setminus\Delta$. On the other hand, as $\epsilon$ is 
a good radius, choosing $V$ small enough we have that $S_\epsilon\pitchfork Y_{i|s}$ 
for any $i\leq l$ and any $s\in V$. This implies that the Whitney stratifications 
$\calY$ and $\{B_\epsilon,S_\epsilon\}$ meet transversely, and that, if $\calZ$ is 
the Whitney stratification formed by pairwise intersection of strata of the previous 
one, the restriction of the projection 
\[\pi:\overline{B}_\epsilon\times (V\setminus\Delta)\to V\setminus\Delta\] 
to any stratum of $\calZ$ is submersive. As $\pi$ is proper we have that $\calZ$ is 
topologically trivial over $V\setminus\Delta$. This proves the independence on 
$s\in V\setminus\Delta$ of the topological type of the tuple~(\ref{tuplefina}).

Let $F:U\times V\to\CC$ and $F':U'\times V'\to\CC$ be two versal unfoldings. 
The mapping $F'':U\times (V\times V')\to\CC$ defined by 
$F_{|(v,v')}:=F_v+F_{v'}-f$ is a $I$-unfolding of $F$ such that 
$F''_{|V\times\{O'\}}=F'$ and $F''_{|\{O\}\times V'}=F'$ (where $O$ and $O'$ are the 
respective origins of $V$ and $V'$). Clearly $F''$ is versal, and moreover its 
discriminant does not contain neither $V\times\{O'\}$ nor $\{O\}\times V'$. The 
independence on $s\in V\times V'\setminus\Delta$ of the topological type of the 
tuple~(\ref{tuplefina}) for the unfolding $F''$ implies the independence on the choice
of versal unfolding. This finishes the proof of Assertion~(2). 

Given $s\in V$, the set of points of $B_{\epsilon+\eta}$ where $F_{|s}$ has extended 
codimension is $\rho_{F_{|s}}^{-1}(C_1)$. For any $i>0$ we let $Z_i$ be the union of the strata among $X_1,...,X_s$ which contain germs of extended codimension
precisely $i$; define $Z:=\cup_{i>0}Z_i$. Assertion $(1)$ implies that for any $s\in V\setminus\Delta$ we have $\rho_{F_{|s}}^{-1}(C_1)=\rho_{F_{|s}}^{-1}(Z)$ 
and that, for any $i>0$, the set of points where $F_{|s}$ has extended codimension $I$ equals $\rho_{F_{|s}}^{-1}(Z_i)$. 
The topological type of the tuple~(\ref{tuplefina}) determine the topological types of the subsets $\rho_{F_{|s}}^{-1}(Z)$ and $\rho_{F_{|s}}^{-1}(Z_i)$ 
for any $i>0$. Consequently, by Assertion~(2), this topological types are independent of $s\in V\setminus\Delta$. 

The independence of the topological type of $\rho_{F_{|s}}^{-1}(Z)$ implies Assertion~(3). As the topological type of $\rho_{F_{|s}}^{-1}(Z_i)$ determines 
the number of points where $F_{|s}$ has extended codimension $i$, Assertion~(2) also implies the first part of Assertion~(4).
The fact that the points of $B_\epsilon\setminus V(I)$ where $F_s$ points of positive extended codimension are of type $A_1$ follows 
from Assertion (3): if $x$ is such a point then $\tilde{I}_x=\calO_{\CC^n,x}$ and, by the morsification theory for isolated singularities,  the only 
codimensionally irreducible singularities in this case are Morse points. 

The fifth assertion is easy taking into account that the strata of $\calX$ are $\calD_{\tilde{I},e}$-invariant.
\end{proof}

The fact that the topological type of the tuple~(\ref{tuplegruesa}) does not depend 
on $s\in V\setminus\Delta$ can be phrased saying that if $g$ and $h$ are any two 
generic close approximations of $f$ in $I$, the partition of $B_\epsilon$
by topological type of $g$ is homeomorphic to the analogous partition for $h$.

\begin{adendum}
The projections to the second factor 
\[\big(\overline{B}_\epsilon\times V\setminus\Delta,\rho_{F|\overline{B}_\epsilon\setminus\Delta}^{-1}(X_1),,...,\rho_{F|\overline{B}_\epsilon\setminus\Delta}^{-1}(X_s)\big)\to V\setminus\Delta\]
\[\big(\overline{B}_\epsilon\times V\setminus\Delta,\rho_{F|\overline{B}_\epsilon\setminus\Delta}^{-1}(B_1),,...,\rho_{F|\overline{B}_\epsilon\setminus\Delta}^{-1}(B_k)\big)\to V\setminus\Delta\]
are locally trivial fibrations of tuples of topological spaces with fibre the tuples~(\ref{tuplefina}) and~(\ref{tuplegruesa}) respectively. The topological 
type of these fibrations is a $\calD_I$-invariant.
\end{adendum}  
 
\begin{adendum}
\label{imagenes}
With the notations of last theorem, let $r$ be the number of $A_1$-points in $B_\epsilon\setminus V(I)$ of a generic deformation $F_{|s}$ of $f$ within $I$. 
There exists a 
closed analytic proper subset $\Delta'\subset V$ such that its complement $V\setminus\Delta'$ is the set of $s\in V$ such that the images by $F_{|s}$ of the 
critical points of $F_{|s}$ in $B_\epsilon\setminus V(f)$ are $r$ different points which are also different to $0$ if $V(f)\neq\emptyset$.  
\end{adendum} 
\begin{proof}
First we show that the set of values of the parameter $s$ with the property of the statement is dense in $V$. Consider $s\in V\setminus\Delta$. Then the only 
singularities of $F_{|s}$ in 
$B_\epsilon\setminus V(I)$ are $r$ Morse points, which we denote it by $p_1,...,p_r$. Suppose that, for any $i\leq r$, we manage to construct an $I$-unfolding 
$G^i:=F_{|s}+th$ such that for any $t$ small enough the critical points of $G^i_{|t}$ in $B_\epsilon\setminus V(I)$ are $p_1,...,p_r$ and
\[G^i_{|t}(p_i)=F_{|s}(p_i)+t\]
\[G^i_{|t}(p_j)=F_{|s}(p_j)\]
for $i\neq j$. 
Then, using the versality of $F$, we deduce easily that for $s'$ close enough to $s$ the critical points of $F_{|s'}$ in $B_\epsilon\setminus V(I)$ are 
$p_1,...,p_r$, their images are pairwise different and different to $0$.

Let us construct the $I$-unfolding $G^1$. As $p_1,...,p_r$ are different and do not belong to $V(I)$ there exist a function $g\in I$ vanishing at 
$p_2$,...,$p_r$ and such that $g(p_1)=1$. Then the function $h:=(2-g^2)g^2$ belongs to $I$, vanish and is singular at $p_2$,...,$p_r$, takes value $1$ and 
is singular at $p_1$. Define $G^1:=F_{|s}+th$.   

Consider the analytic mapping $\sigma:V\setminus\Delta\to\CC^r$ whose $i$-th component $\sigma_i$ assigns to $s$ the $i$-th symmetric function evaluated in the
images of the critical points of $F_{|s}$ in $B_\epsilon\to\CC\setminus V(f)$. Let $\beta:\CC^r\to\CC$ be the analytic mapping associating to $a_1,...,a_r$ the
discriminant of the polynomial $T^r+\sum_{i=1}^ra_iT^{r-i}$. Define $\varphi:V\setminus\Delta\to\CC$ to be the composition $\varphi:=\beta\comp\sigma$.
The function $F$ is defined in $\overline{B}_\epsilon\times V$, therefore, if $V$ is small enough then $F(B_\epsilon\times V)$ is a bounded set in $\CC$. 
Hence
$\sigma$, and consequently $\varphi$ are bounded functions. By Riemann extension theorem we can suppose that the analytic function $\varphi$ is defined in the
whole $V$. It is easy to check that defining $\Delta':=\varphi^{-1}(0)$ if $V(f)=\emptyset$ and 
$\Delta':=\varphi^{-1}(0)\cup\sigma_r^{-1}(0)$ if $V(f)\neq\emptyset$ the statement of the proposition holds.
\end{proof}

\begin{definition}
\label{numinv} 
The minimal closed analytic subset $\Delta\subset V$ with the properties of Theorem~\ref{morsificacion} is called the {\em discriminant} of $F$ in $V$. The 
subset $\Delta'$ introduced in Proposition~\ref{imagenes} is called the {\em bifurcation variety} of $F$ in $V$. The union $\Delta\cup\Delta'$ is called the 
{\em big discriminant} of $F$ in $V$.
\end{definition}
 
The bifurcation variety defined here generalises the bifurcation variety for isolated singularities and the one studied in~\cite{Za2}.

\begin{definition}
\label{morsif}
Let $f\in I$ be a function of finite extended codimension and $\epsilon$ a good radius for it. An $I$-unfolding $F:(\CC^n\times\CC,(O,0))\to\CC$ is called a 
{\em Morsification} if for any small enough $s\neq 0$ the jet-extension $\rho_{F_{|s}}:(\CC^n,O)\to J^\infty(U,\tilde{I})$ is transversal $n$-canonical Whitney
stratification in a neighbourhood of the origin containing $\overline{B}_\epsilon$, and the images of the isolated critical points of $F_{|s}$ are pairwise 
different and, if $V(I)\neq\emptyset$, also different from $0$. 
\end{definition}

A consequence of our result is that any generic uniparametric $I$-unfolding is a Morsification.

\section{Applications} 

In this section we explain two applications of our Relative Morsification Theorem

\subsection{Numerical invariants} 

In Example~\ref{ejem} we proved that in general the extended codimension does not behave in an conservative way with respect to deformation, unlike the Milnor 
number in the case of isolated singularities. We introduce now some numerical $\calD_I$-invariants for functions of $I$ which are conservative and 
related to the extended codimension.

\begin{definition}
\label{splittingfunction} 
Let $f\in I$ be a function such that $c_{I,e}(f)<\infty$. Let $\epsilon$ be small enough so that the origin is the only point of $B_\epsilon$ where $f$ has
positive extended codimension. Define the {\em splitting function}   
\begin{equation}
\sigma_I[f]:\NN\to\ZZ_{\geq 0}
\end{equation}
of $f$ with respect to $I$ imposing that $\sigma_f(k)$ is the number of points of $B_\epsilon$ where $F_{|s}$ has extended codimension $k$, being $F$ a versal
$I$-unfolding and $s$ a parameter not contained in the discriminant. Define the {\em corrected extended codimension} $\tilde{c}_{I,e}(f)$ of $f$ with respect
to $I$ by the formula 
\begin{equation}
\tilde{c}_{I,e}(f):=\sum_{k\in\NN}k\sigma_I[f](k).
\end{equation}

Define the {\em Morse number} $\calM_{I}(f)$ of $f$ with respect to $I$ as the number of points of $B_\epsilon\setminus V(I)$ where $F_{|s}$ has an $A_1$ 
singularity.
\end{definition}

Clearly $\tilde{c}_{I,e}(f)\leq c_{I,e}(f)$. Example~\ref{ejem} shows that the inequality may be strict. By formula~(\ref{ineqcod}) we have 
$\sigma_I[f](k)\leq c_{I,e}(f)/k$. The Morse number, the extended codimension and the corrected extended codimension coincides with the Milnor number when
$I=\calO_{\CC^n,O}$. 

Using versal unfoldings it is easy to check the conservativity of the invariants defined above: consider $f\in I$, a function 
such that $c_{I,e}(f)<\infty$; let $\epsilon$ be small enough so that the origin is the only point of $B_\epsilon$ where $f$ has
positive extended codimension. Given any $I$-unfolding $F:\CC^n\times\CC^k\to\CC$ and any small enough $s\in\CC^k$ of the parameter we have 
\begin{equation}
\label{aditividad}
\sum_{x\in B_\epsilon}\alpha_{\tilde{I}_x}(f)=\alpha_I(f_O),
\end{equation}  
where $\alpha_{\tilde{I}_x}(f_x)$ stands for $\sigma_{\tilde{I}_x}[f_x]$, $\tilde{c}_{\tilde{I}_x}(f_x)$ or $\calM_{\tilde{I}_x}(f_x)$. 

\begin{remark}
Obviously, more refined numerical invariants can be defined taking into account the distribution of the points of a fixed extended codimension of a generic 
deformation $F_{|s}$ in the partitions given by the tuples~(\ref{tuplefina})~and~(\ref{tuplegruesa}), or even considering numerical topological invariants 
of the tuples.
\end{remark}

\subsection{Consequences on the topology of the Milnor fibre}

In this section we deduce topological properties of Milnor fibres with the help of Morsifications.

Let $V=V(I)$ be the analytic germ defined by $I$. Suppose $V\neq\emptyset$. Due to the conical structure of analytic germs there exist $\epsilon_0>0$ such that
$V\cap B_\epsilon$ is contractible for any $\epsilon<\epsilon_0$. 

Recall that we have fixed generators $g_1,...,g_r$ of $\tilde{I}$ in the open subset $U$. 
Consider the real analytic function $\kappa:\overline{B}_\epsilon\to\RR$ by the formula $\kappa(x):=\sum_{i=1}^r|g_i(x)|^2$. We claim that $0$ is 
an isolated critical value of $\kappa$: otherwise, by the Curve Selection Lemma, there exists a germ of analytic path $\gamma:(\RR,0)\to\overline{B}_\epsilon$ 
such that $\gamma^{-1}(V)=\{0\}$ and $(\kappa\comp\gamma)'(t)=0$ for any $t$. This is a contradiction, and hence our claim is true. An analogous reasoning 
yields that $0$ is an isolated critical value of the restriction $\kappa_{|S_\epsilon}$. As a consequence we can choose a positive $\xi_0$ such that  
$0$ is the only critical value of $\kappa$ in $[0,\xi_0]$ and $\kappa^{-1}(\xi)\pitchfork S_\epsilon$ for any $\xi\in(0,\xi_0]$. Then, for any 
$\xi\in (0,\xi_0]$, the set $N_{\epsilon,\xi}:=\kappa^{-1}[0,\xi]$ is a compact neighbourhood of $V\cap\overline{B}_\epsilon$ in $\overline{B}_\epsilon$, whose
boundary is a manifold with corners. Moreover, by Ehresmann fibration theorem the mapping 
\begin{equation}
\label{fib1}
\kappa_{|N_{\epsilon,\xi}\setminus\kappa^{-1}(0)}:N_{\epsilon,\xi}\setminus\kappa^{-1}(0)\to (0,\xi]
\end{equation}
is a locally trivial fibration. 

Consider $f\in I$ with $c_{I,e}(f)<\infty$. Let $F$ be a $I$-morsification of $f$, and $h:=F_{|s}$ for a certain $s$ to be chosen later. We can assume that 
$\epsilon$ is small enough so that there is $\eta>0$ such that $(\epsilon,\eta)$ is a good system of radii for $f$. Take a small enough positive number 
$\delta$ so that the conclusions of Theorem~\ref{topconst} are fulfilled; choose $s\in D_\delta$. Then, for any $t\in D_\eta$ we have
$h^{-1}(t)\pitchfork S_\epsilon$ (the transversality is meant in a stratified sense when $t=0$). 

Define 
\[Y_{\epsilon,\xi,\eta}:=N_{\epsilon,\xi}\cap h^{-1}(D_\eta)\quad\quad \dot{Y}_{\epsilon,\xi,\eta}:=Y_{\epsilon,\xi,\eta}\setminus h^{-1}(0).\]
We claim that there is a positive $\xi\leq\xi_0$ small enough so that there exists $\eta$ satisfying that 
\begin{equation}
\label{fibracioncentral}
h_{|\dot{Y}_{\epsilon,\xi,\eta}}:\dot{Y}_{\epsilon,\xi,\eta}\to D_\eta\setminus\{0\}
\end{equation}
is a locally trivial fibration, and, moreover, if $(\xi',\eta')$, with $\xi'\leq\xi$, is a pair of radii defining an analogous fibration, then this fibration
is equivalent to the one associated with $(\xi,\eta)$.  

As the critical points of $h$ in $\overline{B}_\epsilon$ are either points contained in $V$ or isolated points of positive extended codimension it is 
clear that choosing $\xi$ small enough we can assume that the only critical points that $h$ has in $N_{\epsilon,\xi}$ are included in $V$. Therefore the 
image of all the critical points is $0$. 

The boundary of $N_{\epsilon,\xi}$ is a manifold with corners admitting the following decomposition in smooth strata 
\[\partial N_{\epsilon,\xi}=(\kappa^{-1}([0,\xi))\cap S_\epsilon)\coprod (\kappa^{-1}(\xi)\cap S_\epsilon)\coprod (\kappa^{-1}(\xi)\cap B_\epsilon).\]

As $h^{-1}(t)\pitchfork S_\epsilon$ for any $t\in D_\epsilon\setminus\{0\}$
and $\kappa^{-1}([0,\xi))\cap S_\epsilon$ is open in $S_\epsilon$ we deduce 
\[h^{-1}(t)\pitchfork \kappa^{-1}([0,\xi))\]
for any $t\in D_\epsilon\setminus\{0\}$. 

On the other hand $h^{-1}(0)$ is smooth outside
$V$. We show that there exists $\xi$ such that
\begin{equation}
\label{auxtran}
h^{-1}(0)\pitchfork (\kappa^{-1}(\xi')\cap B_\epsilon)\quad\quad h^{-1}(0)\pitchfork (\kappa^{-1}(\xi')\cap S_\epsilon)
\end{equation}
for any $\xi'<\xi$.
Suppose that for any $\xi>0$ the set of points where $h^{-1}(0)$ is not transversal to $(\kappa^{-1}(\xi)\cap B_\epsilon)$ is not void. Then by the 
Curve Selection Lemma there is a germ of analytic path $\gamma:(\RR,0)\to h^{-1}(0)$ such that $\gamma^{-1}(V)=\{0\}$ and 
$(\kappa\comp\gamma)'(t)=0$ for any $t\in (-\delta,\delta)$. Then the function $\kappa\comp\gamma$ has zero derivative, but is not constant.
this gives a contradiction which shows $h^{-1}(0)\pitchfork (\kappa^{-1}(\xi')\cap B_\epsilon)$ for $\xi'$ small. The transversality 
$h^{-1}(0)\pitchfork (\kappa^{-1}(\xi')\cap S_\epsilon)$ when $\xi'$ is small enough is proven analogously. Using the compactness of 
$(\kappa^{-1}(\xi)\cap S_\epsilon)\coprod (\kappa^{-1}(\xi)\cap B_\epsilon)$ it is easy to show that if $\eta$ is small enough then 
\[h^{-1}(t)\pitchfork (\kappa^{-1}(\xi)\cap B_\epsilon)\quad\quad h^{-1}(t)\pitchfork (\kappa^{-1}(\xi)\cap S_\epsilon)\]
for any $t\in D_\eta$.

Summarising, we have found $\xi$ and $\eta$ such that $h^{-1}(t)$ is transversal to the manifold with corners $\partial N_{\epsilon,\xi}$ for any 
$t\in\dot{D}_\eta$, and such that the only critical points of $h$ at $N_{\epsilon,\xi}$ are at $h^{-1}(0)$. By Ehreshmann fibration theorem the 
mapping~(\ref{fibracioncentral}) is a locally trivial fibration. Using the fact that~(\ref{auxtran}) holds for any $\xi'<\xi$ it is easy to check that the 
fibration analogous to~(\ref{fibracioncentral}) associated with any other suitable pair $(\xi',\eta')$ satisfying  $\xi'\leq\xi$ is  equivalent to the one 
associated with $(\xi,\eta)$. This shows the claim.

Using both the transversality conditions that we have checked up to now and the fibrations~(\ref{fib1}) and~(\ref{fibracioncentral}), it is easy to 
show that the set $\mathrm{Int}(Y_{\epsilon,\xi,\eta})$ of interior points of $Y_{\epsilon,\xi,\eta}$ is a dilation neighbourhood of $V\cap B_\epsilon$ in the 
sense of~\cite{Ma}. On the other hand, considering a Whitney stratification on $V\cap B_\epsilon$ and taking a controlled tube system for it 
(in the sense of~\cite{GWPL}) we obtain a tubular neighbourhood $T$ of $V\cap B_\epsilon$ which admits a deformation retract to it. By the Uniqueness
Theorem of dilation neighbourhoods of~\cite{Ma} we deduce that $T$ is diffeomorphic to $\mathrm{Int}(Y_{\epsilon,\xi,\eta})$. Therefore the contractibility
of $V\cap B_\epsilon$ implies that $\mathrm{Int}(Y_{\epsilon,\xi,\eta})$, and hence $Y_{\epsilon,\xi,\eta}$, is contractible.
 
The following Theorem shows how to use our theory to extract topological information of the Milnor fibre of $f$. 

Let $F$ denote the Milnor fibre of $f$. Denote by $F_0$ the fibre of the fibration~(\ref{fibracioncentral}). 
\begin{theo}
\label{bouquet}
Let $F$ denote the Milnor fibre of $f$. Denote by $F_0$ the fibre of the fibration~(\ref{fibracioncentral}). Recall that $\calM(f)$ denotes the Morse number of
$F$. Then
\begin{equation}
\label{homology}
\mathrm{H}_k(F,\ZZ)=\mathrm{H}_k(F_0,\ZZ)\bigoplus [\mathrm{H}_k(S^{n-1},\ZZ)]^{\calM(f)}
\end{equation}
for any $k$. 

Moreover, if either $n=2$ or both $F_0$ and $F$ have trivial fundamental group then $F$ has the homotopy type of the bouquet
\begin{equation}
\label{homotopy}
F_0\bigvee[\bigvee_{\calM(f)}S^{n-1}],
\end{equation}
that is, the bouquet of $F_0$ with $\calM(f)$ spheres of dimension $n-1$.
\end{theo}
\begin{proof}
By Theorem~\ref{topconst} the Milnor fibre $F$ is equal to the generic fibre of the mapping~(\ref{qwer}). After this essentially the same arguments 
considered in~\cite{Si3} (for which the contractibility of $Y_{\epsilon,\xi,\eta}$ is needed) allow to conclude.
\end{proof}
 
\begin{adendum}
If the critical locus of $f$ is at least of codimension $3$ then both $F_0$ and $F$ have trivial fundamental group, and, hence, the homotopy decomposition of 
last Theorem holds. This happens, in particular, when $\codim(V(I))\geq 3$.
\end{adendum}
\begin{proof}
Suppose that the codimension of the critical locus of $f$ is at least $3$. Let $S\subset J^\infty(U,\tilde{I})$ be the closed analytic subset 
formed by the singular germs. Then any irreducible component of $S$ which is met by the image of the jet extension $\rho_f:B_\epsilon\to J^\infty(U,\tilde{I})$
has at least codimension $3$. The transversality properties of a Morsification $F_{|s}$ imply $\codim(\rho_{F_{|s}}^{-1}(S))\geq 3$, that is, the codimension
of the critical locus of $F_{|s}$ is at least $3$. Then the results of~\cite{KM} imply the simple connectivity of $F_0$ and $F$.
\end{proof}

Theorem~\ref{bouquet} reduces the study of many properties the Milnor fibre of $f$ to the study of the Milnor fibre of $F_0$. In the study of $F_0$
the fibration~(\ref{fibracioncentral}) our Relative Morsification Theorem becomes very important: it tells us that the only singularities that $F_{|s}$ can have 
belong to strata of the topological partition which are unavoidable by $f$ (this imposes conditions for example in the codimension of such sets of 
singularities). This has been already used successfully for certain classes of ideals (see~\cite{Si1}, \cite{Si2}, \cite{Jo}, \cite{Za}, \cite{Ne}).
The transversality of the jet extension $\rho_{F_{|s}}$ to the relevant strata of $J^\infty(U,\tilde{I})$ relates the relative 
position of the different singularity types appearing in $\rho_{F_{|s}}$ with the relative position of the relevant strata in $J^\infty(U,\tilde{I})$. 

\section{Numerical invariants and intersection multiplicities}
     
In this section we give a characterisation of the invariants introduced in Definition~\ref{numinv} in terms of intersection multiplicities in the generalized 
jet space. We will use certain intersection theoretic
constructions of~\cite{Fu}, that can be checked easily that remain valid in the analytic setting, at least in the very restricted degree of generality that we 
need them.

Let $V\subset J^\infty(U,\tilde{I})$ be an irreducible closed analytic subset of codimension $n$. Given any $f\in\tilde{I}_x$, 
its jet extension $\rho_f:(\CC^n,x)\to J^\infty(U,\tilde{I})$ is a germ of analytic mapping. Suppose that $x$ is an isolated point of 
$\rho_f^{-1}(V)$; as we work with germs at $x$ we can actually assume $\rho_f^{-1}(V)=\{x\}$. Express $f=\sum_{i=1}^rf_ig_i$ where 
$g_1,...,g_r$ is our fixed set of generators and $f_1,...,f_r$ are holomorphic in $x$. We have the associated local lifting 
$\tilde{\rho}_f:(\CC^n,x)\to J^\infty(U,\CC^r)$ given by $\tilde{\rho}_f(y)=(j^\infty f_1(y),...,j^\infty f_r(y))$. Choose an integer $k$ such that $V$ is 
$k$-determined; then $V':=j^k\varphi^{-1}(\pi^\infty_k(V))$ is an irreducible variety in $J^k(U,\CC^r)$. Define 
$\tilde{\rho}^k_f:=\pi^\infty_k\comp\tilde{\rho}_f$; let 
$x':=\tilde{\rho}^k_f(x)$. Let $Z_{\dim(V')}(V')$ and $Z_0(\{x'\})$ be respectively the groups of analytic cycles of dimension $\dim(V')$ and $0$ of 
$V'$ and $\{x'\}$. These groups are obviously isomorphic to $\ZZ$ with respective generators $[V']$ and $\{x'\}$. 
By Definition~8.1.2. of~\cite{Fu} there is a refined Gysin homomorphism 
\[(\tilde{\rho}^k_f)^!:Z_{\dim(V')}(V')\to Z_0(\{x'\}.\] 

\begin{definition}
\label{intmult}
We define the {\em intersection multiplicity} of $\rho_f$ and $V$ at $x$ to be the integer $i_x(\rho_f,V)$ characterised by 
\[(\tilde{\rho}^k_f)^!([V'])=i_x(\rho_f,V)[x'].\]
\end{definition}

We have to prove that the last definition is independent on $k$ and on the functions $f_1,...,f_r$ giving rise to the local lifting. Moreover we want to give a
formula to compute $i_x(\rho_f,V)$. For this we need to recall from~\cite{Fu} how the intersection product $(\tilde{\rho}^k_f)^!([V'])$ is 
defined. 

Consider the subvariety $V'':=\CC^n\times V'$ of the product $\CC^n\times J^k(U,\CC^r)$; denote by $\sigma_k$ the projection to the first factor. 
Let $\gamma^k_f:\CC^n\to\CC^n\times J^k(U,\CC^r)$ be defined by 
$\gamma_f(y):=(y,\tilde{\rho}^k_f(y))$; its image $\Gamma^k_f$ is the graph of $\tilde{\rho}^k_f$; let $x'':=\gamma^k_f(x)$. Observe that $J^k(U,\CC^r)$ is 
isomorphic to $U\times\CC^N$ for a certain $N$; let $pr_1$ and $pr_2$ be the projections to the first and second factor. Recall that we have fixed coordinates
$(x_1,...,x_n)$ for $\CC^n$. Let $z_1,...,z_N$ be a coordinate system for $\CC^N$; 
consider the coordinate system $(y_1,...,y_{n+N})$ of $J^k(U,\CC^r)$ defined by $y_i:=x_i\comp pr_1$ for $1\leq i\leq n$, and $y_{n+i}:=z_i\comp pr_2$ for 
$1\leq i\leq N$; then $\{x_1,..,x_n,y_1,...,y_{n+N}\}$ is a coordinate system for $\CC^n\times J^k(U,\CC^r)$. Define $h_i:=y_{n+i}\comp\gamma^k_f\comp\sigma_k$
for $i\leq N$. Then the subvariety $\Gamma^k_f\subset \CC^n\times J^k(U,\CC^r)$ is defined by the 
regular sequence 
\[(y_1-x_1,...,y_n-x_n,y_{n+1}-h_1,...,y_{n+N}-h_N).\]
Composing with the natural ring epimorphism 
$\calO_{\CC^n\times J^k(U,\CC^n),x''}\to\calO_{V'',x''}$ we
obtain a sequence $\mathbf{s}=(s_1,...,s_{n+N})$ of elements of $\calO_{V'',x''}$. Let $\mathrm{K}_{\bullet}(\mathbf{s})$ be the Koszul complex associated to 
$\mathbf{s}$; denote by $\mathrm{H}_i(\mathrm{K}_{\bullet}(\mathbf{s}))$ its i-th homology module. Then unwinding the definition of $(\tilde{\rho}^k_f)^!$ and
applying Example~7.1.2 of~\cite{Fu} we obtain:

\begin{equation}
\label{intmufor}
i_x(\rho_f,V):=\sum_{i=1}^{n+N}(-1)^i\dim_{\CC}(\mathrm{H}_i(\mathrm{K}_{\bullet}(\mathbf{s}))).
\end{equation}

\begin{remark}
If $\mathbf{s}$ is a regular sequence then $\mathrm{H}_i(\mathrm{K}_{\bullet}(\mathbf{s})=0$ for $i>0$. Then 
\[i_x(\rho_f,V)=\dim_\CC(\mathrm{H}_0(\mathrm{K}_{\bullet}(\mathbf{s})))=\dim_\CC(\calO_{\CC^n,c}/\tilde{\rho}^{k*}_f J),\]
where $J$ is the ideal sheaf of $V'$. This happens when $V'$ is a Cohen-Macaulay variety.
\end{remark}
   
Suppose that $l>k$. The projection $\pi^l_k:J^l(U,\CC^r)\to J^k(U,\CC^r)$ is a trivial fibration with fibre $\CC^L$, for a certain $L$; hence we have a product
decomposition 
\begin{equation}
\label{pppddd}
J^l(U,\CC^r)=J^k(U,\CC^r)\times\CC^L,
\end{equation}
let $q_1$ and $q_2$ be the projections to the first and second factors. Let $w_1,...,w_L$ be a 
coordinate system for $\CC^L$. Define a coordinate system $Y_1,...,Y_{n+N+L}$ of $J^l(U,\CC^r)$ by $Y_i:=q_1^*y_i$ for $1\leq i\leq n+N$ and 
$y_{n+N+i}:=q_2^*w_i$ for $1\leq i\leq L$. The subvariety $W':=(j^l\varphi)^{-1}(\pi^\infty_l(V))$ is clearly equal to $(\pi^l_k)^{-1}(V')$, which by the 
product decomposition~(\ref{pppddd}) is equal to $V'\times\CC^L$. Therefore the subvariety $W'':=\CC^n\times W'$ of the product $\CC^n\times J^l(U,\CC^r)$ is 
equal to the product $\CC^n\times W'\times\CC^L$. Observe that $Y_i\comp\gamma^l_f=y_i\comp\gamma^k_f$ for $1\leq i\leq n+N$. Define 
$H_{i}:=Y_{n+i}\comp\gamma^l_f\comp\sigma_l$ (where $\sigma_l$ is the projection of $\CC^n\times J^l(U,\CC^r)$ to the first factor) for $1\leq i\leq N+L$.
Consider the projection $\beta:=(\mathrm{Id}_{\CC^n},\pi^l_k):\CC^n\times J^l(U,\CC^r)\to\CC^n\times J^k(U,\CC^r)$; observe that $H_i=h_i\comp\beta$ for 
any $i\leq N$. Define $x''':=\gamma^l_f(x)$; observe that $\beta(x''')=x''$.
Let $\mathbf{s}'=(s'_1,....,s'_{n+N+L})$ be the sequence of elements of $\calO_{W'',x'''}$ obtained by projecting the regular sequence 
\[(Y_1-x_1,...,Y_n-x_n,Y_{n+1}-H_1,...,Y_{n+N+L}-H_{N+L}).\]
Using the product structure in $W''$ it is easy to check that $(s'_{n+N+1},...,s'_{n+N+L})$ is a 
regular sequence. Moreover if $Z$ is the analytic subspace of $W''$ defined by the ideal generated by $(s'_{n+N+1},...,s'_{n+N+L})$, and $s''_i$ denotes the 
class of $s'_i$ in $\calO_{Z,x'''}$, then the restriction
$\sigma_{|Z}:Z\to V''$ is an isomorphism satisfying $\sigma_{|Z}^*s_i=s''_i$.   

Let $(s_1,...,s_r)$ be elements of a ring $A$. Let $\overline{s}_2,...,\overline{s}_r$ be the classes of $s_2,...,s_r$ in $A/(s_1)$. We have the decomposition
$\mathrm{K}_\bullet(s_1,...,s_r)=\mathrm{K}_\bullet(s_1)\otimes\mathrm{K}_{\bullet}(s_2,...,s_r)$, in which 
\begin{equation}
\label{ddeecc}
\mathrm{K}_p(s_1,...,s_r)=[\mathrm{K}_{0}(s_1)\otimes\mathrm{K}_p(s_2,...,s_r)]\oplus [\mathrm{K}_{1}(s_0)\otimes\mathrm{K}_{p-1}(s_2,...,s_r)].
\end{equation}
Denote by $\beta_p:\mathrm{K}_p(s_2,...,s_r)\to\mathrm{K}_p(\overline{s}_2,...,\overline{s}_r)$ be the natural epimorphism. 
Consider the morphism of complexes $\alpha_\bullet:\mathrm{K}_\bullet(s_1,...,s_r)\to\mathrm{K}_\bullet(\overline{s}_2,...,\overline{s}_r)$ defined by
$\alpha_p=\beta_p\oplus 0$ in terms of the decomposition~(\ref{ddeecc}).

\begin{lema}
\label{koszul}
If $s_1$ is not a zero divisor of $A$ then $\alpha_\bullet$ is a quasi-isomorphism.
\end{lema}
\begin{proof}
The homomorphism $\alpha$ is clearly surjective in each level. It is easy to check that the complex formed by the kernels is acyclic.
\end{proof}

Using the last lemma repeatedly we obtain 
\begin{equation}
\label{equa1}
\sum_{i=1}^{n+N+L}(-1)^i\dim_{\CC}(\mathrm{H}_i(\mathrm{K}_{\bullet}(\mathbf{s}')))=\sum_{i=1}^{n+N}(-1)^i\dim_{\CC}(\mathrm{H}_i(\mathrm{K}_{\bullet}(s''_1,...,s''_{n+N})));
\end{equation}
due to the fact that $\sigma_{|Z}$ is an isomorphism satisfying $\sigma_{|Z}^*s_i=s''_i$ the last quantity equals 
\[\sum_{i=1}^{n+N}(-1)^i\dim_{\CC}(\mathrm{H}_i(\mathrm{K}_{\bullet}(\mathbf{s}))).\]
This proves the independence on $k$ of Definition~\ref{intmult}.

Let $X$ be a smooth analytic variety and $Y$ an irreducible closed analytic subset of codimension $n$ in $X$. Consider an analytic mapping 
$G:B_\epsilon\times D_\delta\to X$. 
Suppose $G_{|0}^{-1}(Y)=\{O\}$. Choosing $\epsilon$ small enough and $0<\delta<<\epsilon$, the restriction $\pi:Z:=G^{-1}(Y)\to D_\delta$
of the projection to the second factor is a finite map. Consequently $Z$ is a $1$-dimensional closed analytic subspace of $B_\epsilon\times D_\delta$.   
For any $t\in D_\delta$ the set $G_{|t}^{-1}(Y)$ is a finite number of points $\{p_1,...,p_d\}$; therefore any cycle in $Z_0(G_{|t}^{-1}(V))$ 
has a unique expression of the form $\sum_{i=1}^sn_ip_i$, where $n_i$ are integers. Define the {\em degree} of a cycle as 
$\mathrm{deg}(\sum_{i=1}^dn_ip_i):=\sum_{i=1}^dn_i$. Proposition~10.2~of~\cite{Fu} tells precisely that

\begin{lema}
\label{intdef} 
In the preceding situation $\mathrm{deg}(G_{|t}^!([V])$ does not depend on $t$.
\end{lema}

Now we prove the independence of Definition~\ref{intmult} on the functions $f_1,...,f_r$. Suppose that we have another expression $f=\sum_{i=1}^rg_ih_i$.
Choose $\epsilon$ so that $f$, the $f_i$'s and the $h_i$'s are defined in $B_\epsilon$. Defining $F_i:B_\epsilon\times\CC\to\CC$ as $F_i:=(1-t)f_i+th_i$ we 
have $f=\sum_{i=1}^rg_iF_{i|t}$ for each $t$. Define an analytic mapping $\phi_{|t}:B_\epsilon\times\CC\to J^\infty(U,\CC^r)$ by the formula 
$\phi(x,t):=j^\infty F_{1|t}(x),...,j^\infty F_{r|t}(x)$. Choose $\epsilon$ small enough so that $B_\epsilon\cap\rho_f^{-1}(V)=\{O\}$. As 
$\phi_{|t}:B_\epsilon\to J^\infty(U,\CC^r)$ lifts $\rho_f$ we have $\phi_{|t}^{-1}(V)=\{O\}$ for any $t$. Therefore $i_x(\rho_f,V)$ computed in terms of 
$G_{1|t},...,G_{r|t}$ is equal to $\mathrm{deg}(\pi^\infty_k\comp\phi_{|t})^!([\pi^\infty_k(V)]))$. Applying Lemma~\ref{intdef} we obtain that this number is 
independent on $t$. Hence Definition~\ref{intmult} is independent on the choice of the functions $f_1,...,f_r$.

\begin{remark}
If $V$ is $\calD_{\tilde{I},e}$-invariant the independence on choices of Definition~\ref{intmult} also could be proved using the versal unfolding to show that 
for a generic deformation $F_{|t}$ of 
$f$, the number of points of $\rho_{F_{|t}}^{-1}(V)$ is equal to $i_x(\rho_f,V)$. In this way it also follows the independence of 
$i_O(\rho_f,V)$ on the choice of the system of generators of $I$. 
\end{remark}

Another easy consequence of Lemma~\ref{intdef} is the following ``Conservation of Number Formula'':

\begin{prop}
\label{dinamicalint} 
Let $V\subset J^\infty(U,\tilde{I})$ be an irreducible closed analytic subset of codimension $n$. Consider any $f\in I$ for such that 
$O$ is an isolated point of $\rho_f^{-1}(V)$. Let $F:\CC^n\times\CC\to\CC$ be any $I$-unfolding of $f$. For any positive and small enough $\epsilon$ there 
exists a positive $\delta$ such that for any $t\in\CC^k$ with $||t||<\delta$ we have
\begin{equation}
\label{consnum}
i_O(\rho_f,V)=\sum_{x\in B_\epsilon} i_x(\rho_{F_{|t}},V).
\end{equation}
\end{prop}
\begin{proof}
Let $\tilde{\rho}_F:\CC^n\times\CC\to J^\infty(U,\CC^r)$ be an analytic local lifting of $\rho_F$. Observe that 
\[\mathrm{deg}(\pi^\infty_k\comp\tilde{\rho}_{F|t})^!([\pi^\infty_k(V)]))=\sum_{x\in B_\epsilon} i_x(\rho_{F_{|t}},V)\]
and apply Lemma~\ref{intdef}.
\end{proof}

This provides an algebraic formula for the splitting function and the Morse number:

\begin{cor}
Consider $f\in I$ with $c_{I,e}(f)<\infty$. Then
\[\sigma_{I}[f](n)=i_O(\rho_f,Z_n),\]
where $Z_n$ is the union of the irreducible components of $C_1$ such that the extended codimension of a generic member of them is precisely $n$.
\end{cor}

If a germ $f_x\in J^\infty(U,\tilde{I}_x)$ has a Morse point at $x$ then 
\begin{enumerate}
\item the projection $pr(x)$ does not belong to $V(I)$,
\item the function $f_x$ has a critical point in $x$, 
\item the function $f_x$ has positive codimension with respect to $\tilde{I}_x$. 
\end{enumerate}
The germs having Morse points are a dense open subset among the germs having these three properties. 
Clearly conditions $(1)$ and $(2)$ hold simultaneously if and only if conditions $(1)$ and $(3)$ hold at the same time. In other words, if 
$Z$ is the set of germs $f_x$ having $x$ as a critical point then 
\begin{equation}
\label{eeqq}
Z\setminus pr^{-1}(V(I))=C_1\setminus pr^{-1}(V(I)). 
\end{equation}

We define the finite-determined closed analytic subset $M\subset J^\infty(U,\tilde{I})$ as the union of the irreducible components of $C_1$ not contained in 
$pr_\infty^{-1}(V(I))$. We have 

\begin{cor}
Consider $f\in I$ with $c_{I,e}(f)<\infty$. Then 
\[\calM(f)=i_O(\rho_f,M).\]
\end{cor}

Suppose that $M$ is $k$-determined. A key step to effectively compute the Morse number is to have an explicit description of the ideal of  
$M':=(j^k\varphi)^{-1}(\pi^\infty_k(M))$ as a subset of $J^k(U,\CC^r)$. Due to Equality~(\ref{eeqq}) we can obtain such a 
description: 

Observe that $Z$ is clearly $1$-determined. We show that it is a closed 
analytic subset and give generators for the ideal defining $(j^1\varphi)^{-1}(\pi^\infty_1(Z))$. 

Recall that $\CC\{x\}$ denotes the space of convergent power series in $n$ variables $x_1,...,x_n$ and $\mm$ denotes its maximal ideal.
If $\alpha=(\alpha_1,...,\alpha_r)$ we set $x^\alpha:=x_1^{\alpha_1}...x_r^{\alpha_r}$. Given any element 
$h=(h_1,...,h_r)\in\CC\{x\}^r$, for any $\alpha\in\ZZ_{\geq 0}^r$ we denote by $a^i_\alpha(g)$ the coefficient in $x^\alpha$ of the power series expansion of 
$h_i$. The set $U\times (\CC\{x\}/\mm^{k+1})^r$ is an affine space with coordinates $x_1,...,x_n$ and 
$a^i_\alpha$ (for $1\leq i\leq n$ and $\alpha=(\alpha_1,...,\alpha_r)$ such that $|\alpha|:=\sum_{i=1}^r\alpha_i\leq k$). We consider the trivialisation 
$\tau_k:U\times (\CC\{x\}/\mm^{k+1})^r\to J^k(U,\CC^k)$ (see Formula~\ref{r1trivi}). Define $\beta_0,...,\beta_n\in\ZZ_{\geq 0}^r$ by 
$\beta_0=(0,...,0)$ and $\beta_i$ the $n$-tuple whose only non-zero component is the $i$-th one and has value $1$. An easy computation shows that 
$(j^k\varphi\comp\tau_k)^{-1}(\pi^\infty_k(Z))$ the subset of $U\times (\CC\{x\}/\mm^{k+1})^r$ given by the set of common zeros of the functions $Q_1,...,Q_n$,
where
\[Q_i(...,z_i,...,a^i_\alpha,...):=\sum_{j=1}^n[a^j_{\beta_0}\frac{\partial f_j}{\partial x_i}(z_1,...,z_n)+a^j_{\beta_i}f_j(z_1,...,z_n)].\]

Define $J_1:=(Q_1,...,Q_r)$ and let $J_2$ be the 
pullback of $I$ by the projection of $U\times (\CC\{x\}/\mm^{k+1})^r$ to its first factor, that is, the ideal generated by 
$\{g_i(z_1,...,z_n):1\leq i\leq r\}$. Then the ideal of functions vanishing at $M'$ is 
\begin{equation}
\label{formideal}
J_{M'}=(\sqrt{J_1}:J_2).
\end{equation}
This shows in particular that $(j^\infty\varphi)^{-1}(M)$ is $1$-determined: we could have worked with $k=1$.

\begin{remark}
In~\cite{Bo1} it was proved a slightly different Morsification result for hypersurface singularities with critical locus an i.c.i.s. and transversal type 
$A_1$. There the critical locus a generic deformation of a function of finite codimension is the Milnor fibre of the i.c.i.s together with a finite number of 
$A_1$ points. The techniques of this section can be easily adapted to compute the number of $A_1$ points in such a generic deformation.
\end{remark}

\subsection{Conservative numerical invariants} 

The following proposition tells us how to associate a numerical invariant to any closed analytic subset $V$ of pure codimension $n$ of $J^\infty(U,\tilde{I})$:
\begin{prop}
\label{invint}
Suppose that $V$ is $\calD_{\tilde{I},e}$-invariant. Then the intersection multiplicity $i_x(\rho_f,V)$ is a $\calD_{\tilde{I},e}$-invariant defined in
$J^\infty(U,\tilde{I})$.
\end{prop}
\begin{proof}
By simplicity we work at the origin.
Let $f\in I$ and $\phi\in\calD_I$. We have to prove the equality $i_O(\rho_f,V)=i_O(\rho_{\phi_*f},V)$. Consider neighbourhoods $U_1$ and $U_2$ of the origin 
in $\CC^n$ such that $\phi$ is defined in $U_1$ and $\phi(U_1)=U_2$. The mapping $\phi$ induces a bijection 
\[\phi_*:J^k(U_1,\tilde{I})\to J^k(U_2,\tilde{I})\] 
and an analytic isomorphism 
\[\phi_*:J^k(U_1,\CC^r)\to J^k(U_2,\CC^r)\]
for any positive $k$ and $r$. A set of generators $\{f_{1,1},...,f_{1,r}\}$ of 
$\tilde{I}_{|U_1}$ induces a set of generators $\{f_{2,1},...,f_{2,r}\}$ of $\tilde{I}_{|U_2}$ defining $f_{2,i}:=\phi_*f_{1,i}$. Consider the associated 
epimorphisms 
\[\varphi_i:\calO_{U_1}^r\to \tilde{I}_{U_1}\]
for $i=1,2$, defined by $\varphi_i(h_1,...,h_r):=\sum_{j=1}^rh_jf_{i,j}$. They induce respectively mappings 
\[j^k\varphi_i:J^k(U_i,\CC^r)\to J^k(U_i,\tilde{I})\]
satisfying $j^k\varphi_2\comp j^k\phi_*=j^k\phi_*\comp j^k\varphi_1$. 

Pick up $h_1,...,h_r$ satisfying $f=\sum_{i=1}^rh_if_{1,i}$. This induces an analytic lifting $\tilde{\rho}^k_f$ of $\rho^k_f$, defined by the formula 
$\tilde{\rho}^k_f(x):=(j^kh_1(x),...,j^kh_2(x))$. Noticing that $\phi_*f=\sum_{i=1}^r\phi_*h_if_{2,i}$ we deduce that $j^k\phi_*\comp\tilde{\rho}^k_f$ is an 
analytic lifting of $\rho^k_{\phi^*f}$. As
$j^k\phi_*:J^k(U_1,\CC^r)\to J^k(U_2,\CC^r)$ is an analytic isomorphism which satisfies 
\[j^k\phi_*[(j^k\varphi_1)^{-1}(\pi^\infty_k(V))]=(j^k\varphi_2)^{-1}(\pi^\infty_k(V))\]
(because of the $\calD_{I}$-invariance of $V$), the equality $i_O(\rho_f,V)=i_O(\rho_{\phi_*f},V)$ is satisfied.
\end{proof}

Any numerical $\calD_{\tilde{I},e}$-invariant $\Xi:J^\infty(U,\tilde{I})\to\ZZ_{\geq 0}\cup\{\infty\}$ constructed by intersection multiplicity with 
$\calD_{\tilde{I},e}$-invariant closed analytic 
subsets of pure codimension $n$ of $J^\infty(U,\tilde{I})$ clearly satisfy the following properties: consider any $f_x\in\tilde{I}_x$ with 
$c_{\tilde{I}_x,e}(f)<\infty$ and $F:(\CC^n,x)\times(\CC^k,O)\to\CC$, any 
$I$-unfolding of $f$. 
\begin{enumerate}
\item Finiteness: $\Xi(f_x)<\infty$.
\item Analyticity: for any integer $N$ the set 
\[\{(x,t)\in\CC^n\times\CC^k:\Xi(F_{|t,x})>0\}\]
is closed analytic.  
\item Discreteness and conservation of number: there exists a sufficiently small neighbourhood $U_x$ of $x$ such that $\Xi(f_y)=0$ for any 
$y\in U_x\setminus\{x\}$, and, for any positive and small enough $\delta$ we have 
\[\Xi(f_x)=\sum_{y\in U_y}\Xi(F_{|t,y})\]
for any $t\in B_\delta$.
\end{enumerate}

\begin{definition}
A numerical invariant satisfying the above properties will be called a {\em conservative invariant}.
\end{definition}

Observe that the property of conservation of number implies the upper semicontinuity of $\Xi$, that is $\Xi(F_{|t,x})\leq\Xi(f_x)$ for any $t$ close enough to
$x$.

The following theorem tells that any conservative invariant is a (locally) finite sum of invariants constructed by intersection multiplicity:

\begin{theo}
\label{invchar}
Let $\Xi:J^\infty(U,\tilde{I})\to\ZZ_{\geq 0}\cup\{\infty\}$ be a numerical $\calD_{\tilde{I},e}$-invariant satisfying properties (1)-(3) given above. Then 
there exist a family of pairs $(V_i,n_i)$, where $V_i$ is a $n$-codimensional $\calD_{\tilde{I},e}$-invariant closed analytic subset of 
$J^\infty(U,\tilde{I})$, and $n_i$ is a positive integer, such that for any $f_x\in J^\infty(U,\tilde{I})$ satisfying  $c_{\tilde{I}_x,e}(f_x)<\infty$ the 
family $\{V_i\}_{i\in\NN}$ is locally finite at $f_x$ and $\Xi(f_x)=\sum_{i\in\NN}n_i i_x(\rho_f,V_i)$ (the sum is finite by the locally finiteness of the 
family at $f_x$).
\end{theo}
\begin{proof}
Recall that $C_n$ denotes the subset of $J^\infty(U,\tilde{I})$ formed by germs of extended codimension at least $n$. Define the sets 
\[Z:=\{f_x\in J^\infty(U,\tilde{I}):\Xi(f_x)>0\}\quad\quad Z_k:=Z\setminus C_k,\]
for any $k>0$. 

We claim that $Z_k$ is $k+\lambda$-determined, where $\lambda$ is the Uniform Artin-Rees constant. Consider 
$f_x\in Z\setminus C_k$ and $g_x\in\mm_x^{k+\lambda+1}\cap\tilde{I}_x$. We have to show that $f_x+g_x$ belongs to $Z_k$. By Lemma~\ref{deter2.5} the subset 
$C_k$ is $k+\lambda$-determined, and hence $f_x+g_x\notin C_k$.     
As $c_{\tilde{I}_x,e}(f_x)<k$, by Lemma~\ref{deter1} we have $\mm_x^{k-1}\tilde{I}_x\subset\Theta_{\tilde{I}_x,e}(f_x)$. 
As $\mm_x\Theta_{\tilde{I}_x,e}\subset\Theta_{\tilde{I}_x}$ we deduce $\mm_x^{k+1}\tilde{I}_x\subset\mm_x\Theta_{\tilde{I}_x}(f_x)$. 
By Uniform Artin-Rees $\mm_x^{k+\lambda+1}\cap\tilde{I}_x\subset\mm_x\Theta_{\tilde{I}_x}(f_x)$. Then 
the finite $I$-determinacy theorem (see Theorem~6.5 of~\cite{Pe2}) tells that $f_x$ is a $k+\lambda$-determined in $I$; in other words, given any 
$g_x\in\mm_x^{k+\lambda+1}\cap\tilde{I}_x$ then $f_x+g_x$ is in the orbit $\calD_{\tilde{I}_x}(f)$. Consequently, as $Z$ is clearly 
$\calD_{\tilde{I},e}$-invariant and $f_x$ belongs to $Z$, also $f_x+g_x$ belongs to $Z$. This shows the claim.

Our next claim is that the topological closure $\overline{Z}_k$ of $Z_k$ in $J^{\infty}(U,\tilde{I})$ is a $k+2\lambda$-determined closed analytic subset of 
$J^\infty(U,\tilde{I})$, which have no irreducible components contained in $C_k$. Moreover, the set $Z_k$ is closed in $J^\infty(U,\tilde{I})\setminus C_k$.
Let $\{f_1,...,f_N\}$ be monomials in $x_1,...,x_n$ forming a basis of $\calO_{\CC^n,O}/\mm_O^{\lambda+k+1}$. Consider 
$S:=(\calO_{\CC^n,)}/\mm_O^{\lambda+k+1})^r$ viewed as an affine space. Consider a system of coordinates $\{s^i_j\}$ for $1\leq i\leq r$ and $1\leq j\leq N$ 
characterised by the property that the point with coordinates $(a^1_1,...,a^r_N)$ represents the $r$-tuple $(\sum_{j=1}^Na^1_jf_j,..., \sum_{j=1}^Na^r_jf_j)$.
Define the $I$-unfolding $F:U\times S\to\CC$ by the formula
\begin{equation}
\label{elanterior}
F(x,s^1_1,...,s^r_N):=\sum_{i=1}^r\sum_{j=1}^Ns^i_jf_j(x)g_j(x),
\end{equation}
where $g_1,...,g_r$ are the fixed generators for $\tilde{I}$ at $U$. According to the second property of $\Xi$ the subset $A\subset U\times S$ formed by the
pairs $(x,s)$ such that $\Xi(F_{|s,x})>0$ is closed analytic. The mapping $\alpha:U\times S\to J^{k+\lambda}(U,\CC^r)$ defined by 
\[\alpha(x,s^1_1,...,s^r_N)=(\sum_{j=1}^Ns^1_j j^{k+\lambda}f_j(x),...,\sum_{j=1}^Ns^r_j j^{k+\lambda}f_j(x))\]
is an analytic vector bundle isomorphism. We have defined $F$ and $\alpha$ so that the compatibility relation 
$\pi^\infty_{k+\lambda}\comp\rho_{F}=j^{k+\lambda}\varphi\comp\alpha$ holds. This implies
\[\alpha(A)\setminus (j^{k+\lambda}\varphi)^{-1}(\pi^\infty_{k+\lambda}(C_k))=(j^{k+\lambda}\varphi)^{-1}(\pi^\infty_{k+\lambda}(Z_k)).\]
The set $B:=\alpha(A)\cup (j^{k+\lambda}\varphi)^{-1}(\pi^\infty_{k+\lambda}(C_k))$ is a closed analytic subset of 
$J^{k+\lambda}(U,\CC^r)$ which is $j^{k+\lambda}\varphi$-saturated for being equal to the union
\[[\alpha(A)\setminus (j^{k+\lambda}\varphi)^{-1}(\pi^\infty_{k+\lambda}(C_k))]\cup (j^{k+\lambda}\varphi)^{-1}(\pi^\infty_{k+\lambda}(C_k)),\]
whose terms are $j^{k+\lambda}\varphi$-saturated closed analytic subsets. Let $B'$ be the union of the irreducible components of $B$ not contained in 
$(j^{k+\lambda}\varphi)^{-1}(\pi^\infty_{k+\lambda}(C_k))$. Due to Lemma~\ref{irraux} we know that $B'':=(\pi^{k+2\lambda}_{k+\lambda})^{-1}(B')$ is 
$j^{k+2\lambda}\varphi$-saturated, and therefore $Z'_k:=\pi^\infty_{k+2\lambda+2}(j^{k+2\lambda}\varphi(B''))$ is a $k+2\lambda$-determined closed analytic 
subset of $J^\infty(U,\tilde{I})$ with no irreducible components contained in $C_k$. By construction,  it is clear that $Z'_k\setminus C_k=Z_k$. 
By Lemma~\ref{topclo} and the fact that no irreducible component of $Z'_k$ is contained in $C_k$ we conclude that the topological closure of 
$Z_k=Z'_k\setminus C_k$ equals $Z'_k$. Obviously, the set  $\overline{Z}_k\setminus Z_k=Z'_k\setminus (Z'_k\setminus C_k)$ is contained in $C_k$. 
Our claim is proved.

For any positive integers $k$ and $m$ we have 
\begin{equation}
\label{zxcvb}
\overline{Z}_k\setminus C_m\subset Z_m.
\end{equation}
Indeed, consider $f_x\in J^\infty(U,\tilde{I})$ such that $c_{\tilde{I},e}(f)<m$. Suppose that $\Xi(f_x)=0$. Then, as $Z_{m}$ is closed in 
$J^\infty(U,\tilde{I})\setminus C_{m}$, there is an open neighbourhood of $f_x$ in $J^\infty(U,\tilde{I})$ where $\Xi$ vanishes. After this,  
inclusion~(\ref{zxcvb}) in case $k<m$ follows easily. If $k\geq m$ then 
\[\overline{Z}_k\setminus C_m\subset\overline{Z}_k\setminus C_k=Z_k,\]
but as $Z_k\setminus C_m$ is clearly equal to $Z_m$, inclusion~(\ref{zxcvb}) follows.

Consider $Y:=\cup_{k>0}\overline{Z}_k$; inclusion~(\ref{zxcvb}) implies 
\begin{equation}
\label{Ygriega}
Y\setminus C_m=Z_m=\overline{Z}_m\setminus C_m. 
\end{equation}
for any $m>0$. Therefore $Y$ is a finite determined closed analytic subset locally around each point of $J^\infty(U,\tilde{I})\setminus C_\infty$. This is 
expressed in other words as follows: let $\{V_i\}_{i\in\NN}$ be the set whose elements are the irreducible components of any of the $\overline{Z}_k$'s. Then 
the family $\{V_i\}_{i\in\NN}$ is locally finite in $J^\infty(U,\tilde{I})\setminus C_\infty$. Moreover, each $V_i$ is $\calD_{\tilde{I},e}$-invariant, for 
being an irreducible component of a certain $\overline{Z}_k$, which is $\calD_{\tilde{I},e}$-invariant for containing the $\calD_{\tilde{I},e}$-invariant
dense Zariski open subset $Z\setminus C_k$.

Pick up a component $V_i$. Consider $f_x\in V_i$ with $c_{\tilde{I}_x,e}(f_x)<\infty$, which, by the local finiteness
of the family $\{V_i\}_{i\in\NN}$ at $J^\infty(U,\tilde{I})\setminus C_\infty$ can be chosen so that there is a neighbourhood $W$ of $f_x$ in 
$J^\infty(U,\tilde{I})$ satisfying $W\cap Z=V_i$. If $\codim(V_i)<n$, an easy argument using Statement~$(\dagger)$ of the 
proof of Lemma~\ref{codimbound} shows that $\dim_x(\rho_f^{-1}(V_i))>0$. As $\Xi(f_y)>0$ at any $y\in \rho_f^{-1}(V_i)$ we are contradicting the third property
of $\Xi$. Hence $\codim(V_i)\geq n$. Consider the canonical Whitney stratification $\calX$ of $V_i$ (see Theorem~\ref{Whitneystr}); the strata are 
$\calD_{\tilde{I},e}$-invariant. Consider a versal $I$-unfolding $F:U_x\times S\to\CC$ of $f_x$ at $x$ (where $U_x$ is a neighbourhood of $x$), let $s_0$ be the
point in $S$ such that $F_{|s_0,x}=f_x$. 
By Proposition~\ref{versalimplytrans}, the mapping $\rho_F$ is transversal to $\calX$. By Theorem~\ref{transpar} the set of $s\in S$ such that
$\rho_{F|s}$ is transversal to $\calX$ is dense in $S$. Hence, if $s$ is generic and $\codim(V_i)>n$, the set $\rho_{F|s}(U_x)$ cannot meet $V_i$. 
Choosing $S$ and $U_x$ small enough we obtain that the image of $\rho_F$ lies in $W$. Consequently $(\rho_{F|s})^{-1}(Z)=(\rho_{F|s})^{-1}(V_i)=\emptyset$.
This means that $\Xi(F_{|s,y})=0$ for any $y\in U_x$, which, if $s$ is close enough to $s_0$, contradicts the third property of $\Xi$, since $\Xi(f_x)>0$. 
We conclude that $\codim(V_i)=n$. 

Suppose that $V_i$ is a component of $\overline{Z}_k$. As $\overline{Z}_k$ is $k+2\lambda$-saturated, then the component $V_i$ is $k+3\lambda$-saturated.
Let $B$ be the union of all the irreducible components of $\overline{Z}_k$ different from $V_i$; 
then $B\cup C_k$ contains all the components $V_j$ different from $V_i$. The set $\Delta_i:=(B\cup C_k\cup\mathrm{Sing}(V_i))\cap V_i$ is a proper 
$k+4\lambda$-determined closed analytic subset of $V_i$. As $V_i$ is smooth outside $\Delta_i$ we can define
$D_i\subset V_i$ to be the subset of $F_x$ such that $\rho_f$ is not transversal to $V_i$ at $x$. The $k+3\lambda$-determinacy of $V_i$ implies easily that 
$D_i$ is $k+4\lambda$-determined: if $f_x$ and $g_x$ have the same $k+4\lambda$-jet, by Uniform Artin-Rees, they differ in an element of 
$\mm_x^{k+3\lambda+1}\tilde{I}_x$, and we can find liftings $\tilde{\rho}_f$ and $\tilde{\rho}_g$ of $\rho_f$ and $\rho_g$ which are equal up to 
$k+3\lambda$-jets.  

The set $X:=(j^{k+4\lambda+1}\varphi)^{-1}(\pi^\infty_{k+4\lambda+1}(V_i\setminus\Delta_i))$ is a smooth irreducible locally closed analytic subset of 
$J^{k+4\lambda+1}(U,\CC^r)$. Let $\calT\to X$ be its tangent bundle. Denote by $\iota:X\hookrightarrow J^{k+4\lambda+1}(U,\CC^r)$ the inclusion mapping
As $\calT$ is a sub-bundle of the restriction to $X$ of the tangent bundle of 
$J^{k+4\lambda+1}(U,\CC^r)$, the mapping $\pi^{k+4\lambda+1}_{k+2\lambda}$ induces a natural homomorphism 
\[q:\calT\to\calT',\]
where $\calT'$ is the pullback by $\pi^{k+4\lambda+1}_{k+2\lambda}\comp\iota$ of the tangent bundle $\calT''$ of 
$J^{k+2\lambda}(U,\CC^r)$. Observe that, as $pr_{k+2\lambda}:J^{k+2\lambda}(U,\CC^r)\to U$ defines a vector bundle, for any 
$p\in J^{k+2\lambda}(U,\CC^r)$
the fibre $\calT''_p$ decomposes naturally as 
\[T_{pr_{k+2\lambda}(p)}U\oplus J^{k+2\lambda}(U,\CC^r)_p.\]
Any vector of $\calT''_p$ is decomposed accordingly in two components; the first is called the base component, 
and the second is called the fibre component. Consider the trivial vector bundle of rank $n$ over $X$ and choose $e_1,...,e_n$ to be global trivialising 
sections. There is an analytic homomorphism 
\[\sigma:\CC^n\times X\to\calT'\]
obtained defining $\sigma(e_i,j^{k+4\lambda+1}f_1(y),...,j^{k+4\lambda+1}f_r(y))$ to be the unique vector of 
\[\calT'_{(j^{k+4\lambda+1}f_1(y),...,j^{k+4\lambda+1}f_r(y))}=\calT''_{(j^{k+2\lambda}f_1(y),...,j^{k+2\lambda}f_r(y))}\]
whose fibre component equals
$j^{k+2\lambda}(\partial f_1/\partial x_i)(y),...,j^{k+2\lambda}(\partial f_r/\partial x_i)(y))$ and whose base component is $\partial /\partial x_i$.
The set 
\[D':=\{p\in X:q(\calT_p)+\sigma(\CC^n)\neq \calT'_p\}\]
is easily shown to be closed analytic in $X$. As $D'$ has been defined so that 
$(j^{k+4\lambda+1}\varphi)^{-1}(\pi^\infty_{k+4\lambda+1}(D_i))=D'$ we conclude that $D_i$ is a proper closed analytic $k+3\lambda$-determined in 
$V_i\setminus\Delta_i$. 

We claim that $\Xi$ is constant in $V_i\setminus (\Delta_i\cup D_i)$.
The set $\overline{X}$ is an irreducible analytic subset of $J^{k+4\lambda+1}(U,\CC^r)$, for being $V_i$ irreducible. Then, by well known properties 
of complex analytic sets (see~\cite{Lo}, Ch.~IV, $\S$~2), the set $X\setminus D'$ is path-connected. Therefore $V_i\setminus (\Delta_i\cup D_i)$ is 
path-connected. Consequently we only have to show that $\Xi$ is locally constant in $V_i\setminus (\Delta_i\cup D_i)$. Consider
$f_x\in V_i\setminus (\Delta_i\cup D_i)$. Then $f_x\in J^\infty(U,\tilde{I})\setminus C_k$. We have seen that any $g_y\in J^\infty(U,\tilde{I})\setminus C_k$
is $k+\lambda$-determined in $\tilde{I}_y$. Therefore, if we consider the $I$-unfolding $G:U_x\times S\to\CC$ defined by $G(y,s):=f_x(y)+F(x,s)$ (where $U_x$ 
is a neighbourhood of $x$ in $U$ and $F$ is 
the $I$-unfolding defined by formula~(\ref{elanterior})), then any germ $g_y$ belonging to $J^\infty(U_x,\tilde{I})\setminus C_k$ is 
$\calD_{\tilde{I},e}$-equivalent to a germ of the form $G_{|s,y}$ for a certain $s\in S$. Consequently, as $f_x=G_{|0,x}$, to prove our claim is enough to show
that $\Xi(G_{|s,y})$ is constant in a neighbourhood of $(x,0)$ in $U_x\times S$. As $f_x$ does not belong to $B$, taking $U_x$ and a neighbourhood $S'$ of $0$ 
in $S$ small enough we can assume 
\[\rho_{G}^{-1}(Z)\cap (U_x\times S')=\rho_{G}^{-1}(V_i)\cap (U_x\times S'),\]
and $\rho_f^{-1}(Z)=\{x\}$. On the other hand, as $f_x$ does not belong to $D_i$, the mapping $\rho_{G_{|0}}=\rho_{f}$ is transversal to $V_i$ at $x$;  
therefore, by the open nature of transversality, if we choose $U_x$ and $S'$ small enough, the mapping $\rho_{F_{|s}}$ is transversal to $V_i$ for any 
$s\in S'$.
Consequently, for any $s\in S'$ the set $\rho_{G_{|s}}^{-1}(V_i)$ consists of a unique point $y(s)$, which is 
the only point $y\in U_x$ where $\Xi_{G|s,y}>0$. By the third property of $\Xi$ we 
conclude $\Xi(G_{|s,y(s)})=\Xi(f_x)$. This shows the claim.

Let $n_i$ the value that $\Xi$ assumes at any point of $V_i\setminus (\Delta_i\cup D_i)$. It only remains to be shown that 
\begin{equation}
\label{loultimo}
\Xi(f_x)=\sum_{i\in\NN}n_i i_x(\rho_f,V_i)
\end{equation}
 for any $f_x$ of finite extended codimension. Given such $f_x$ defined in $U_x$, consider the associated mapping 
$\rho_f:U_x\to J^\infty(U,\tilde{I})$. The local finiteness of the family $\{V_i\}_{i\in\NN}$ implies that $U_x$ can be taken small enough so that 
$\rho_f^{-1}(Z)=\{x\}$ and the image of $\rho_f$ meets only finitely many components $V_{i_1},...,V_{i_m}$. We consider a versal 
$I$-unfolding $F:U_x\times S\to\CC$ of $f_x$. An (analogous to previous ones) transversality reasoning using Proposition~\ref{versalimplytrans} and 
Theorem~\ref{transpar} shows that the subset of $s\in S$ such that $\rho_{F_{|s}}$ only meets $Z$ in $V_{i_j}\setminus (\Delta_{i_j}\cup D_{i_j})$ 
(for $1\leq j\leq m$), and it does it
transversally, is dense in $S$. Using this, the third property of $\Xi$, and Proposition~\ref{dinamicalint}, the proof of equation~(\ref{loultimo}) is 
straightforward.
\end{proof}

\section{Examples}

In this section we illustrate our the Morsification Theory spelling it out for certain classes of ideals. 

\subsection{Classical case: $I=\calO_{\CC^n,O}$}
This is the case of isolated singularities. Here our theory recovers the classical Morsification Theorem. Any conservative invariant is 
a multiple of the Milnor number, for being the set of singular germs $M$ (which is the closure of the set of Morse germs) the only $n$-codimensional closed 
analytic subset of $J^\infty(U,\calO_{\CC^n})$.

\subsection{The analytic subspace $V(I)$ is smooth}
Let $k$ be the codimension of $V(I)$. We can give coordinate functions $y_1,...,y_k,x_1,...,x_p$ (with $k+p=n$) of $\CC^n$ such that $I=(y_1,...,y_k)$. 
We set $V:=V(I)$. Let $\tilde{I}$ be the ideal sheaf generated by $y_1,...,y_k$. We study the $\calD_{\tilde{I},e}$-invariant analytic subspaces of 
$J^\infty(\CC^n,\tilde{I})$ which are of codimension at most $n$; let $Z$ be such a subspace. Consider the projection mapping 
$pr_\infty:J^\infty(\CC^n,\tilde{I})\to\CC^n$.

If $pr_\infty(Z)$ is not contained in $V$ then $Z_{|\CC^n\setminus V}$ is a non-empty $\calD_{\tilde{I},e}$-invariant subset of 
$J^\infty(\CC^n\setminus V,\tilde{I})$. As $\tilde{I}_{|\CC^n\setminus V}=\calO_{\CC^n\setminus V}$ then, by the case $I=\calO_{\CC^n,O}$, the subspace 
$Z_{|\CC^n\setminus V}$ is either the total space $J^\infty(\CC^n\setminus V,\tilde{I})$ or the set of singular germs. Therefore $Z$ is either equal to 
$J^\infty(\CC^n,\tilde{I})$ or to the $n$-codimensional subvariety $M$ which is the closure of the set of Morse points. 

Suppose that $pr(Z)\subset V$. Let $A$ be the set of integer multi-indexes $\alpha\in\ZZ_{\geq 0}^k$ such that $|\alpha|\geq 2$. Any function $f\in I$ can be 
written as 
\begin{equation}
\label{representacion}
f(x_1,...,x_p,y_1,...,y_k)=\sum_{i=1}^kf_i(x_1,...,x_p)y_i+\sum_{\alpha\in A}a_\alpha(x_1,...,x_n)y^{\alpha},
\end{equation}
with the $f_i$'s and $a_\alpha$'s convergent power series in $x_1,...,x_p$. 
It is easy to check that if $f$ is non-singular (that is $f_i(0,...,0)\neq 0$ for a certain $i\leq k$) then $f$ is $\calD_I$-equivalent to $f=y_1$. 

Let $Z_1$ be the subset of $J^\infty(\CC^n,\tilde{I})_{|V}$
formed by germs $f_x$ having $x$ as a critical point. This condition, in terms of the expression~(\ref{representacion}), means precisely 
$f_i(0,...,0)=0$ for any $i\leq k$. Therefore $Z_1$ is a $k$-codimensional $1$-determined closed analytic subset of $J^\infty(\CC^n,\tilde{I})_{|V}$. As the 
codimension of $J^\infty(\CC^n,\tilde{I})_{|V}$ in $J^\infty(\CC^n,\tilde{I})$ is $k$, the set $Z_1$ is $2k$-codimensional in $J^\infty(\CC^n,\tilde{I})$.
Consequently if $p<k$ then all the $\calD_{\tilde{I},e}$-invariant proper analytic subsets of $J^\infty(\CC^n,\tilde{I})$ which are of codimension at most $n$ 
are $M$ and $J^\infty(\CC^n,\tilde{I})_{|V}$.

Assume $p\geq k$. Let $f\in Z_1$, express it as in~(\ref{representacion}). If the differentials 
\[\{df_i(0,...,0):1\leq i\leq k\}\]
are linearly independent then $f$ is easily seen to be $\calD_I$-equivalent to 
\begin{equation}
\label{tilded}
\sum_{i=1}^kx_iy_i.
\end{equation}
We say that a function of $I$ is of type $\tilde{D}(n-2k,0)$ if it is $\calD_I$-equivalent to~(\ref{tilded}). Functions of type $\tilde{D}(n-2k,0)$
have a $2k$-codimensional critical locus and transversal type $A_1$, they are $\calD$-equivalent to a function of type $D(n-2k,0)$ 
(see Definition~\ref{defdkp}).

The subset $Z_2\subset Z_1$  
consisting of germs for which $\{df_i(0,...,0)\}_{1\leq i\leq k}$ are not linearly independent is a $2$-determined analytic subset of codimension
$p-k+1$ in $Z_1$, and hence of codimension $2k+p-k+1=n+1$ in $J^\infty(\CC^n,\tilde{I})$. Therefore if $p\geq k$ all the $\calD_{\tilde{I},e}$-invariant 
analytic subsets of $J^\infty(\CC^n,\tilde{I})$ which are of codimension at most $n$ are $M$, $J^\infty(\CC^n,\tilde{I})_{|V}$ and $Z_1$, and their 
codimensions are $n$, $k$ and $2k$.

Our Morsification theory tells us that any function $f$ with $c_{I,e}<\infty$ can be approximated (preserving the geometry at the boundary of a Milnor ball) 
by a function whose singularities are finitely many Morse points and a smooth $2k$-codimensional set of points of type $\tilde{D}(n-2,0)$.
Moreover any conservative invariant $\Xi$ is of the form $\Xi=n_1\calM$ if $n\neq 2k$ and $\Xi=n_1\calM +n_2\Delta$ if $n=2k$, where 
$\calM$ is the Morse number, $\Delta$ is the number of $\tilde{D}(n-2,0)$ appearing in a generic deformation of $f$.

\begin{remark}
Let $I$ be arbitrary and consider $f\in I$ with $c_{I,e}(f)<0$. The above reasoning tells that the only singularities that a generic $I$-deformation $g$ close 
enough to $f$ can have at a smooth point of $V(I)$ are of type $\tilde{D}(n-2k,0)$, where $k=\codim_x(V(I),\CC^n)$. Moreover if such a singularity appear then 
there is a $2k$-codimensional locally closed subset of $\CC^n$ where $g$ has this singularity type.
\end{remark}
 
\subsection{Transversal type $A_1$}

In \cite{Bo1}, \cite{Pe3}, \cite{Ne}, \cite{Si1}, \cite{Si2}, \cite{Za}, in order to study functions which are singular with transversal type $A_1$ and a 
certain fixed singular locus, the following point of view was taken: consider the ideal $J$ defining the singular locus (in all these works the singular locus 
is asked to be either an i.c.i.s. or a low dimensional isolated singularity at the origin) with its reduced structure and study the singularities 
appearing in generic deformations of functions which are finite codimensional with respect to the primitive ideal $\int J$. 

Suppose that $J$ is a radical ideal. Consider the variety $V=V(J)$ and the stratification $\CC^n=W\coprod V_0\coprod V_1$, where $U=\CC^n\setminus V$, 
$V_0=V\setminus\mathrm{Sing}(V)$ and $V_1=\mathrm{Sing}(V)$. Let $\tilde{J}$ be the ideal sheaf of functions vanishing at $V$. Denote the primitive ideal 
$\int\! J$ by $I$ and let $\tilde{I}$ be the ideal sheaf associated to it. 

We study the $\calD_{\tilde{I},e}$-invariant subspaces of $J^\infty(U,\tilde{I})$ of codimension at most $n$ and whose image by the projection 
$pr:J^\infty(U,\tilde{I})\to U$ is not contained in $V_1$. Let $Z_1$ be the closed $1$-determined subset formed by germs $f_x\in J^\infty(U,\tilde{I})$ having
a critical point at $x$. Clearly $Z_1=M\cup \pi^{-1}(V)$, where $M$ is the closure of the set of germs $f_x\in J^\infty(U,\tilde{I})$ having an 
isolated singularity at $x$; the set of germs $f_x$ which have an $A_1$-critical point at $x$ are a dense open subset in $M$.

Given any $f_x\in pr^{-1}(V)$ we define $rk(f_x)$ to be the rank of its Hessian matrix at $x$. The set $K_r:=\{f_x\in pr^{-1}(V):rk(f_x)\leq r\}$ is a 
$2$-determined $\calD_{\tilde{I},e}$-invariant closed subset for any integer $r$. Let $C_{1|V}$ be the closed analytic subset consisting of germs of 
$pr^{-1}(V)$ of extended codimension at least $1$.

Consider any $x\in V_0$; if $d=\dim_x(V)$ there is a coordinate system 
\[(x_1,...,x_d,y_1,...,y_m)\]
of $\CC^n$ at $x$ such that $\tilde{J}_x=(y_1,...,y_m)$. In Pellikaan~\cite{Pe2} it is proved that $\tilde{I}_x=\int\!\tilde{J}_x=\tilde{J}_x^2$. Therefore, 
we can express any $g\in\tilde{I}_x$ as $g=\sum_{i,j\leq m}h_{i,j}y_iy_j$, where $h_{i,j}\in\calO_{\CC^n,x}$ and $h_{i,j}=h_{j,i}$; moreover the $h_{i,j}$ are 
unique modulo the ideal $\tilde{J}_x$, hence the matrix $(h'_{i,j})_{i,j\leq m}$ (where 
$h'_{i,j}:=h_{i,j}(O)+\sum_{k=1}^d\frac{\partial h_{i,j}}{\partial x_k}(O)x_k$) is well defined.

\begin{definition}[Pellikaan,~\cite{Pe2}]
\label{defdkp}
Let $f\in\calO_{\CC^n,x}$. We say that $f$ is of type $D(d,k)$ if there is a coordinate system 
$\{x_1,...,x_d,y_1,...,y_m\}$ of $\CC^n$ at $x$ such that 
\[f(x,y)=\sum_{i,j\leq k}l_{i,j}y_iy_j+\sum_{i=k+1}^my_i^2,\] 
where $\{l_{i,j}:i,j\leq k\}$ is a collection of linearly independent linear forms in $x_1,...,x_d$.
\end{definition}

Let $f_x\in\tilde{I}_x$. In~\cite{Pe2} it is proved that $c_{\tilde{I},e}(f_x)=0$ if and only if $f$ is of type $D(d,k)$ for a certain $k$. Consequently, if 
$f_x\in K_r$ then $c_{\tilde{I},e}(f_x)=0$ if and only if $f$ is of type $D(d,d-r)$. 

In~\cite{Bo1} it was shown that a generic 
deformation within $\tilde{I}_x$ of any $f_x\in\tilde{I}_x$ with $c_{\tilde{I}_x,e}<\infty$ only has $A_1$ and $D(d,k)$ points as critical points 
(for $k\leq r)$; moreover the locus where the deformation have $D(d,k)$ points is a smooth subvariety of codimension $k(k+1)/2$ in $V_x$; in~\cite{Bo1} the 
singular locus $V_x$ of $f_x$ is smoothed while deforming $f_x$, but in our case it is already smooth (for being $x\in V_0$).

Let $C'$ be an irreducible component of $C_{1|V}$ such that its image under $pr$ is not contained in $V_1$. We claim that its codimension is strictly 
bigger than $n$: suppose that the codimension of $C'$ is smaller or equal than $n$, consider $f_x\in C'$; if the codimension $c$ of $C'$ equals $n$, by the 
Conservation of Number Formula for intersection multiplicities, any generic deformation of $f_x$ within $\tilde{I}_x$ must contain points of positive extended 
codimension within $V_0$. This is a contradiction with the fact that in a generic deformation only points of type $D(d,k)$ arise in $V_0$. If the codimension 
of $C'$ is strictly smaller than $n$ we consider a subvariety $C''$ of $C'$ of codimension $n$ containing $f_x$ and repeat the argument. 

Define $C$ to be the 
union of the components of $C_{1|V}$ with image under $pr$ not contained in $V_1$. Let $C'$ be the union of the other components. 

Let $V'$ be a connected component of $V_0$ of dimension $d$. Consider $x\in V'$ and $f_x\in\tilde{I}_x$ with $rk(f_x)=r$. Either $f_x$ is of type
$D(d,d-r)$ or $f_x$ belongs to $C$. Hence any irreducible component $Y$ of $K_{r|V}$ not contained in $C$ can be expressed as $Y_0\coprod (Y\cap C)$, where 
$Y_0$ is the open subset containing points of type $D(d,d-r)$. As the locus of $D(d,d-r)$-points of a generic function $g_x$ within $\tilde{I}_x$ has 
codimension $(d-r)(d-r+1)/2$ in $V_0$ then the codimension of $Y$ in $J^\infty(U,\tilde{I})$ is $n-d+(d-r)(d-r+1)/2$.

We have the following splittings in locally closed $\calD_{\tilde{I},e}$-invariant subsets of $J^\infty(U,\tilde{I})$: 
\[J^\infty(U,\tilde{I})=Z_1\coprod (J^\infty(U,\tilde{I})\setminus Z_1)\]
\[Z_1=M\setminus pr^{-1}(V)\coprod pr^{-1}(V)\]
\[pr^{-1}(V)=pr^{-1}(V_1)\coprod (C\setminus pr^{-1}(V_1))\coprod (\coprod_{r\in\ZZ_{\geq 0}}[K_r\setminus (C\cup K_{r-1}\cup pr^{-1}(V_1))]).\]

Applying the Relative Morsification Theorem to them we obtain 

\begin{prop}
\label{tipomorse}
Let $I\subset\calO_{\CC^n,O}$ be a radical ideal and $V$ be the subvariety defined by it. Let $V_1,...,V_r$ be the connected components of 
$V\setminus\mathrm{Sing}(V)$; let $d_i$ be the dimension of $V_i$. If $f\in\int I$ is such that $c_{\int I,e}<\infty$, then given a small
neighbourhood $U$ of the origin, any generic deformation of $f$ sufficiently close to $f$ 
\begin{itemize}
\item has only $A_1$ singularities in  $U\setminus V$.
\item only has singularities of type $D(d_i,k)$ at $V_i$, for $k\leq d_i$. Moreover the locus of points of type $D(d_i,k)$ is a smooth subvariety of 
codimension $k(k+1)/2$ in $V_i$. 
\end{itemize}
Suppose that in addition $I$ has an isolated singularity at the origin. Then the topological type at the origin of any generic deformation of $f$ is the 
generic topological type of a function in $\int I$ (such a generic topological type exists by the results of~\cite{Bo2}). Moreover any conservative invariant 
is the sum an integer multiple of the Morse number, and integer multiples of the number of $D(d_i,k)$ points in a generic deformation, 
for $i\leq r$ and $k$ such that $k(k+1)/2=d_i$. 
\end{prop}

\subsection{Line singularities with simple transversal type}

In~\cite{Jo} line singularities with transversal types $A_1$, $A_2$, $A_3$, $D_4$, $E_6$, $E_7$ and $E_8$ where studied from a topological point of view 
using a Morsification result. For any transversal type as above it is constructed an ideal $I(S)$ such that any singularity with transversal type 
$S$ has a right representative in $I(S)$, and if two functions $f,g\in I(S)$ are $R$-equivalent, then they are $\calD_{I(S)}$-equivalent. For any $S$ 
the orbits of $\calD_{I(S)}$ of codimension $1$ in $J^\infty(V(I),\tilde{I(S)})$ are determined, and the singularity types determined by them are called 
$F_iS$, for $1\leq i\leq 3$. 
It is proved that any function $f\in I(S)$ can be deformed within $I(S)$ to any function having only $A_1$-points outside the singular line $L$, a finite 
number $h_i(f)$ of points of type $F_iS$ in $L$, for $1\leq i\leq 3$, and the generic singularity in $I(S)$ along the rest of the points of $L$.

Our theory recovers the morsification result, interprets $h_i(f)$ as intersection multiplicities, and shows that any conservative invariant 
$\Xi$ is of the form 
\[\Xi(f)=\sum_{i=1}^3n_ih_i(f)+n_4\calM(f).\]

\subsection{Other examples}

\begin{example}{\em 
Take $I=\mm^k$; then $\Theta_{I,e}=\mm\Theta$. Give weight $1$ to the variables $x_i$'s and $-1$ to the derivations $\partial/\partial x_i$'s.
Then the $0$-weight graded piece of $\Theta_{I,e}$ is generated by $x_i\partial/\partial x_j$, with $1\leq i,j\leq n$. Therefore for any function $f$, the 
module $I/\tau_{I,e}(f)+\mm I$ has complex dimension at least $N_k-2n$, where $N_k$ is the dimension of the space of homogeneous polynomials in $n$ variables 
of degree $k$. This provides examples of ideals $I$ for which the support of $\sigma_I[f]$ contains arbitrarily high integers for any $f\in I$.} 
\end{example}

\begin{example}{\em
Choose a function $g\in\CC\{x_1,...,x_{n-1}\}$ with an isolated singularity at the origin. View $g$ as an element of $\calO_{\CC^n,O}$; its zero-set
is singular along the line $L$ defined by $x_1=...=x_{n-1}=0$. Define $I=(g^2)$. Any $X\in\Theta_{I,e}$ admits a unique decomposition as $X=X_1+X_2$, where
\[X_1=a_n\partial/\partial x_n\quad\quad X_2=\sum_{i=1}^{n-1}a_i\partial/\partial x_i,\]
for $a_1,...,a_n\in\calO_{\CC^n,O}$. There is a unique expression  $X_2=\sum_{k=0}^\infty x_n^kX_{2,k}$, where each $X_{2,k}$ belongs to
$\CC\{x_1,...,x_{n-1}\}(\partial/\partial X_1,...,\partial/\partial x_{n-1})$. Then, each $X_{2,k}$ belongs to $\Theta_{I',e}$ where $I'$ is the 
ideal generated by $g$ in $\CC\{x_1,...,x_{n-1}\}$. Define $Y_2:=X_{2,0}$ and $Y_3=\sum_{k=1}^\infty x_n^{k-1}\times X_{2,k}$. We have the decomposition
\begin{equation}
\label{dection}
X=X_1+Y_2+x_nY_3.
\end{equation}

Consider $f\in I$ of the form $f=pg^2$, where $p$ is a polynomial in $x_n$ of degree $k\geq 1$. Decomposition~(\ref{dection}) implies 
\begin{equation}
\label{ecu1}
(\frac{dp}{dx_n}g^2)+\tau_{I',e}(g^2)\calO_{\CC^n,O}\subset\tau_{I,e}(f)
\end{equation}
As $g^2J(g)\subset\tau_{I',e}(g^2)$ (where $J(g)$ is the Jacobian ideal of $g$) and $(x_n^{k-1}g^2)\subset (p'(x_n)g^2)$ we have 
\begin{equation}
\label{cotasup}
c_{I,e}(f)\leq (k-1)\mu(g),
\end{equation}
where $\mu(g)$ is the Milnor number of $g$. 

If $dp/dx_n$ vanishes at $0$ decomposition~(\ref{dection}) implies 
\begin{equation}
\label{ecu2}
\tau_{I,e}(f)\subset (x_ng^2)+\tau_{I',e}(g^2)\calO_{\CC^n,O}.
\end{equation}
In this case we have  
\begin{equation}
\label{cotainf}
c_{I,e}(f)\geq\dim_\CC(\frac{(g^2)}{\tau_{I',e}(g^2)})=\dim_\CC(\frac{(g)}{\tau_{I',e}(g)}). 
\end{equation}

We choose $n=3$, $g(x_1,x_2):=(x_2^2+x_1^3)^2+x_2^5$ and $p(x_3):=1+x_3^k$ with $k>2$. We weight the variables by $wt(x_1)=2$, $wt(x_2)=3$, $wt(x_3)=0$.
A computation shows that for any $X\in\Theta_{I',e}$ we have $X\in (x_2\partial/\partial x_1)+ (x_1,x_2)^2(\partial/\partial x_1,\partial/\partial x_2)$ and
$wt(X(g))\geq 15$. As $X(g)=hg$ for a certain $h\in\CC\{x_1,x_2\}$ we have 
that $wt(h)\geq 3$. Summarising
\begin{equation}
\label{estimacion}
\Theta_{I',e}\subset (x_2\partial/\partial x_1)+ (x_1,x_2)^2(\partial/\partial x_1,\partial/\partial x_2)\quad\quad
\tau_{I',e}(g)\subset (x_1^2g,x_2g).
\end{equation}

Consider $h_1,...,h_m$ in $\CC\{x_1,x_2\}$ forming a basis of the complex vector space $\CC\{x_1,x_2\}/J(g)$. By~(\ref{ecu1}), the $I$-unfolding 
\[F=f(x_1,x_2,x_3)+\sum_{i=1}^r\sum_{j=0}^{k-1}t_{i,j}h_i(x_1,x_2)x_3^jg^2(x_1,x_2)\]
depending on parameters $t_{i,j}$ is versal. We can choose $h_1=1$ and $h_i\in (x_1,x_2)$ if $i>1$. Then for any value $t$ in the parameter space 
$F_{|t}$ is of the form 
\[F_{|t}=[1+h_{t}(x_3)+q(x_1,x_2,x_3)]g^2,\]
with $q(x_1,x_2,x_3)\in (x_1,x_2)$ and where $h_{t}(x_n)$ is a polynomial in $x_3$ of degree $k$. Choose $t$ generic; we study the points $x$ in which 
$F_{|t,x}$ has positive extended codimension. Let $L$ be the line defined by $x_1=x_2=0$. The set $V(g)$ is smooth outside $L$, by Proposition~\ref{tipomorse} 
the only points in which $F_{|t,x}$ cam have extended codimension outside $L$ are $A_1$-points. Using that $\partial/\partial x_3\in\Theta_{I,e}$ we deduce 
that when $p'(x_n)$ is not zero then $(0,0,x_n)$ is a point of extended codimension $0$. We can assume that the derivative $p'(x_n)$ has $k-1$ simple roots. 
Let $a$ be a root; consider $x=(0,0,a)$. A computation taking into account inclusions~(\ref{estimacion}) yield 
\[\tau_{\tilde{I}_x,e}(F_{|t,x})\subset (x_1^2g^2,x_2g^2,(x_3+\mu x_1)g^2),\]
where $\mu$ is the coefficient of $x_1x_3$ in $q$. Therefore $\tau_{\tilde{I}_x,e}(F_{|t,x})\geq 2$. 

This example shows how a point of positive extended codimension splits in a generic deformation in several points of positive extended codimension, such 
that several (k-1) of them have extended codimension at least $2$.}
\end{example}

\end{document}